\documentclass{amsart}
\usepackage{amssymb}

\usepackage{graphicx}
\usepackage{psfig}
\usepackage{amscd}

\theoremstyle{plain}
\newtheorem{theorem}{Theorem}[section]
\newtheorem{corollary}[theorem]{Corollary}
\newtheorem{lemma}[theorem]{Lemma}
\newtheorem{proposition}[theorem]{Proposition}
\newtheorem{summary}{\mdseries \scshape Summary}

\theoremstyle{definition}

\newtheorem{scholium}[theorem]{Scholium}
\newtheorem{example}[theorem]{Example}

\theoremstyle{remark}
\newtheorem{remark}[theorem]{Remark\upshape }

\newtheorem{acknowledgements}{Acknowledgements}

\numberwithin{equation}{section}

\newcommand{\at}{{@}}

\newcommand{\frogkern}{\kern0.5em}

\newlength{\customskipamount}
\newlength{\qedskip}
\newlength{\qedadjust}
\newbox\ipbox
\newcommand{\ip}[2]{\left\langle #1\mathrel{\mathchoice
{\setbox\ipbox=\hbox{$\displaystyle \mathstrut #1#2$}
\vrule height\ht\ipbox width0.25pt depth\dp\ipbox}
{\setbox\ipbox=\hbox{$\textstyle \mathstrut #1#2$}
\vrule height\ht\ipbox width0.25pt depth\dp\ipbox}
{\setbox\ipbox=\hbox{$\scriptstyle \mathstrut #1#2$}
\vrule height\ht\ipbox width0.25pt depth\dp\ipbox}
{\setbox\ipbox=\hbox{$\scriptscriptstyle \mathstrut #1#2$}
\vrule height\ht\ipbox width0.25pt depth\dp\ipbox}
} #2\right\rangle}

\newbox\crossbox
\def\bigcross{\mathop{\mathchoice
{\setbox\crossbox=\hbox{\Huge \bfseries $\times$}
\raise\fontdimen22\textfont2%
\hbox{\raise0.5\dp\crossbox%
\hbox{\lower0.5\ht\crossbox%
\box\crossbox}}}
{\setbox\crossbox=\hbox{\huge \bfseries $\times$}
\raise\fontdimen22\textfont2%
\hbox{\raise0.5\dp\crossbox%
\hbox{\lower0.5\ht\crossbox%
\box\crossbox}}}
{\setbox\crossbox=\hbox{\LARGE \bfseries $\times$}
\raise\fontdimen22\scriptfont2%
\hbox{\raise0.5\dp\crossbox%
\hbox{\lower0.5\ht\crossbox%
\box\crossbox}}}
{\setbox\crossbox=\hbox{\Large \bfseries $\times$}
\raise\fontdimen22\scriptscriptfont2%
\hbox{\raise0.5\dp\crossbox%
\hbox{\lower0.5\ht\crossbox%
\box\crossbox}}}}}

\def\func#1{\mathop{\rm #1}}%
\def\limfunc#1{\mathop{\rm #1}}%
\long\def\TeXButton#1#2{#2}%

\begin{document}
\title[Dense analytic subspaces in fractal $L^{2}$-spaces]{Dense
analytic subspaces in fractal $L^{2}$-spaces}
\author{Palle E. T. Jorgensen}
\address{Department of Mathematics\\
The University of Iowa\\
Iowa City, IA 52242\\
U.S.A.}
\email{jorgen%
\at%
%
math.uiowa.edu}
\author{Steen Pedersen}
\address{Department of Mathematics\\
Wright State University\\
Dayton, OH 45435\\
U.S.A.}
\email{steen%
\at%
%
math.wright.edu}
\thanks{Work supported by the National Science Foundation}

\begin{abstract}
We consider self-similar measures $\mu $ with support in the interval $0\leq
x\leq 1$ which have the analytic functions $\left\{ e^{i2\pi nx}:n=0,1,2,%
\dots%
%
\right\} $ span a dense subspace in $L^{2}\left( \mu \right) $. Depending on
the fractal dimension of $\mu $, we identify subsets $P\subset \mathbb{N}%
_{0}=\left\{ 0,1,2,%
\dots%
%
\right\} $ such that the functions $\left\{ e_{n}:n\in P\right\} $ form an
orthonormal basis for $L^{2}\left( \mu \right) $. We also give a
higher-dimensional affine construction leading to self-similar measures $\mu 
$ with support in $\mathbb{R}^{\nu }$. It is obtained from a given expansive 
$\nu $-by-$\nu $ matrix and a finite set of translation vectors, and we show
that the corresponding $L^{2}\left( \mu \right) $ has an orthonormal basis
of exponentials $e^{i2\pi \lambda \cdot x}$, indexed by vectors $\lambda $
in $\mathbb{R}^{\nu }$, provided certain geometric conditions 
(involving the Ruelle transfer operator)
hold for the
affine system.
\end{abstract}

\maketitle
\tableofcontents
\listoffigures

\section{\label{S1}Introduction: Fractal Measures}

We use the notation $e_{n}\left( x\right) :=e^{i2\pi nx}$, or in complex
form $e_{n}=z^{n}$ where $z=e^{i2\pi x}$, $x\in \mathbb{R}$. Recall that if $%
\mu $ is Lebesgue measure on $I=\left[ 0,1\right] \simeq \mathbb{R}\diagup 
\mathbb{Z}=\mathbb{T}$, then the functions $\left\{ e_{n}:n\in \mathbb{N}%
_{0}\right\} $ span the Hardy space $H_{2}$ of analytic functions on $%
\mathbb{T}$, and $\left\{ e_{n}:n\in \mathbb{Z}\right\} $ is an orthonormal
basis for $L^{2}\left( I\right) $ with normalized Lebesgue measure$.$
Functions of the form $F\left( z\right) =\sum_{n=0}^{\infty }A_{n}z^{n}$
with $\sum_{n=0}^{\infty }\left| A_{n}\right| ^{2}<\infty $ may be viewed
also as defined on $\,\mathcal{U}:=\left\{ z\in \mathbb{C}:\left| z\right|
\leq 1\right\} $ and $F\left( z\right) $ may be identified with its values
on the boundary $\partial \,\mathcal{U}=\mathbb{T}$, represented by $%
z=e^{i2\pi x}$, $x\in \mathbb{R}$. Writing $f\left( x\right) :=F\left(
e^{i2\pi x}\right) $, and defining $E_{z}\left( x\right) =\sum_{n=0}^{\infty
}\bar{z}^{n}e^{i2\pi nx}$, we get the familiar reproducing kernel formula 
\begin{equation}
F\left( z\right) =%
\ip{E_{z}}{f}%
%
=\int_{0}^{1}\overline{E_{z}\left( x\right) }f\left( x\right)
\,dx=\sum_{n=0}^{\infty }z^{n}A_{n},\quad \left| z\right| <1.  \label{eq1}
\end{equation}
Since we show that, for some fractal measures $\mu $ on $I=\left[ 0,1\right] 
$, there are associated subsets $P\subset \mathbb{N}_{0}$ with orthogonal
and total $\left\{ e_{n}:n\in P\right\} $ in $L^{2}\left( \mu \right) $, it
follows that this property then carries over to $L^{2}\left( \mu \right) $,
i.e., the corresponding Hardy space $H_{2}\left( P,\mu \right) $ is all of $%
L^{2}\left( \mu \right) $.

In the present paper, we will consider fractal measures $\mu $ supported on
compact sets (typically totally disconnected) in $\mathbb{R}^{\nu }$ which
arise from iteration algorithms that are defined from a given affine system
of maps in $\mathbb{R}^{\nu }$. These maps in turn will be defined from an
expansive real $\nu $-by-$\nu $ matrix $R$, and a given ``sparse'' set of
translation vectors in $\mathbb{R}^{\nu }$. The fractal nature of these
constructions arises from choosing the number of translation vectors to be
strictly smaller than $\lvert \det R\rvert$ where $R$ is the given expansive
matrix. Our aim will be to understand the possible \emph{orthogonal}
harmonic function systems in $L^{2}\left( \mu \right) $, and, in particular,
to give conditions for $L^{2}\left( \mu \right) $ to have an orthonormal
harmonic basis. (This continues earlier work of the co-authors \cite{JoPe97}%
.) Many of the ideas going into the general case are present already in a
relatively simple example in a single dimension, and this will be introduced
and discussed in Section \ref{S2} below.

While there has been a great amount of recent interest in measures which
arise from geometric iteration algorithms, see, e.g., 
\cite{Ban91}, \cite{Ban96}, \cite{BaGe94},
\cite{BoTa87}, \cite
{Hof95}, \cite{LaWa96c}, \cite{LaWa96d}, \cite{KSS95}, \cite{Rue88}, and 
\cite{PoSi95}, there are relatively few results on the corresponding
harmonic analysis of these measures. But recently (see, e.g., \cite{BoTa87}
and \cite{Hof95}) diffraction lines in quasicrystals have lent themselves to
such a harmonic analysis. In our earlier paper \cite{JoPe97}, we noted a
connection between certain diffraction spectra and the spectral duality
considerations which go into the harmonic analysis of affine self-similarity
(for the measures $\mu $ which we sketched above).

While previous harmonic-analysis approaches to self-affine measures (see,
e.g., \cite{Str89}) have stressed continuous transforms and asymptotic
summation methods, our present approach is discrete and directly based on
Fourier series, for those $\mu $ which have the orthogonality property under
discussion.

A central tool in our work
is a certain \emph{double duality:}
first the usual duality of
Fourier analysis, corresponding
to the dual variables on
either side of the spectral
transform; and secondly
a duality which derives
from our use of \emph{matrix
scaling}. Small scales
correspond to \emph{compact}
attractors of fractal Hausdorff
dimension, while large
scales (``fractals in the large'')
correspond to a \emph{discrete}
set of frequencies (in $\nu $
dimensions),
$\lambda =\left( \lambda _{1},\dots ,\lambda _{\nu }\right) 
\in \mathbb{R}^{\nu }$
which label our Fourier basis of orthogonal exponentials
$e_{\lambda }\left( x\right) :=e^{i2\pi \lambda \cdot x}$
where $x$ is restricted to the (``small scale'') fractal.
In our setup, both scales, small and large, are
\emph{finitely generated}, referring to two given finite subsets
$B$ and $L$ in $\mathbb{R}^{\nu }$
(one on each side of the duality) which
are paired in a certain
unitary matrix $U\left( B,L\right) $, defined from
the two sets, and to be described in (\ref{eq12}) below.
The unitary matrix $U\left( B,L\right) $
is related to one studied by Hadamard.
It turns out that not all configurations
of sets $B,L$ allow such a unitary
pairing, and there is a further constraint
from the dimension $\nu $
of the ambient Euclidean space.

We now turn to the theory, following
an illustrative example in a single
dimension, i.e., $\nu =1$.

\section{\label{S2}An Example: Scale Four}

It is known (see, e.g., \cite{JoPe93b}, \cite{JoPe96}) that there is a
unique probability measure $\mu $ on $\mathbb{R}$ of compact support such
that 
\begin{equation}
\int f\,d\mu =\frac{1}{2}\left( \int f\left( \frac{x}{4}\right) \,d\mu
\left( x\right) +\int f\left( \frac{x}{4}+\frac{1}{2}\right) \,d\mu \left(
x\right) \right)   \label{eq2}
\end{equation}
for all continuous $f$. In fact, the support $K$ of $\mu $ is the Cantor set
obtained by dividing $I=\left[ 0,1\right] $ into four equal subintervals,
and retaining only the first and third. (See Figure 1 below.) 
\begin{figure}[tbp]
\begin{center}
\ 
\psfig{file=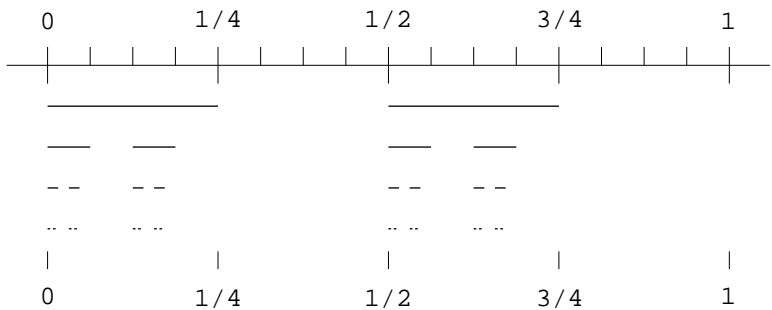,height=116bp,width=288bp}
\end{center}
\caption{Support of $\protect\mu$}
\label{cantor1}
\end{figure}
%
The Hausdorff dimension of $\mu $ is $d_{H}=\frac{\ln 2}{\ln 4}=\frac{1}{2}$.

This is a special case of a more general construction in $\nu $ dimensions ($%
\nu \geq 1$) corresponding to some given real matrix $R$, and a finite
subset $B\subset \mathbb{R}^{\nu }$. It is assumed that 
\begin{equation}
R\text{ has eigenvalues }\xi _{i}\text{ all satisfying }\left| \xi
_{i}\right| >1.  \label{eq3}
\end{equation}
The subset $B$ is required to satisfy an open-set condition: Introduce 
\begin{equation}
\sigma _{b}x=R^{-1}x+b,\quad x\in \mathbb{R}^{\nu }.  \label{eq4}
\end{equation}
It is assumed that there is a nonempty, bounded open set $V$ such that 
\begin{equation}
\bigcup_{b\in B}\sigma _{b}V\subset V  \label{eq5}
\end{equation}
with the union disjoint corresponding to distinct points in $B$. Our present 
$\left\{ \sigma _{b}\right\} $ systems (see below for the axioms) are
special cases of \emph{iterated function systems} (i.f.s.) considered in 
\cite{Hut81}, see also \cite{Fal86}. There are many interesting more general
i.f.s., and that context also leads to measures $\mu $ which satisfy a
general version of the invariance property (\ref{eq6}), and there is then a
corresponding ``open-set assumption''. But for our present \emph{affine}
systems, the splitting property (\ref{eq5}), for some open subset $V$ in $%
\mathbb{R}^{\nu }$, can be shown in fact to be automatic, see \cite{JoPe96}.
In concrete examples, e.g., Figures \ref{cantor1}--\ref{eiffel} below, the
choice of such an open set $V$ is often apparent from the geometry of the
affine maps $\sigma _{b}$. We do not directly use the open-set condition (%
\ref{eq5}) in our present proofs%
, except in the application of
Theorem \ref{Thm8.3} (Remark \ref{RemNew8.4}) where
the interior of a particular simplex is
identified for the open set $V$, and the property
gets used in (\ref{eq8ins3})%
. For general i.f.s., the property is needed
for the computation of geometric invariants, such as the Hausdorff
dimension, see \cite{Fal86}, for the particular iteration limits; and, for
our present systems, it forces the invariant fractal measure $\mu $ to be
locally translation invariant, see \cite[Appendix]{JoPe97}. If $N=\#\left(
B\right) $, then the corresponding measure $\mu $ on $\mathbb{R}^{\nu }$
(depending on $R$ and $B$) has compact support, and satisfies 
\begin{equation}
\int f\,d\mu =\frac{1}{N}\sum_{b\in B}\int f\left( \sigma _{b}\left(
x\right) \right) \,d\mu \left( x\right)  \label{eq6}
\end{equation}
for all continuous $f$. (For more details on the ``open-set condition'' and
affinely generated fractal measures, we give the following background
references: \cite{JoPe94}, \cite{Str94}, \cite{Str95}, and \cite{Ped97}.)
Define, for $t\in \mathbb{R}^{\nu }$, the Fourier transform 
\begin{equation}
\hat{\mu}\left( t\right) =\int e^{i2\pi t\cdot x}\,d\mu \left( x\right)
\label{eq7}
\end{equation}
with $t\cdot x=\sum_{i=1}^{\nu }t_{i}x_{i}$, we then get 
\begin{equation}
\hat{\mu}\left( t\right) =\chi _{B}\left( t\right) \hat{\mu}\left(
R^{*\,-1}t\right)  \label{eq8}
\end{equation}
where 
\begin{equation}
\chi _{B}\left( t\right) :=\frac{1}{N}\sum_{b\in B}e^{i2\pi b\cdot t}
\label{eq9}
\end{equation}
and $R^{*}$ is the transposed matrix. For the above example, this amounts to 
\begin{equation}
\hat{\mu}\left( t\right) =\frac{1}{2}\left( 1+e^{i\pi t}\right) \hat{\mu}%
\left( \frac{t}{4}\right) ,\quad t\in \mathbb{R}.  \label{eq10}
\end{equation}

\section{\label{S3}Orthogonal Frequencies and Fractal Hardy Spaces}

Assume that the matrix $R$ in (\ref{eq4}) has integral entries, and that 
\begin{equation}
RB\subset \mathbb{Z}^{\nu },\quad 0\in B,  \label{eq11}
\end{equation}
but that none of the differences $b-b^{\prime }$ is in $\mathbb{Z}^{\nu }$
when $b,b^{\prime }\in B$ are different. Furthermore, assume that some
subset $L\subset \mathbb{Z}^{\nu }$ satisfies $0\in L$, $\#\left( L\right) =N
$ ($=\#\left( B\right) $), and the matrix 
\begin{multline}
H_{BL}:=N^{-\frac{1}{2}}\left( e^{i2\pi b\cdot l}\right) \text{is unitary as
an }N\times N\text{ complex matrix, }  \label{eq12} \\
\text{i.e., }H_{BL}^{\ast }H_{BL}^{{}}=I_{N}^{{}}\quad \text{(}\,^{\ast }=%
\text{transposed conjugate).}
\end{multline}
In fact, it can be checked that the assumed non-integrality of the
differences $b-b^{\prime }$ (when $\neq 0$) follows from assuming that $%
H_{BL}$ is unitary for some $L$ as described. For our purposes, the
assumptions $L\subset \mathbb{Z}^{\nu }$ and $RB\subset \mathbb{Z}^{\nu }$
may actually be weakened as follows: 
\begin{equation}
\left( R^{n}b\right) \cdot l\in \mathbb{Z},\text{\quad }\forall n\in \mathbb{%
N},\;b\in B,\;l\in L.  \label{eq12bis}
\end{equation}
Details on this more general setup will be given in Section  \ref{S9} below.

\begin{lemma}
\label{Lem3.1}With the assumptions, set 
\begin{equation}
P:=\left\{ l_{0}+R^{\ast }l_{1}+%
\dots%
%
:l_{i}\in L\text{, finite sums}\right\} .  \label{eq13}
\end{equation}
Then the functions $\left\{ e_{\lambda }:\lambda \in P\right\} $ are
mutually orthogonal in $L^{2}\left( \mu \right) $ where 
\begin{equation}
e_{\lambda }\left( x\right) :=e^{i2\pi \lambda \cdot x}.  \label{eq14}
\end{equation}
\end{lemma}

\begin{proof}%
%
Let $\lambda =\sum R^{\ast \,i}l_{i}$, $\lambda ^{\prime }=\sum R^{\ast
\,i}l_{i}^{\prime }$ be points in $P$, and assume $\lambda \neq \lambda
^{\prime }$. Then 
\begin{align*}
\ip{e_{\lambda }}{e_{\lambda ^{\prime }}}_{\mu }& 
=\int \overline{e_{\lambda }}e_{\lambda ^{\prime }}\,d\mu \\
& =\int e^{i2\pi \left( \lambda ^{\prime }-\lambda \right) \cdot x}\,d\mu
\left( x\right) \\
& =\hat{\mu}\left( \lambda ^{\prime }-\lambda \right) \\
& =\hat{\mu}\left( l_{0}^{\prime }-l_{0}^{{}}+R^{\ast }\left( l_{1}^{\prime
}-l_{1}^{{}}\right) +\cdots \right) \\
& =\chi _{B}\left( l_{0}^{\prime }-l_{0}^{{}}\right) \hat{\mu}\left(
l_{1}^{\prime }-l_{1}^{{}}+R^{\ast }\left( l_{2}^{\prime }-l_{2}^{{}}\right)
+\cdots \right) .
\end{align*}%
%
If $l_{0}^{\prime }\neq l_{0}^{{}}$ then $\chi _{B}\left( l_{0}^{\prime
}-l_{0}^{{}}\right) =0$ by (\ref{eq12}). If not, there is a first $n$ such
that $l_{n}^{\prime }\neq l_{n}^{{}}$, and then 
\begin{align*}
\hat{\mu}\left( \lambda ^{\prime }-\lambda \right) & =\hat{\mu}\left(
R^{\ast \,n}\left( l_{n}^{\prime }-l_{n}^{{}}\right) +R^{\ast \,n+1}\left(
l_{n+1}^{\prime }-l_{n+1}^{{}}\right) +\cdots \right) \\
& =\chi _{B}\left(
l_{n}^{\prime }-l_{n}^{{}}\right) \hat{\mu}\left( l_{n+1}^{\prime
}-l_{n+1}^{{}}+\cdots \right) \\
& =0
\end{align*}%
%
since $\chi _{B}\left( l_{n}^{\prime }-l_{n}^{{}}\right) =0$.%
\end{proof}%
%

\begin{corollary}
\label{Cor3.2}Let $\mu $ be the measure on the line $\mathbb{R}$ given by 
\textup{(\ref{eq2})} and with Hausdorff dimension $d_{H}=\frac{1}{2}$. 
\textup{(}We have $R=4$, $B=\left\{ 0,\frac{1}{2}\right\} $ and $L=\left\{
0,1\right\} $.\textup{)} Then 
\begin{equation}
P=\left\{ l_{0}+4l_{1}+4^{2}l_{2}+%
\dots%
%
:l_{i}\in \left\{ 0,1\right\} \text{, finite sums}\right\} ,  \label{eq15}
\end{equation}
and $\left\{ e_{\lambda }:\lambda \in P\right\} $ is an \emph{orthonormal}
subset of $L^{2}\left( \mu \right) $.
\end{corollary}

\begin{proof}%
%
Immediate from the lemma.%
\end{proof}%
%

\begin{lemma}
\label{Lem3.3}Let the subsets $B,L\subset \mathbb{R}^{\nu }$, and the matrix 
$R$ be as described before Lemma \textup{\ref{Lem3.1}}. Let 
\begin{equation}
Q_{1}\left( t\right) :=\sum_{\lambda \in P}\left| \hat{\mu}\left( t-\lambda
\right) \right| ^{2},\quad t\in \mathbb{R}^{\nu }.  \label{eq16}
\end{equation}
Then $\left\{ e_{\lambda }:\lambda \in P\right\} $ is an orthonormal basis
for $L^{2}\left( \mu \right) $ if and only if $Q_{1}\equiv 1$ on $\mathbb{R}%
^{\nu }$.
\end{lemma}

\begin{proof}%
%
If $\left\{ e_{\lambda }:\lambda \in P\right\} $ is an orthogonal basis for $%
L^{2}\left( \mu \right) $, the Bessel inequality is an identity when applied
to $e_{t}$; that is, 
\begin{equation*}
1=\left\| e_{t}\right\| _{\mu }^{2}=\sum_{\lambda }\left|%
\ip{e_\lambda}{e_t}%
%
_{\mu }\right|^{2}=\sum_{\lambda }\left|\hat{\mu}\left( t-\lambda \right)
\right|^{2}.
\end{equation*}
Conversely, if this holds, and if $f\in L^{2}\left( \mu \right) \ominus
\left\{ e_{\lambda }:\lambda \in P\right\} $, then $%
\ip{e_t}{f}%
%
_{\mu }=0$ for all $t\in \mathbb{R}^{\nu }$, or equivalently $\int e^{-i2\pi
t\cdot x}f\left( x\right) \,d\mu \left( x\right) =0$ for all $t\in \mathbb{R}%
^{\nu }$. This implies $f=0$ by Stone--Weierstrass applied to the compact
support $\limfunc{supp}\left( \mu \right) $.%
\end{proof}%
%

\begin{theorem}
\label{Thm3.4}Let $H_{2}\left( P,\mu \right) $ be the closed span in $%
L^{2}\left( \mu \right) $ of the functions $\{e^{i2\pi
nx}:n=0,1,4,5,16,17,20,21,%
\dots%
%
\}$ \textup{(}i.e., $P=\{l_{0}+4l_{1}+4^{2}l_{2}+%
\dots%
%
:l_{i}\in \left\{ 0,1\right\} $, finite sums$\}$\textup{).} Then 
\begin{equation}
H_{2}\left( P,\mu \right) =L^{2}\left( \mu \right) .  \label{eq17}
\end{equation}
\end{theorem}

\begin{remark}
\label{RemNew3.5}While $L^{2}\left( \mu \right) $ is spanned by $\left\{
z^{n}:n\in P\right\} $ ($z=e^{i2\pi x}$), and $P\subset \mathbb{N}_{0}$, it
is of course \emph{not} the case that 
\begin{equation}
\int z^{n}\,d\mu =0\text{\quad for \emph{all} }n\in \mathbb{N}.
\label{Exaggeration}
\end{equation}
(The coefficients do vanish for $n\in \left\{ 1,4,5,16,17,%
\dots%
%
\right\} =P\diagdown \left\{ 0\right\} $, even when $n$ is of the form $%
n=4^{k}\left( 2\mathbb{Z}+1\right) $, $k\in \mathbb{N}_{0}$, but not for $n$
in a bigger subset of $\mathbb{N}$.) First of all, we showed in \cite{JoPe96}
that $\left\{ z^{n}:n\in P\right\} $ is maximally orthogonal; and secondly,
if $\mu $ is viewed as a measure on $\mathbb{T}$, then the validity of (\ref
{Exaggeration}) would imply (by the F. and M. Riesz theorem) absolute
continuity of $\mu $ with respect to Lebesgue measure $\lambda $, i.e., $%
\frac{d\mu }{d\lambda }\in L^{1}\left( \lambda \right) $, which is clearly
false. It is known (e.g., \cite{JoPe94}) that $\mu $ is purely singular.
\end{remark}

\begin{corollary}
\label{CorNew3.6}There is a canonical
isometric embedding $\Phi $ of $L^{2}\left( \mu
\right) $ into the subspace $H_{2}\left( z^{4}\right) +zH_{2}\left(
z^{4}\right) $ of $H_{2}$ where $H_{2}\left( z^{4}\right) :=\left\{ f\left(
z^{4}\right) :f\in H_{2}\right\} $; and it is given by
\begin{equation}
\Phi \left( \vphantom{\sum}\smash{\sum_{\lambda \in P}}
c_{\lambda }e_{\lambda }\right) =\sum_{n\in
P}c_{4n}z^{4n}+z\sum_{n\in P}c_{4n+1}z^{4n}.  \label{eqPhican}
\end{equation}
\end{corollary}

\begin{proof}%
%
Since $\sum_{\lambda \in P}\lvert c_{\lambda }\rvert^{2}<\infty $, and $%
P=\bigcup_{l\in \left\{ 0,1\right\} }l+4P$, with $4P\cap \left( 1+4P\right)
=\varnothing $, the representation (\ref{eqPhican}) is well defined. Note
that $\Phi $ is everywhere defined on $L^{2}\left( \mu \right) $ by Theorem 
\ref{Thm3.4}, and the two functions $f_{0}\left( z\right) =\sum_{n\in
P}c_{4n}z^{4n}$ and $f_{1}\left( z\right) =\sum_{n\in P}c_{4n+1}z^{4n}$ are
in $H_{2}\left( z^{4}\right) $.%
\end{proof}%
%

\begin{remark}
\label{RemNew3.7}(\emph{Fractal Hardy spaces\/}) An iteration
of the argument from the proof of
the corollary yields, for each
$n\in \mathbb{N}$, a natural \emph{isometric}
embedding $\Phi _{n}$ of $L^{2}\left( \mu \right) $
into the subspace of $H_{2}$
characterized as $n$ increases by:
\begin{multline*}
H_{2}\left( z^{4^{n}}\right)
 +zH_{2}\left( z^{4^{n}}\right)
 +z^{4}H_{2}\left( z^{4^{n}}\right)
 +z^{5}H_{2}\left( z^{4^{n}}\right)  \\
 +z^{16}H_{2}\left( z^{4^{n}}\right)
 +z^{17}H_{2}\left( z^{4^{n}}\right)%
\dots%
%
 +z^{\frac{4^{n}-1}{3}}H_{2}\left( z^{4^{n}}\right) .
\end{multline*}
Specifically, let $n\in \mathbb{N}$ be fixed, and
let $P_{n}=\{l_{0}+4l_{1}+%
\dots%
%
+4^{n-1}l_{n-1}:l_{i}\in \left\{ 0,1\right\} \}$.
Then the functions in
$\Phi _{n}\left( L^{2}\left( \mu \right) \right) $ ($\subset H_{2}$) have
the following characteristic module representation:
\begin{equation*}
\vphantom{\sum_{p\in P_{n}}}\left\{ \vphantom{\sum}\smash{\sum_{p\in
P_{n}}}z^{p}f_{p}\left( z^{4^{n}}\right) :f_{p}\in H_{2}\right\} .
\end{equation*}
For each $n$, $\Phi _{n}$ maps \emph{into} this
space, and not onto.
\end{remark}

\begin{remark}
\label{Rem3.5}(\emph{Spectral pairs\/}) Our interest in the problem of
finding dense analytic subspaces in $L^{2}\left( \mu \right) $, for
probability measures, grew out of our earlier work on \emph{spectral pairs},
see, e.g.,\ \cite{JoPe91}, \cite{JoPe92}, \cite{JoPe93a}. Consider subsets $%
\Omega $ and $\Lambda $ in $\mathbb{R}^{\nu }$, with $\Omega $ of finite
positive Lebesgue measure, and let $L^{2}\left( \Omega \right) $ be the
corresponding Hilbert space from the $\Omega $-restricted and normalized $%
\nu $-dimensional Lebesgue measure. Let $e_{\lambda }\left( x\right)
:=e^{i2\pi \lambda \cdot x}$ on $\Omega $. We say that $\left( \Omega
,\Lambda \right) $ is a \emph{spectral pair} if $\left\{ e_{\lambda
}:\lambda \in \Lambda \right\} $ is an orthonormal basis in $L^{2}\left(
\Omega \right) $. It is known \cite{Fug74} that, for $\nu =2$, the case when 
$\Omega $ is either the triangle, or the disk, does \emph{not} admit any
sets $\Lambda $ such that $\left( \Omega ,\Lambda \right) $ is a spectral
pair. On the other hand, our work in \cite{JoPe94} showed that, when $\left(
\Omega ,\Lambda \right) $ \emph{is} a spectral pair in $\nu $ dimensions,
then $\Omega $ is often ``generated'' by some amount of self-affine
structure, as described by a system of affine transforms $\sigma _{b}$ as in
(\ref{eq4}).

Since there is a lack of symmetry of the two sets $\Omega $ and $\Lambda $
in a spectral pair, we introduced in \cite{JoPe97} a generalized spectral
pair formulation for two measures $\mu $ and $\rho $. The context was
locally compact abelian groups: Let $G$ be a locally compact abelian group
with dual group $\Gamma $. Let $\mu $ be a Borel measure on $G$, and $\rho $
one on $\Gamma $. For $f$ of compact support and continuous, introduce 
\begin{equation*}
F_{\mu }f\left( \xi \right) =\int_{G}\overline{\left\langle \xi
,x\right\rangle }f\left( x\right) \,d\mu \left( x\right) 
\end{equation*}
where $\left\langle \xi ,x\right\rangle $ denotes the pairing between points 
$\xi $ in $\Gamma $ and $x$ in $G$. If $f\mapsto F_{\mu }f$ extends to an
isomorphic isometry (i.e., unitary) of $L^{2}\left( \mu \right) $ \emph{onto}
$L^{2}\left( \rho \right) $, then we say that $\left( \mu ,\rho \right) $ is
a \emph{spectral pair}. (The case when $\mu $ is a restriction of Lebesgue
measure was studied in \cite{Ped87}.) It is clear how the earlier definition
of spectral pairs is a special case, even when $G$ is restricted to the
additive group $\mathbb{R}^{\nu }$. But it is not immediate that there are
examples $\left( \mu ,\rho \right) $ of the new spectral pair type which
cannot be reduced to the old one.

Theorem \ref{Thm3.4} shows that this is indeed the case (i.e., that there
are examples): Let $G=\mathbb{R}$, and let $\mu $ be the fractal measure in
Theorem \ref{Thm3.4}. Let 
\begin{equation*}
P=\left\{ l_{0}+4l_{1}+4^{2}l_{2}+%
\dots%
%
:l_{i}\in \left\{ 0,1\right\} ,\text{ finite sums}\right\} =\left\{
0,1,4,5,16,17,%
\dots%
%
\right\} ,%
\notag%
%
\end{equation*}
and let $\rho =\rho _{P}$ be the counting measure of $P$. Then the
conclusion in Theorem \ref{Thm3.4} may be restated to the effect that $%
\left( \mu ,\rho _{P}\right) $ is a spectral pair. This is perhaps
surprising as earlier work on Fourier analysis of fractal measures, see
e.g.\ \cite{Str90c}, \cite{Str93}, and \cite{JoPe95}, suggested a continuity
in the Fourier transform, and also the presence of asymptotic estimates,
rather than exact identities.

It can be shown, as a consequence of \cite[Corollary A.5]{JoPe97} that if $%
\left( \mu ,\rho \right) $ satisfies our spectral-pair property for any
measure $\rho $, then $\rho =\rho _{P}$ for some subset $P\subset \mathbb{R}%
^{\nu }$, i.e., $L^{2}\left( \mu \right) $ has an orthonormal basis of the
form $\left\{ e_{\lambda }:\lambda \in P\right\} $. The basis for this
argument is the \emph{finiteness} of $\mu $, when generated from \textup{(%
\ref{eq6})}.
\end{remark}

\begin{remark}
\label{Rem3.6}(The triadic Cantor measure) The significance of the
assumptions (\ref{eq11})--(\ref{eq12}) on the pair $R,B$ lies in the
identity (\ref{eq20}) below, and also in orthogonality. If, for example, we
work with the more traditional triadic Cantor set, then the results in
Lemmas \ref{Lem3.1} and \ref{Lem3.3} no longer are valid. To see this, take $%
R=3$ and $B=\left\{ 0,\frac{2}{3}\right\} $. Let $\mu _{3}$ denote the
corresponding measure on $\mathbb{R}$ with support equal to the triadic
Cantor set (see Figure \ref{cantor2}). 
\begin{figure}[tbp]
\begin{center}
\ 
\psfig{file=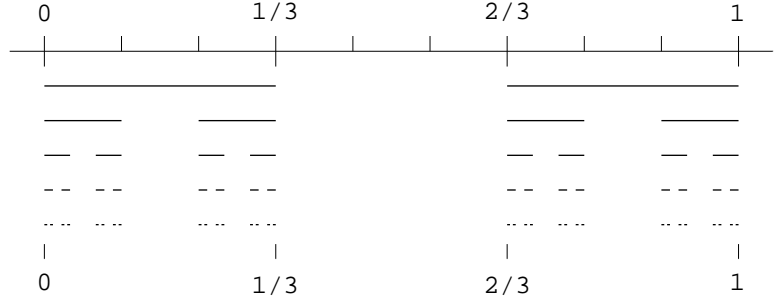,height=114bp,width=288bp}
\end{center}
\caption{Support of $\protect\mu_3$}
\label{cantor2}
\end{figure}
%
It is determined by 
\begin{equation*}
\int f\,d\mu _{3}=\frac{1}{2}\left( \int f\left( \frac{x}{3}\right) \,d\mu
_{3}\left( x\right) +\int f\left( \frac{x}{3}+\frac{2}{3}\right) \,d\mu
_{3}\left( x\right) \right) ,\quad \forall f\in C_{c}\left( \mathbb{R}%
\right) ,%
\notag%
%
\end{equation*}
has $d_{H}=\frac{\ln 2}{\ln 3}$,$\ $and satisfies 
\begin{equation*}
\widehat{\mu _{3}}\left( t\right) =\frac{1}{2}\left( 1+e^{i\frac{4}{3}\pi
t}\right) \widehat{\mu _{3}}\left( \frac{t}{3}\right) ,\quad t\in \mathbb{R}.
\end{equation*}

Choose $L=\left\{ 0,\frac{3}{4}\right\} $ so that (\ref{eq12}) is valid;
then the subset $P_{3}$ (which corresponds to $P=P\left( L\right) $ in Lemma 
\ref{Lem3.3}), is 
\begin{equation*}
P_{3}=\left\{ \frac{3}{4}\left( l_{0}+3l_{1}+3^{2}l_{2}+\cdots \right)
:l_{i}\in \left\{ 0,1\right\} ,\text{ finite sums}\right\} ,
\end{equation*}
but the corresponding exponentials $\left\{ e_{\lambda }:\lambda \in
P_{3}\right\} $ are now \emph{not} mutually orthogonal in $L^{2}\left( \mu
_{3}\right) $. Take for example the two points $\lambda =\frac{3}{4}$ and $%
\lambda ^{\prime }=\frac{9}{4}$ both in the set $P_{3}$. The corresponding
exponentials $e_{\lambda }$ and $e_{\lambda ^{\prime }}$ are both orthogonal
to $e_{0}$, but they are not mutually orthogonal, i.e., $%
\ip{e_{\lambda }}{e_{\lambda ^{\prime }}}%
%
_{\mu _{3}}\neq 0$. In fact, for the $\mu _{3}$-inner product: 
\begin{equation*}
\ip{e_{\lambda }}{e_{\lambda ^{\prime }}}%
%
_{\mu _{3}}=\widehat{\mu _{3}}\left( \frac{3}{2}\right) =\widehat{\mu _{3}}%
\left( \frac{1}{2}\right) =\frac{1}{4}\widehat{\mu _{3}}\left( \frac{1}{6}%
\right) \neq 0.
\end{equation*}
It can further be shown (see \cite{JoPe97}
and Section \ref{SNew6} below%
) that there is \emph{no subset} $%
P\subset \mathbb{R}$ such that, if $\rho _{P}$ denotes the corresponding
counting measure on $\mathbb{R}$, then $\left( \mu _{3},\rho _{P}\right) $
is a spectral pair in the above more general sense (\emph{a fortiori,}
fractions don't provide a basis either). Similarly, it can be checked that
the identity (\ref{eq20}) in Lemma \ref{Lem4.1} below fails for this pair $%
\mu _{3},P_{3}$, i.e., the triadic one with $d_{H}=\frac{\ln 2}{\ln 3}$.
\end{remark}

\section{\label{S4}Proof of Theorem \ref{Thm3.4}: Ruelle's Transfer Operator}

We need only verify the assumption in Lemma \ref{Lem3.3}. But when applied
to the example in Section \ref{S2}, the issue becomes showing that $%
\sum_{n\in P}\left| \hat{\mu}\left( t-n\right) \right| ^{2}=1$, $t\in 
\mathbb{R}$, where the summation in $n$ is over 
\begin{equation}
P=P\left( L\right) =\left\{ l_{0}+4l_{1}+4^{2}l_{2}+%
\dots%
%
:l_{i}\in \left\{ 0,1\right\} \text{, finite sums}\right\} .  \label{eq18}
\end{equation}
More generally, let the sets $B,L\subset \mathbb{R}^{\nu }$, and the matrix $%
R$, be as in Lemma \ref{Lem3.1}.

\begin{lemma}
\label{Lem4.1}The function 
\begin{equation}
Q_{1}\left( t\right) :=\sum_{\lambda \in P}\left| \hat{\mu}\left( t-\lambda
\right) \right| ^{2}  \label{eq19}
\end{equation}
\textup{(}where $P=\{l_{0}+R^{\ast }l_{1}+%
\dots%
%
:l_{i}\in L$, finite sums$\}$\textup{)} satisfies the functional identity 
\textup{(}$t\in \mathbb{R}^{\nu }$\textup{)} 
\begin{equation}
Q\left( t\right) =\sum_{l\in L}\left| \chi _{B}\left( t-l\right) \right|
^{2}\,Q\left( R^{\ast \,-1}\left( t-l\right) \right) .  \label{eq20}
\end{equation}
\end{lemma}

\begin{proof}%
%
Let $t\in \mathbb{R}^{\nu }$. Then 
\begin{align*}
Q_{1}\left( t\right) & =\sum_{\lambda \in P}\left|\hat{\mu}\left( t-\lambda
\right) \right|^{2} \\
& =\sum_{\lambda \in P}\left|\chi _{B}\left( t-\lambda \right)
\right|^{2}\left| \hat{\mu}\left( R^{\ast \,-1}\left( t-\lambda \right)
\right) \right| ^{2} \\
& =\sum_{l\in L}\sum_{\lambda \in P}\left|\chi _{B}\left( t-l-R^{\ast
}\lambda \right) \right|^{2}\left| \hat{\mu}\left( R^{\ast \,-1}\left(
t-l\right) -\lambda \right) \right| ^{2} \\
& =\sum_{l\in L}\left|\chi _{B}\left( t-l\right) \right|^{2}\sum_{\lambda
\in P}\left| \hat{\mu}\left( R^{\ast \,-1}\left( t-l\right) -\lambda \right)
\right| ^{2} \\
& =\sum_{l\in L}\left|\chi _{B}\left( t-l\right) \right|^{2}\,Q_{1}\left(
R^{\ast \,-1}\left( t-l\right) \right) . 
\settowidth{\qedskip}{$\displaystyle 
               Q_{1}\left( t\right) =
               \sum_{l\in L}\left|\chi _{B}\left( t-l\right) 
\right|^{2}Q_{1}\left( R^{\ast
               \,-1}\left( t-l\right) \right) .$}
               \settowidth{\qedadjust}{$\displaystyle 
               Q_{1}\left( t\right) =
               \sum_{l\in L}\sum_{\lambda \in P}
\left|\chi _{B}\left( t-l-R^{\ast
               }\lambda \right) \right|^{2}
\left| \hat{\mu}\left( R^{\ast \,-1}\left(
               t-l\right) -\lambda \right) \right| ^{2}$}
               \addtolength{\qedadjust}{\textwidth}
               \addtolength{\qedskip}{-0.5\qedadjust}
               \rlap{\hbox to-\qedskip{\hfil $\qedsymbol $}\hss }
\end{align*}%
%
\renewcommand{\qed}{\relax} 
\end{proof}%
%

\begin{remark}
\label{Rem4.2}We shall need the following operator $C$,
\emph{Ruelle's transfer operator},
defined on functions 
$Q$: 
\begin{equation}
C\left( Q\right) \left( t\right) :=\sum_{l\in L}\left| \chi _{B}\left(
t-l\right) \right| ^{2}\,Q\left( R^{\ast \,-1}\left( t-l\right) \right) .
\label{eq21}
\end{equation}
\end{remark}

We also note

\begin{lemma}
\label{LemNew4.3}\ 
\begin{enumerate}
\item \label{LemNew4.3(1)}The function $Q_{1}$,
defined on $\mathbb{R}^{\nu }$ by \textup{%
(\ref{eq19})}, has an entire analytic extension to $\mathbb{C}^{\nu }$
which is of linear exponential growth in the imaginary direction, i.e., with 
$\mathbb{C}^{\nu }=\mathbb{R}^{\nu }+i\mathbb{R}^{\nu }$,
\begin{equation}
\left| Q_{1}\left( t+is\right) \right| \leq e^{4\pi m\left\| s\right\| _{2}}
\label{eqLemNew4.3(1)}
\end{equation}
\textup{(}$t,s\in \mathbb{R}^{\nu }$\textup{)} where $m$ is the diameter of
the support of $\mu $ inside $\mathbb{R}^{\nu }$.

\item \label{LemNew4.3(2)}
The operator $C$ in \textup{(\ref{eq21})}
extends naturally to the space $\mathcal{E}%
_{m}$ of entire analytic functions $f$ on $\mathbb{C}^{\nu }$ which satisfy
\textup{(\ref{eqLemNew4.3(1)})}
and $f\left( 0\right) =1$, and maps this space into
itself, i.e., $C\colon \mathcal{E}_{m}\rightarrow \mathcal{E}_{m}$.
\end{enumerate}
\end{lemma}

\begin{proof}%
%
For $t\in \mathbb{R}^{\nu }$, 
\begin{equation*}
Q_{1}\left( t\right) 
=\sum_{\lambda \in P}\hat{\mu}\left( t-\lambda \right) \hat{%
\mu}\left( \lambda -t\right) =\sum_{\lambda \in P}%
\ip{e_{\lambda }}{e_{t}}%
%
_{\mu }%
\ip{e_{-\lambda }}{e_{-t}}%
%
_{\mu }.
\end{equation*}
If the right-hand side is truncated, summing only over the finite subsets in 
$P$ given by the restriction $\left| \lambda \right| \leq n$, $n=1,2,%
\dots%
%
$, then each of the resulting finite sums $\sum_{\lambda \in P,\left|
\lambda \right| \leq n}\cdots $ defines an entire analytic function in $t$
(i.e., $t$ is now allowed to range over $\mathbb{C}^{\nu }$, with the
convention $e_{t}\left( x\right) =e^{i2\pi t\cdot x}$, and $t\cdot
x=\sum_{j=1}^{\nu }t_{j}x_{j}$, $t_{j}\in \mathbb{C}$, $x_{j}\in \mathbb{R}$%
). These functions are also in $L^{2}\left( \mu \right) $ since the measure $%
\mu $ from (\ref{eq6}) has compact support $\limfunc{supp}\left( \mu \right) 
$ in $\mathbb{R}^{\nu }$. Hence 
\begin{align*}
\left| 
\vphantom{\sum}\smash{\sum_{
\substack{\lambda \in P \\
\left| \lambda \right| \leq n}}}
\cdots \right|
\leq 
\vphantom{\sum}\smash{\sum_{
\substack{\lambda \in P \\
\left| \lambda \right| \leq n}}}
\!\left| 
\,\cdots \right|
& \leq \sum_{\lambda }\left| 
\ip{e_{\lambda }}{e_{t}}
\right| \left| 
\ip{e_{-\lambda }}{e_{-t}}
\right| \\
& \leq \vphantom{\sum_{\lambda }}
\left( \vphantom{\sum}\smash{\sum_{\lambda }}\left| 
\ip{e_{\lambda }}{e_{t}}
\right| ^{2}\right) ^{\frac{1}{2}}
\vphantom{\sum_{\lambda }}
\left( \vphantom{\sum}\smash{\sum_{\lambda }}\left| 
\ip{e_{-\lambda }}{e_{-t}}
\right| ^{2}\right) ^{\frac{1}{2}} \\
& \leq \left\| e_{t}\right\| _{L^{2}\left( \mu \right) }\left\|
e_{-t}\right\| _{L^{2}\left( \mu \right) }.
\end{align*}
Note this holds for the $P$ summation, for all $n$. Letting $n\rightarrow
\infty $, we then also get the estimate for the unrestricted summation over $%
\lambda \in P$. This may be formulated as a pointwise approximation 
\begin{equation*}
\lim_{n\rightarrow \infty }Q_{n}\left( t\right) =\tilde{Q}\left( t\right)
,\quad t\in \mathbb{C}^{\nu },
\end{equation*}
the limit now providing an extension of $Q\left( \,\cdot \,\right) $,
initially given only on $\mathbb{R}^{\nu }$ by (\ref{eq19}). Suppose the
support of $\mu $ is contained in $\left\{ x\in \mathbb{R}^{\nu }:\left|
x\right| \leq m\right\} $; then 
\begin{equation*}
\left\| e_{t}\right\| _{L^{2}\left( \mu \right) }^{2}=\int \left| e^{i2\pi
t\cdot x}\right| ^{2}\,d\mu \left( x\right) =\int_{
\rlap{$\scriptstyle \left| x\right| \leq m$}
}\;e^{-4\pi \left( \func{Im}t\right) \cdot x}\,d\mu \left( x\right) \leq
e^{4\pi m\left\| \func{Im}t\right\| }
\end{equation*}
where 
\begin{equation*}
\left\| \func{Im}t\right\| =\vphantom{\sum_{j=1}^{\nu }}
\left( \vphantom{\sum}\smash{\sum_{j=1}^{\nu }}\left| \func{Im}%
t_{j}\right| ^{2}\right) ^{\frac{1}{2}},\quad t\in \mathbb{C}^{\nu }.
\end{equation*}
This means that the functions $Q_{n}\left( \,\cdot \,\right) $ are uniformly
bounded in $t\in \mathbb{C}^{\nu }$, restricted to each strip $\left\| \func{%
Im}t\right\| \leq \limfunc{constant}$. We then conclude by the theorems of
Montel and Vitali (see \cite[p.\ 143]{Neh75}) that the limit $\tilde{Q}%
\left( t\right) $ is entire analytic, $t\in \mathbb{C}^{\nu }$, and serves
as an analytic extension of
$Q_{1}$ from (\ref{eq19}), given only on $\mathbb{R}%
^{\nu }$.

The technique which we used above in estimating the extended $Q\left(
t+is\right) $ in $\mathbb{C}^{\nu }=\mathbb{R}^{\nu }+i\mathbb{R}^{\nu }$
also applies, \emph{mutatis mutandis}, on the other term, to yield:
\begin{equation}
\sum_{l\in L}\left| \chi _{B}\left( t+is-l\right) \right| ^{2}\leq
N^{-1}\sum_{b\in B}e^{4\pi \left\| b\right\| _{2}\left\| s\right\| _{2}}
\label{mutatismutandis}
\end{equation}
by virtue of (\ref{eq12}). When combined with (\ref{eqLemNew4.3(1)}),
defining $\mathcal{E}_{m}$, and (\ref{eq21}), the invariance, $C\colon 
\mathcal{E}_{m}\rightarrow \mathcal{E}_{m}$ follows.
We include the details of proof for estimate (\ref{mutatismutandis})
(using unitarity of the matrix in (\ref{eq12})):
\begin{align*}
\sum _{l\in L}\left| \chi _{B}\left( t+is-l\right) \right| ^{2}
&=N^{-1}\sum _{l\in L}\left| \vphantom{\sum}\smash{\sum _{b\in B}}
N^{-\frac{1}{2}}e^{-i2\pi b\cdot l}
e^{i2\pi b\cdot \left( t+is\right) }\right| ^{2}  \\
&=N^{-1}\sum _{b\in B}\left| 
e^{i2\pi b\cdot \left( t+is\right) }\right| ^{2}  \\
&=N^{-1}\sum _{b\in B}e^{-4\pi b\cdot s}  \\
&\leq N^{-1}\sum _{b\in B}
e^{4\pi \left\| b\right\| _{2}\left\| s\right\| _{2}}.
\settowidth{\qedskip}{$\displaystyle 
\sum _{l\in L}\left| \chi _{B}\left( t+is-l\right) \right| ^{2}
\leq N^{-1}\sum _{b\in B}
e^{4\pi \left\| b\right\| _{2}\left\| s\right\| _{2}}.$}
\settowidth{\qedadjust}{$\displaystyle 
\sum _{l\in L}\left| \chi _{B}\left( t+is-l\right) \right| ^{2}
=N^{-1}\sum _{l\in L}\left| \vphantom{\sum}\smash{\sum _{b\in B}}
N^{-\frac{1}{2}}e^{-i2\pi b\cdot l}
e^{i2\pi b\cdot \left( t+is\right) }\right| ^{2}$}
\addtolength{\qedadjust}{\textwidth}
\addtolength{\qedskip}{-0.5\qedadjust}
\rlap{\hbox to-\qedskip{\hfil $\qedsymbol $}\hss }
\end{align*}
\renewcommand{\qed}{\relax}\end{proof}

\begin{lemma}
\label{Lem4.3}Let $\mu $ be a probability measure with compact support on $%
\mathbb{R}^{\nu }$, \linebreak and let $P\subset \mathbb{R}^{\nu }$ be such
that $\left\{ e_{\lambda }:\lambda \in P\right\} $ is orthogonal in $%
L^{2}\left( \mu \right) $. Set
$Q_{1}\left( t\right) =\sum_{\lambda \in P}\left| 
\hat{\mu}\left( t-\lambda \right) \right| ^{2}$. Then $Q_{1}$ is $C^{\infty }$
on $\mathbb{R}^{\nu }$; in fact, it has an entire analytic extension. Let $%
A\colon L^{2}\left( \mu \right) \rightarrow L^{2}\left( \mu \right) $ be the
orthogonal projection onto $H_{2}\left( P,\mu \right) $. Then

\begin{enumerate}
\item  \label{Lem4.3(1)}$Q_{1}\left( t\right) =\left\| Ae_{t}\right\|
_{L^{2}\left( \mu \right) }^{2}$,

\item  \label{Lem4.3(2)}$\frac{\partial \,}{\partial t_{j\mathstrut }}%
Q_{1}\left( t\right) |_{t=0}=0$,

\item  \label{Lem4.3(3)}$\frac{\partial ^{2\mathstrut }\,}{\partial
t_{j\mathstrut }^{2\mathstrut }}Q_{1}\left( t\right) |_{t=0}=\left\|
Ax_{j}\right\| _{\mu }^{2}-\left\| x_{j}\right\| _{\mu }^{2}$, where $%
\left\| x_{j}\right\| _{\mu }^{2}=\int x_{j}^{2}\,d\mu \left( x\right) $, and

\item  \label{Lem4.3(4)}$\frac{\partial ^{2n\mathstrut }\,}{\partial
t_{j}^{2n\mathstrut }}Q_{1}\left( t\right) |_{t=0}=\left( -1\right) ^{n}\left(
\left\| x_{j}^{n}\right\| _{\mu }^{2}-\left\| Ax_{j}^{n}\right\| _{\mu
}^{2}\right) $.
\end{enumerate}
\end{lemma}

\begin{proof}%
%
Compute!%
\end{proof}%
%

\begin{corollary}
\label{Cor4.4}Let $\left\{ e_{\lambda }:\lambda \in P\right\} $ be a set of
orthogonal frequencies in $L^{2}\left( \mu \right) $ \linebreak and let $%
Q_{1}\left( t\right) =\sum_{\lambda \in P}\left| \hat{\mu}\left( t-\lambda
\right) \right| ^{2}$, $t\in \mathbb{R}^{\nu }$. 
Assume $L$ is not contained in a hyperplane in $\mathbb{R}^{\nu }$.   
Then the set $P$ is total if
and only if 
\begin{equation}
\left( \frac{\partial \,}{\partial t_{j}}\right) ^{2n}Q_{1}\left( t\right)
|_{t=0}=0,\quad j=1,%
\dots%
%
,\nu ,\quad n=1,2,%
\dots%
%
.  \label{eq22}
\end{equation}
\end{corollary}

\begin{proof}%
%
By Lemmas \ref{Lem3.3} and \ref{Lem4.3}, we only need to note that the
monomials are dense in $L^{2}\left( \mu \right) $, by Stone--Weierstrass,
and that a vector $f\in L^{2}\left( \mu \right) $ is in the subspace $%
AL^{2}\left( \mu \right) $ if and only if $\left\| Af\right\| _{\mu
}=\left\| f\right\| _{\mu }$.
(For more details, see Section \ref{S9} below.)%
\end{proof}%
%

\section{\label{S5}The Example Revisited}

Let $\nu =1$, $B=\left\{ 0,\frac{1}{2}\right\} $, $L=\left\{ 0,1\right\} $,
and $R=4$. Let $\mu $ be the corresponding measure ($d_{H}\left( \mu \right)
=\frac{1}{2}$), and let 
\begin{equation}
P=P\left( L\right) =\left\{ l_{0}+4l_{1}+4^{2}l_{2}+%
\dots%
%
:l_{i}\in \left\{ 0,1\right\} \text{, finite sums}\right\} .  \label{eq23}
\end{equation}

\begin{lemma}
\label{Lem5.1}The function 
\begin{equation*}
Q_{1}\left( t\right) =\sum_{\lambda \in P}\left| \hat{\mu}\left( t-\lambda
\right) \right| ^{2}
\end{equation*}
satisfies the following identity: 
\begin{equation}
Q_{1}\left( t\right) =\cos ^{2}\left( \frac{\pi t}{2}\right) 
Q_{1}\left( \frac{t}{4}%
\right) +\sin ^{2}\left( \frac{\pi t}{2}\right) 
Q_{1}\left( \frac{t-1}{4}\right)
.  \label{eq24}
\end{equation}
Let $C$ be the operator defined by the right-hand side in \textup{(\ref{eq24}%
)}. On the convex space $W_{0}:=\left\{ Q\in C^{1}\left( \left[ -\frac{1}{3}%
,0\right] \right) :Q\left( 0\right) =1\right\} $, introduce the metric 
\begin{equation}
d\left( Q_{1},Q_{2}\right) :=\sup \left\{ \left| Q_{1}^{\prime }\left(
t\right) -Q_{2}^{\prime }\left( t\right) \right| :t\in \left[ \textstyle-%
\frac{1}{3},0\right] \right\} .  \label{eq25}
\end{equation}
Then $\left( W_{0},d\right) $ is a complete metric space, and
\begin{equation}
d\left( C\left( Q_{1}\right) ,C\left( Q_{2}\right) \right) \leq \left( \frac{%
1}{4}+\frac{\pi \sqrt{3}}{16}\right) d\left( Q_{1},Q_{2}\right) ,\quad
\forall Q_{1},Q_{2}\in W_{0}.  \label{eq26}
\end{equation}
\end{lemma}

\begin{proof}%
%
First note that the functional
identity (\ref{eq24}) in the lemma is
a special case of (\ref{eq20}) in
Lemma \ref{Lem4.1} above. It is clear
that $C\mathbf{1}=\mathbf{1}$ where $\mathbf{1}$ denotes
the constant function. In the general
case of  Lemma \ref{Lem4.1}, that follows
from (\ref{eq12}).

For $Q_{1},Q_{2}\in W_{0}$, we have $Q_{1}\left( t\right) -Q_{2}\left(
t\right) =\int_{0}^{t}\left( Q_{1}^{\prime }\left( s\right) -Q_{2}^{\prime
}\left( s\right) \right) \,ds$, and therefore $\lvert Q_{1}\left( t\right)
-Q_{2}\left( t\right) \rvert\leq \frac{1}{3}d\left( Q_{1},Q_{2}\right) $. In
particular, $\left( W_{0},d\right) $ is a complete metric space. We have $C$
inducing an operator in the convex set $W_{0}$, and 
\begin{align*}
C\left( Q\right) ^{\prime }\left( t\right) & =\frac{1}{4}\left( \cos
^{2}\left( \frac{\pi t}{2}\right) Q^{\prime }\left( \frac{t}{4}\right) +\sin
^{2}\left( \frac{\pi t}{2}\right) Q^{\prime }\left( \frac{t-1}{4}\right)
\right) \\
& \qquad +\frac{\pi }{2}\left( -\sin \pi t\right) Q\left( \frac{t}{4}\right)
+\frac{\pi }{2}\left( \sin \pi t\right) Q\left( \frac{t-1}{4}\right) \\
& =\frac{1}{4}\left( \cos ^{2}\left( \frac{\pi t}{2}\right) Q^{\prime
}\left( \frac{t}{4}\right) +\sin ^{2}\left( \frac{\pi t}{2}\right) Q^{\prime
}\left( \frac{t-1}{4}\right) \right) \\
& \qquad +\frac{\pi }{2}\left( \sin \pi t\right) \int_{\frac{t}{4}}^{\frac{
t-1}{4}}Q^{\prime }\left( s\right) \,ds.
\end{align*}%
%
For $Q_{1},Q_{2}\in W_{0}$, the first term is estimated above by $\frac{1}{4}%
\sup_{t\in \left[ -\frac{1}{3},0\right] }\left| Q_{1}^{\prime }\left(
t\right) -Q_{2}^{\prime }\left( t\right) \right| $ since the two affine
transformations $t\mapsto \frac{t}{4}$ and $t\mapsto \frac{t-1}{4}$ leave
the interval $\left[ -\frac{1}{3},0\right] $ invariant. The second term is
bounded by 
\begin{equation*}
\frac{\pi }{2}\cdot \frac{\sqrt{3}}{2}\int_{\frac{t-1}{4}}^{\frac{t}{4}%
}\left| Q_{1}^{\prime }\left( s\right) -Q_{2}^{\prime }\left( s\right)
\right| \,ds
\end{equation*}
where we estimate $\left| \sin \pi t\right| $ by $\frac{\sqrt{3}}{2}$ in the
interval $t\in \left[ -\frac{1}{3},0\right] $. Since the integral is bounded
by $\frac{1}{4}\sup_{s\in \left[ -\frac{1}{3},0\right] }\left| Q_{1}^{\prime
}\left( s\right) -Q_{2}^{\prime }\left( s\right) \right| $, the result
follows.%
\end{proof}%
%

\begin{corollary}
\label{Cor5.2}Since the constant function $Q\equiv 1$ in $\left[ -\frac{1}{3}%
,0\right] $ satisfies $C\left( Q\right) =Q$ by \textup{(\ref{eq24})}, we
conclude that $Q_{1}$ from \textup{(\ref{eq19})} must be constant in $\left[ -%
\frac{1}{3},0\right] $ by the Banach fixed-point principle. Theorem \textup{%
\ref{Thm3.4}} now follows by Lemma \textup{\ref{Lem5.1}}.
\end{corollary}

\section{\label{SNew6}One Dimension: Even and Odd Scales}

If a subset $P\subset \mathbb{R}$ is \emph{orthogonal} in the sense that $%
\left\{ e_{\lambda }:\lambda \in P\right\} $ is orthogonal in $L^{2}\left(
\mu \right) $, for some fixed probability measure $\mu $ on $\mathbb{R}$,
then so is the reflected set $-P$, and the translates $t+P$ when $t\in 
\mathbb{R}$ is fixed. In reviewing orthogonal subsets $P$ we shall therefore
impose the condition $0\in P$ for simplicity.

We now examine the following three cases for $\mu $ as they illustrate the
nature of the assumptions made above, as well as the variety of the possible
cases: For the subset $B=\left\{ 0,\frac{1}{2}\right\} $, we have considered
the two scales $R=2$ and $R=4$, and also $R=3$ separately in (\ref{eqj})
below. The resulting measures $\mu _{2}$, $\mu _{4}$, and $\mu _{3}$ are
determined by their Fourier transforms, as follows, (\ref{eqa}), (\ref{eqb}%
), and (\ref{eqj}): 
\begin{align}
\widehat{\mu _{2}}\left( t\right) & =\prod_{n=0}^{\infty }\frac{1}{2}\left(
1+e^{i\frac{\pi t}{2^{n}}}\right) =\prod_{n=1}^{\infty }e^{i\frac{\pi t}{
2^{n}}}\cos \left( \frac{\pi t}{2^{n}}\right)   \label{eqa} \\
& =e^{i\pi t\sum_{n=1}^{\infty }\frac{1}{2^{n}}}\prod_{n=1}^{\infty }\cos
\left( \frac{\pi t}{2^{n}}\right) =e^{i\pi t}\frac{\sin \left( \pi t\right) 
}{\pi t},  \notag
\end{align}%
%
where in the last step we used a familiar infinite product formula for $\cos
\left( \frac{x}{2^{n}}\right) $. The infinite product coincides with the
integral $\int_{0}^{1}e^{i2\pi xt}\,dx=\frac{e^{i2\pi t}-1}{i2\pi t}$ since $%
\mu _{2}$ is Lebesgue measure on $I$. Similarly, 
\begin{equation}
\widehat{\mu _{4}}\left( t\right) =\prod_{n=0}^{\infty }\frac{1}{2}\left(
1+e^{i\frac{\pi t}{4^{n}}}\right) =e^{i\pi \frac{2t}{3}}\prod_{n=0}^{\infty
}\cos \left( \frac{\pi t}{2\cdot 4^{n}}\right) .  \label{eqb}
\end{equation}
(The $\widehat{\mu _{3}}$ product will be discussed separately.) We picked
the set $L=\left\{ 0,1\right\} $ ($\subset \mathbb{Z}$) to satisfy condition
(\ref{eq12}) above. Then from (\ref{eqa})--(\ref{eqb}) we get: 
\begin{align*}
P_{2}& =\left\{ l_{0}+2l_{1}+2^{2}l_{2}+
\dots
:l_{i}\in \left\{ 0,1\right\} \text{, finite sums}\right\} 
=\mathbb{N}_{0}=\left\{ 0,1,2,
\dots
\right\}  \\
\intertext{and}
P_{4}& =\left\{ l_{0}+4l_{1}+4^{2}l_{2}+
\dots
:l_{i}\in \left\{ 0,1\right\} \text{, finite sums}\right\} 
=\left\{ 0,1,4,5,16,17,
\dots
\right\} .
\end{align*}%
%
We get explicit formulas for the zero sets of the transforms, 
\begin{equation}
\mathbf{Z}\left( \widehat{\mu _{2}}\right) =\mathbb{Z}\diagdown \left\{
0\right\} \text{\quad (= the nonzero integers),}  \label{eqc}
\end{equation}
and 
\begin{multline}
\mathbf{Z}\left( \widehat{\mu _{4}}\right) =\left\{ 4^{n}\cdot \left( 1+2%
\mathbb{Z}\right) :n=0,1,2,%
\dots%
%
\right\}   \label{eqd} \\
\text{(= powers of }4\text{ times odd integers)}
\end{multline}
where $\mu _{2}$ is Lebesgue measure, supported on $I=\left[ 0,1\right] $,
while $\mu _{4}$ is the fractal measure with $d_{H}=\frac{1}{2}$ and support
on the Cantor set from Figure \ref{cantor1}. In the first case, it is
classical that $H_{2}\left( P_{2},\mu _{2}\right) $ is the familiar Hardy
space in $L^{2}\left( \mu _{2}\right) $. That is because $L^{2}\left( \mu
_{2}\right) $ may be identified with $L^{2}$ of the circle, and we are back
to classical Fourier series, i.e., $H_{2}$ spanned by $1,z,z^{2},%
\dots%
%
$ in $L^{2}\left( \mathbb{T}\right) $. The contrast between the two cases is
made clear from a comparison of the respective functions $Q_{2}$ and $Q_{4}$
(derived from (\ref{eq16})) as follows: 
\begin{equation*}
Q_{2}\left( t\right) =\sum_{n\in P_{2}}\left| \widehat{\mu _{2}}\left(
t-n\right) \right| ^{2}=\sum_{n=0}^{\infty }\left| \frac{\sin \pi \left(
t-n\right) }{\pi \left( t-n\right) }\right| ^{2}=\frac{\sin ^{2}\left( \pi
t\right) }{\pi ^{2}}\sum_{n=0}^{\infty }\frac{1}{\left( t-n\right) ^{2}},
\end{equation*}
where we provide the following interpretation of the last expression (see 
\cite{Art64} and \cite{Car95}): Let $\Gamma \left( t\right) $ be the gamma
function, with its meromorphic extension resulting from the functional
equation 
\begin{equation*}
\Gamma \left( t\right) \Gamma \left( 1-t\right) =\frac{\pi }{\sin \left( \pi
t\right) },
\end{equation*}
recalling 
\begin{equation*}
\Gamma \left( t\right) =\int_{0}^{\infty }e^{-x}x^{t-1}\,dx
\end{equation*}
for positive $t$. When $t>0$, 
\begin{equation*}
\frac{d\,}{dt}\left( \frac{\Gamma ^{\prime }\left( t\right) }{\Gamma \left(
t\right) }\right) =\sum_{n=0}^{\infty }\frac{1}{\left( t+n\right) ^{2}}%
=:\psi \left( t\right) 
\end{equation*}
is well defined, and negative $t$ is from the meromorphic extension. It
follows that 
\begin{equation}
Q_{2}\left( t\right) =\left\| A_{2}e_{t}\right\| _{L^{2}}^{2}=\frac{\sin
^{2}\left( \pi t\right) }{\pi ^{2}}\psi \left( -t\right) ,  \label{eqe}
\end{equation}
where $A_{2}$ denotes the projection in $L^{2}$ onto the Hardy space $H_{2}$%
. (Note that it is not immediate from this formula (\ref{eqe}) that $%
Q_{2}\left( t\right) $ is entire analytic.) In contrast, for $Q_{4}\left(
t\right) $ we have 
\begin{equation*}
Q_{4}\left( t\right) =\left\| A_{4}e_{t}\right\| _{L^{2}\left( \mu
_{4}\right) }^{2}\equiv 1\text{;}
\end{equation*}
i.e., as noted in Theorem \ref{Thm3.4}, $A_{4}=I$, or, equivalently, the $%
\mu _{4}$-Hardy space is all of $L^{2}\left( \mu _{4}\right) $, and embedded
as a subspace of $H_{2}$ via 
\begin{equation*}
H_{2}\left( z^{4}\right) \oplus zH_{2}\left( z^{4}\right) .
\end{equation*}
(See Corollary \ref{CorNew3.6} for details.)

In the $\mu _{2}$ case, the operator $C$, which is analogous to $C$ for $\mu
_{4}$ in (\ref{eq24}), is 
\begin{equation}
C\left( Q\right) \left( t\right) =\cos ^{2}\left( \frac{\pi t}{2}\right)
Q\left( \frac{t}{2}\right) +\sin ^{2}\left( \frac{\pi t}{2}\right) Q\left( 
\frac{t-1}{2}\right) .  \label{eqf}
\end{equation}
In the $\mu _{2}$ case, the convex $W_{0}$ is 
\begin{equation}
W_{0}=\left\{ Q\in C^{1}\left( \left[ -1,0\right] \right) :Q\left( 0\right)
=1\right\} .  \label{eqg}
\end{equation}
It can be checked that $C$ is then \emph{not} strictly contractive in $W_{0}$
relative to the metric 
\begin{equation}
d\left( Q_{1},Q_{2}\right) =\sup_{t\in \left[ -1,0\right] }\left|
Q_{1}^{\prime }\left( t\right) -Q_{2}^{\prime }\left( t\right) \right| ,
\label{eqh}
\end{equation}
and, as noted, $C$ has fixed points in $W_{0}$ other than the constant
function $Q\equiv 1$. The contractivity constant can be checked in this case
(i.e., (\ref{eqf}), (\ref{eqg}), (\ref{eqh}) corresponding to $R=2$) to be $%
\frac{1}{2}+\frac{\pi }{4}$ which is $>1$ in contrast to the $R=4$ case in (%
\ref{eq26}) where the contractivity constant is $\frac{1}{4}+\frac{\pi \sqrt{%
3}}{16}$ ($<1$). This illustrates the assumption 
\begin{equation}
\#\left( B\right) <\left| \det R\right|  \label{eqi}
\end{equation}
which will be placed on the general systems from Section \ref{S6} below.

Set $z=e^{i2\pi t}$, and $g\left( z\right) =\cos ^{2}\left( \pi t\right) =%
\frac{1}{2}+\frac{1}{4}z+\frac{1}{4}\bar{z}$. Then (\ref{eqf}) may be
rewritten in complex form as 
\begin{equation*}
C\left( Q\right) \left( z\right) =\sum_{%
\substack{w\\w^{2}=z}%
%
}g\left( w\right) Q\left( w\right) ,
\end{equation*}
which is a special case of the familiar Perron--Frobenius--Ruelle operator 
\cite{Rue94}. The presence of the above mentioned two independent solutions
(in $W_{0}$) to $C\left( Q\right) =Q$ is \emph{not} predicted by Ruelle's
theory. The Perron--Frobenius eigenvalue problem $C\left( Q\right) =Q$ is
actually studied in wavelet theory, see, e.g., \cite[Chapter 6]{Dau92}.
Since $\left\{ z\in \mathbb{T}:g\left( z\right) =1\right\} $ is just a
singleton there is only one periodic solution $Q$ in $W_{0}$. This is
consistent with the second solution from (\ref{eqe}) above being
non-periodic.

We also considered $B=\left\{ 0,\frac{2}{3}\right\} $ and $R=3$, where the
measure $\mu _{3}$ has $d_{H}=\frac{\ln 2}{\ln 3}$ (see Figure \ref{cantor2}%
), and is given by 
\begin{equation}
\widehat{\mu _{3}}\left( t\right) =\prod_{n=1}^{\infty }\frac{1}{2}\left(
1+e^{i\frac{4\pi t}{3^{n}}}\right) =e^{i\pi t}\prod_{n=1}^{\infty }\cos
\left( \frac{2\pi t}{3^{n}}\right)   \label{eqj}
\end{equation}
with 
\begin{equation}
\mathbf{Z}\left( \widehat{\mu _{3}}\right) =\left\{ \frac{3^{n}}{4}\left( 1+2%
\mathbb{Z}\right) :n=1,2,3,%
\dots%
%
\right\} .  \label{eqk}
\end{equation}
The contrast to the earlier case is that the sets $L$ which make (\ref{eq12}%
) valid, e.g., $L=\left\{ 0,\frac{3}{4}\right\} $, are \emph{not} contained
in $\mathbb{Z}$. It can be checked (see Theorem \ref{ThmNew6.1} below) that $%
L^{2}\left( \mu _{3}\right) $ does \emph{not} contain orthogonal families $%
\left\{ e_{\lambda }:\lambda \in P\right\} $ when $P$ ($\subset \mathbb{R}$)
has three, or more, distinct elements. In particular, there are no infinite
orthogonal sets of exponentials in $L^{2}\left( \mu _{3}\right) $. It can be
checked, from our formula for $\mathbf{Z}\left( \widehat{\mu _{4}}\right) $,
that $P_{4}$ ($=\{0,1,4,5,16,17,%
\dots%
%
\}$) is a unique maximal, orthogonal subset for $\mu _{4}$, containing $0$,
and contained in $\mathbb{N}_{0}$.

We finally mention that $\mu _{3,4}:=\mu _{3}\ast \mu _{4}$ (convolution)
has $P_{4}$ as a maximal orthogonal set for $L^{2}\left( \mu _{3,4}\right) $%
. (Recall $\mathbf{Z}\left( \widehat{\mu _{3}}\right) \cap \mathbf{Z}\left( 
\widehat{\mu _{4}}\right) =\varnothing $, from (\ref{eqd}) and (\ref{eqk})
above, and $\widehat{\mu _{3,4}}=\widehat{\mu _{3}}\cdot \widehat{\mu _{4}}$%
.) But $\left\{ e_{\lambda }:\lambda \in P_{4}\right\} $ is \emph{not} total
in $L^{2}\left( \mu _{3,4}\right) $; or, equivalently, 
\begin{equation*}
Q\left( t\right) :=\sum_{n\in P_{4}}\left| \widehat{\mu _{3,4}}\left(
t-n\right) \right| ^{2}
\end{equation*}
is \emph{not} constant on $\mathbb{R}$.

The possibilities described above may be summarized in the following two
theorems ($\nu =1$) for measures on $\mathbb{R}$ with affine
self-similarity. For both results, we have $R,b$ determined as follows: $%
R\in \mathbb{Z}\diagdown \left\{ 0\right\} $, $b\in \mathbb{R}\diagdown
\left\{ 0\right\} $; defining $\sigma _{0}\left( x\right) =R^{-1}x$, $\sigma
_{b}\left( x\right) =R^{-1}x+b$, recall that the corresponding probability
measure $\mu $ on $\mathbb{R}$ depends on both $R$ and $b$, being determined
uniquely by 
\begin{equation}
\label{determinemu}
\mu =\frac{1}{2}\left( \mu \circ \sigma _{0}^{-1}+\mu \circ \sigma
_{b}^{-1}\right)
\end{equation}
when $\left| R\right| >1$.

\begin{theorem}
\label{ThmNew6.1}Let $R$ be an odd integer, $R\neq \pm 1$, and let $\mu $ be
the corresponding measure. Then any set of $\mu $-orthogonal exponentials
contains at most two elements.
\end{theorem}

\begin{proof}%
%
Recall that 
\begin{equation*}
\hat{\mu}\left( t\right) =\prod_{n=0}^{\infty }\frac{1}{2}\left( 1+e^{i\frac{%
2\pi bt}{R^{n}}}\right) =e^{i\frac{\pi btR}{R-1}}\prod_{n=0}^{\infty }\cos
\left( \frac{\pi bt}{R^{n}}\right) 
\end{equation*}
and therefore 
\begin{equation*}
\mathbf{Z}\left( \hat{\mu}\right) =\left\{ \frac{R^{n}}{2b}\left( 2\mathbb{Z}%
+1\right) :n=0,1,%
\dots%
%
\right\} .
\end{equation*}
If $\gamma _{j}$, $j=1,2,3$, are such that the $e_{\gamma _{j}}$'s are
mutually orthogonal in $L^{2}\left( \mu \right) $, then the differences $%
\gamma _{i}-\gamma _{j}$ ($i\neq j$) are in $\mathbf{Z}\left( \hat{\mu}%
\right) $. Let $\lambda _{1}=\gamma _{1}-\gamma _{2}$, $\lambda _{2}=\gamma
_{2}-\gamma _{3}$, $\lambda _{0}=\gamma _{1}-\gamma _{3}$, and $\lambda _{j}=%
\frac{R^{n_{j}}}{2b}\left( 2z_{j}+1\right) $, $z_{j}\in \mathbb{Z}$. Since 
\begin{equation*}
\lambda _{1}+\lambda _{2}=\lambda _{0},
\end{equation*}
we get 
\begin{equation*}
R^{n_{1}}\left( 2z_{1}+1\right) +R^{n_{2}}\left( 2z_{2}+1\right)
=R^{n_{0}}\left( 2z_{0}+1\right) ,
\end{equation*}
which is a contradiction: the left-hand side is even while the right-hand
side is odd.%
\end{proof}%
%

\begin{theorem}
\label{ThmNew6.something}Let $R$ be an even \textup{(}$\neq 0$\textup{)}
integer, $\left| R\right| \geq 4$, and let $b\in \frac{1}{2}\mathbb{Z}%
\diagdown \left\{ 0\right\} $. Let $z_{0}=\frac{1}{2b}$, and set $%
P:=z_{0}\left\{ l_{0}+Rl_{1}+R^{2}l_{2}+%
\dots%
%
:l_{i}\in \left\{ 0,1\right\} \text{, finite sums}\right\} $. Then $\left\{
e_{\lambda }:\lambda \in P\right\} $ is an orthonormal basis for $%
L^{2}\left( \mu \right) $.
\end{theorem}

\begin{proof}%
%
The only part of the proof which is not included in the previous discussion
is the strict contractivity of the operator $C\colon Q\rightarrow C\left(
Q\right) $, from (\ref{eq21}). It specializes to 
\begin{equation*}
C\left( Q\right) \left( t\right) =\cos ^{2}\left( \pi bt\right) Q\left( 
\frac{t}{R}\right) +\sin ^{2}\left( \pi bt\right) Q\left( \frac{t-z_{0}}{R}%
\right) .
\end{equation*}
Let $J$ denote the closed interval between $0$ and $\frac{z_{0}}{1-R}$.
Depending on the signs, we can have $0$ as an endpoint to the left or the
right. Set 
\begin{equation*}
W_{0}=\left\{ Q\in C^{1}\left( J\right) :Q\left( 0\right) =1\right\}
\end{equation*}
and 
\begin{equation*}
d\left( Q_{1},Q_{2}\right) :=\sup_{t\in J}\left| Q_{1}^{\prime }\left(
t\right) -Q_{2}^{\prime }\left( t\right) \right| .
\end{equation*}
Then the argument from the proof of Theorem \ref{Thm3.4} leads to the
following strict contractivity for $Q\mapsto C\left( Q\right) $ with an
explicit formula for the contractivity constant ($\gamma <1$): 
\begin{equation*}
d\left( C\left( Q_{1}\right) ,C\left( Q_{2}\right) \right) \leq \gamma
d\left( Q_{1},Q_{2}\right) \quad \forall Q_{1},Q_{2}\in W_{0},
\end{equation*}
where 
\begin{equation*}
\gamma =\gamma \left( R\right) :=\frac{\pi }{2R}\sin \left( \frac{\pi }{%
\left| R\right| -1}\right) +\frac{1}{R}.
\end{equation*}
(The details of the proof of this contractivity property are the same, \emph{%
mutatis mutandis}, as the corresponding ones in Section \ref{S4} for the
special case $R=4$, $b=\frac{1}{2}$.) Specifically, one of the terms in the
estimation of $d\left( C\left( Q_{1}\right) ,C\left( Q_{2}\right) \right) $
will be 
\begin{equation*}
\pi b\left| \sin \left( 2\pi bt\right) \right| \int_{\frac{t-z_{0}}{R}}^{%
\frac{t}{R}}\left| Q_{1}^{\prime }\left( s\right) -Q_{2}^{\prime }\left(
s\right) \right| \,ds
\end{equation*}
for the case when $\frac{t-z_{0}}{R}\leq \frac{t}{R}$, and $t\in J$. The $%
\sin $-term may be estimated by $\left| \sin \left( 2\pi b\frac{z_{0}}{%
\left| R\right| -1}\right) \right| =\sin \left( \frac{\pi }{\left| R\right|
-1}\right) $ and the integral by $\left| \frac{z_{0}}{R}\right| d\left(
Q_{1},Q_{2}\right) $, and the result follows.

Note that $\gamma =\gamma \left( R\right) $ is a function only of $R$, i.e.,
independent of the second constant $b$. Also note that $R\mapsto \left|
\gamma \left( R\right) \right| $ depends only on $\left| R\right| $. In the
range of $\left| R\right| $, i.e., $\left| R\right| \geq 4$, $\left| \gamma
\left( R\right) \right| $ is decreasing, starting at $\frac{\pi }{16}\sqrt{3}%
+\frac{1}{4}$ ($<1$) when $\left| R\right| =4$. Finally, note that the
interval $J$ is chosen to be invariant under the two transformations, $%
t\mapsto \frac{t}{R}$ and $t\mapsto \frac{t-z_{0}}{R}$. This completes the
proof.%
\end{proof}%
%

\section{\label{S6}Three Dimensions}

Let $B$ and $L$ be two finite subsets in $\mathbb{R}^{\nu }$, both
containing $0$, having the same number of elements $N$, i.e., $\#\left(
B\right) =\#\left( L\right) =N$; and assume further that the matrix $H_{BL}$
in (\ref{eq12}) is unitary. Also suppose $L\subset \mathbb{Z}^{\nu }$. A $%
\nu $-by-$\nu $ matrix $R$ will be given satisfying the following four
conditions:

\begin{enumerate}
\item  \label{S6Intro(1)}$R=\left( a_{ij}\right) _{i,j=1}^{\nu }$, $%
a_{ij}\in \mathbb{Z}$,

\item  \label{S6Intro(2)}$R\left( B\right) \subset \mathbb{Z}^{\nu }$,

\item  \label{S6Intro(3)}the eigenvalues $\xi _{i}$ of $R$ satisfy $\left|
\xi _{i}\right| >1$,

\item  \label{S6Intro(4)}$N<\left| \det R\right| $.
\end{enumerate}

Define 
\begin{align}
P=P\left( L\right) & =\left\{ l_{0}+R^{\ast }l_{1}+R^{\ast \,2}l_{2}+
\dots
:l_{i}\in L\text{, finite sums}\right\} , \label{eqNew7.1}  \\
\sigma _{b}\left( x\right) & =R^{-1}x+b, \notag  \\
\tau _{l}^{{}}\left( x\right) & =R^{\ast }x+l, \notag  \\
\rho _{l}^{{}}& =\tau _{l}^{-1}=
\sigma _{l}^{\ast },\text{\quad and} \notag  \\
\chi _{B}\left( t\right) & =N^{-1}\sum_{b\in B}e^{i2\pi b\cdot t}.\notag  
\end{align}%
%
If $\mu $ is the measure (depending on $B$ and $R$) satisfying (\ref{eq6}),
then 
\begin{equation*}
\hat{\mu}\left( t\right) =\chi _{B}\left( t\right) \hat{\mu}\left(
R^{*\,-1}t\right) ,\quad t\in \mathbb{R}^{\nu }.
\end{equation*}
Clearly then 
\begin{equation*}
P=\bigcup_{l\in L}\left( l+R^{*}P\right) ,
\end{equation*}
and it follows from the unitarity of $H_{BL}$ and (\ref{S6Intro(1)})--(\ref
{S6Intro(2)}) that this union is disjoint. This is equivalent to the
assertion that the representation 
\begin{equation*}
\lambda =l_{0}+R^{*}l_{1}+R^{*\,2}l_{2}+\cdots \quad \left( l_{i}\in L\right)
\end{equation*}
of points in $P\left( L\right) $ is \emph{unique}. Then we proved in Lemmas 
\ref{Lem4.1} and \ref{Lem4.3} that 
\begin{equation*}
Q\left( t\right) :=\sum_{\lambda \in P}\left| \hat{\mu}\left( t-\lambda
\right) \right| ^{2}
\end{equation*}
is entire analytic, and that it satisfies (\ref{eq20}).

The attractor $\mathbf{X}\left( L\right) $ of $L$ is the (unique) compact
subset $\mathbf{X}$ of $\mathbb{R}^{\nu }$ satisfying 
\begin{equation*}
R^{*}\mathbf{X=}\bigcup_{l\in L}\left( \mathbf{X}-l\right) .
\end{equation*}

To prove that 
$P\left( L\right) $,
the ``fractal in the large'',
which is constructed from
some given
scaling
system $\left( R,B,L\right) $ in $\mathbb{R}^{\nu }$,
provides an orthogonal basis of
exponentials
$\left\{ e_{\lambda }:\lambda \in P\left( L\right) \right\} $
for
$L^{2}\left( \mu \right) $, using the Ruelle transfer
operator, we must exclude the
special case when $L$ is contained in
a hyperplane in $\mathbb{R}^{\nu }$. But this
condition is in fact necessary as we
show in Section \ref{S9} below. Moreover,
this restriction on $L$, together with the
Hadamard matrix property (\ref{eq12}) for
the two sets $B$ and $L$, is
stringent, and we showed in \cite{JoPe96}
that it can hold only in
dimensions $\nu =1,3$, and higher.
While we can get examples for
$\nu =2$ (the plane) where the
axioms, including (\ref{eq12}), are
satisfied, we showed in
\cite{JoPe96} that those planar examples will
``collapse'' down into a line: The
possible examples in the plane
will either be essentially
one-dimensional, or they will
not have $H_{2}\left( P\left( L\right) \right) $ equal
to $L^{2}\left( \mu \right) $, but only a proper
subspace; see Section \ref{S9} for more details.

The following scholium for $\mathbb{R}^{\nu }$ will serve as a guiding
principle for the results in Sections \ref{S8} and  \ref{S9}
below. We will also there
discuss two metrics on the testing functions $Q$, both involving the
gradient $\nabla Q$, but one giving a better contractivity constant than the
other. The proofs, and the evaluation of contractivity constants, will be
postponed to Sections \ref{S8} and \ref{S9}.

\begin{scholium}
\label{Lem6.1}The operator 
\begin{equation*}
C\left( Q\right) \left( t\right) :=\sum_{l\in L}\left| \chi _{B}\left(
t-l\right) \right| ^{2}Q\left( \rho _{l}\left( t\right) \right) 
\end{equation*}
acts boundedly on functions defined on the convex hull of $\mathbf{X}\left(
L\right) $.
\end{scholium}

\begin{proof}%
%
Details to follow: see the proof of Theorem \ref{Thm8.3} below,
and Section  \ref{S9}.%
\end{proof}%
%

The contractivity property of the scholium refers to the $C^{1}$ metric on
the convex hull of $\mathbf{X}\left( L\right) $ applied to the convex set of 
$C^{1}$ functions $Q$ satisfying $Q\left( 0\right) =1$, and the proof
requires all the listed assumptions, including (\ref{S6Intro(1)})--(\ref
{S6Intro(4)}) for the matrix $R$.

But with weaker assumptions, we can still get orthogonality of $\left\{
e_{\lambda }:\lambda \in P\left( L\right) \right\} $ and uniqueness for the $%
L$-expansion:

\begin{lemma}
\label{Lem6.2}With the following assumptions: 
\begin{align}
L& \subset \mathbb{Z}^{\nu },  \label{eqalpha} \\
R\left( B\right) & \subset \mathbb{Z}^{\nu },  \label{eqbeta} \\
H_{BL}^{{}}H_{BL}^{\ast }& =I_{N}^{{}},  \label{eqgamma}
\end{align}%
%
we have that the representation $\lambda =l_{0}+R^{\ast }l_{1}+R^{\ast
\,2}l_{2}+\cdots $ for points in $P\left( L\right) $ is unique, and the
exponentials $\left\{ e_{\lambda }:\lambda \in P\left( L\right) \right\} $
are mutually orthogonal in $L^{2}\left( \mu \right) $ where $\mu $ is the
probability measure on $\mathbb{R}^{\nu }$ subject to 
\begin{equation*}
\mu =N^{-1}\sum_{b\in B}\mu \circ \sigma _{b}^{-1}.
\end{equation*}
\end{lemma}

\begin{proof}%
%
The crucial step in the proof of both conclusions of the lemma is the
observation that 
\begin{equation*}
\chi _{B}\left( l-l^{\prime }+R^{\ast }\left( p-p^{\prime }\right) \right) =0
\end{equation*}
whenever $l$ and $l^{\prime }$ are \emph{distinct} points in $L$, and $%
p,p^{\prime }\in P\left( L\right) $. First notice that for $b\in B$: 
\begin{align*}
e_{b}\left( l-l^{\prime }+R^{\ast }\left( p-p^{\prime }\right) \right) &
=e_{b}\left( l-l^{\prime }\right) e_{b}\left( R^{\ast }\left( p-p^{\prime
}\right) \right) \\
& =e_{b}\left( l-l^{\prime }\right) e_{Rb}\left( p-p^{\prime }\right) \\
& =e_{b}\left( l-l^{\prime }\right) .
\end{align*}%
%
We have $Rb\in \mathbb{Z}^{\nu }$ by (\ref{eqbeta}), $p-p^{\prime }\in 
\mathbb{Z}^{\nu }$, in fact, $P\left( L\right) \subset \mathbb{Z}^{\nu }$,
by (\ref{eqalpha}), and $\sum_{b\in B}e_{b}\left( l-l^{\prime }\right) =0$
by (\ref{eqgamma}).

To prove the uniqueness assertion, it is enough to show that, if 
\begin{equation*}
l-l^{\prime }+R^{\ast }\left( p-p^{\prime }\right) =0
\end{equation*}
for some $l,l^{\prime }\in L$ and $p,p^{\prime }\in P\left( L\right) $, then 
$l=l^{\prime }$. (The uniqueness will then follow by induction.) Suppose
not, i.e., $l\neq l^{\prime }$; then we saw that 
\begin{equation*}
\chi _{B}\left( l-l^{\prime }+R^{\ast }\left( p-p^{\prime }\right) \right)
=0.
\end{equation*}
Since $\chi _{B}\left( 0\right) =1$, we have a contradiction.

The orthogonality of $e_{\lambda }$ and $e_{\lambda ^{\prime }}$ for
distinct $\lambda $ and $\lambda ^{\prime }$ in $P\left( L\right) $ follows
from the representation 
\begin{equation*}
\hat{\mu}\left( t\right) =\prod_{n=0}^{\infty }\chi _{B}\left( R^{\ast
\,-n}t\right)
\end{equation*}
applied to $t=\lambda -\lambda ^{\prime }$. The first nonzero term in the $L$%
-expansion of $\lambda -\lambda ^{\prime }$ yields a factor in the product
of the form 
\begin{equation*}
\chi _{B}\left( l-l^{\prime }+R^{\ast }\left( p-p^{\prime }\right) \right)
\end{equation*}
and we noted that this vanishes whenever $l\neq l^{\prime }$ in $L$. It
follows that $%
\ip{e_{\lambda ^{\prime }}}{e_{\lambda }}%
%
_{\mu }=\hat{\mu}\left( \lambda -\lambda ^{\prime }\right) =0$ which is the
desired orthogonality.%
\end{proof}%
%

Consider subsets $B,L\subset \mathbb{R}^{\nu }$ and an integer matrix $R$
with the properties mentioned in Lemma \ref{Lem3.1}. We then have dual
transforms 
\begin{equation*}
\sigma _{l}^{\ast }\left( t\right) :=R^{\ast \,-1}\left( t-l\right) ,\quad
t\in \mathbb{R}^{\nu }
\end{equation*}
indexed by $l\in L$. (Compare to the associated transformations $\left\{
\sigma _{b}:b\in B\right\} $ in (\ref{eq4}).) Assume that $\left\{ \sigma
_{l}^{\ast }:l\in L\right\} $ satisfies the ``open-set condition'' (\ref{eq5}%
) for some open set $V^{\ast }\subset \mathbb{R}^{\nu }$ which contains $0$
in its closure. 
We now turn to basis properties of the exponentials
$\left\{e_{\lambda }:\lambda \in P\right\} $
when $P=\left\{ l_{0}+R^{\ast }l_{1}+%
\dots%
%
:l_{i}\in L\right\} $. 
Note that Lemmas 
\ref{Lem3.1} and \ref{Lem3.3} will apply provided we can get 
\begin{equation*}
Q\longmapsto C\left( Q\right) \left( t\right) :=\sum_{l\in L}\left| \chi
_{B}\left( t-l\right) \right| ^{2}Q\left( \sigma _{l}^{\ast }\left( t\right)
\right)
\end{equation*}
contractive in the convex set 
\begin{equation}
W_{0}:=\left\{ Q\in C^{1}\left( \overline{V^{\ast }}\right) :Q\left(
0\right) =1\right\} .  \label{eq27}
\end{equation}
Note that some $\overline{V^{\ast }}$ may possibly be scaled down to get $%
C\colon V_{0}\rightarrow V_{0}$ strictly contractive.
This will be studied in detail in Sections \ref{S8} and  \ref{S9}.


\begin{figure}[tbp]
\begin{center}
\textit{Aus dem Paradies, das Cantor uns geschaffen, 
soll uns niemand \\vertrieben k\"onnen.} 
--- D. Hilbert \cite[p.\ 170]{Hil26}\\[\baselineskip]
\ \psfig{file=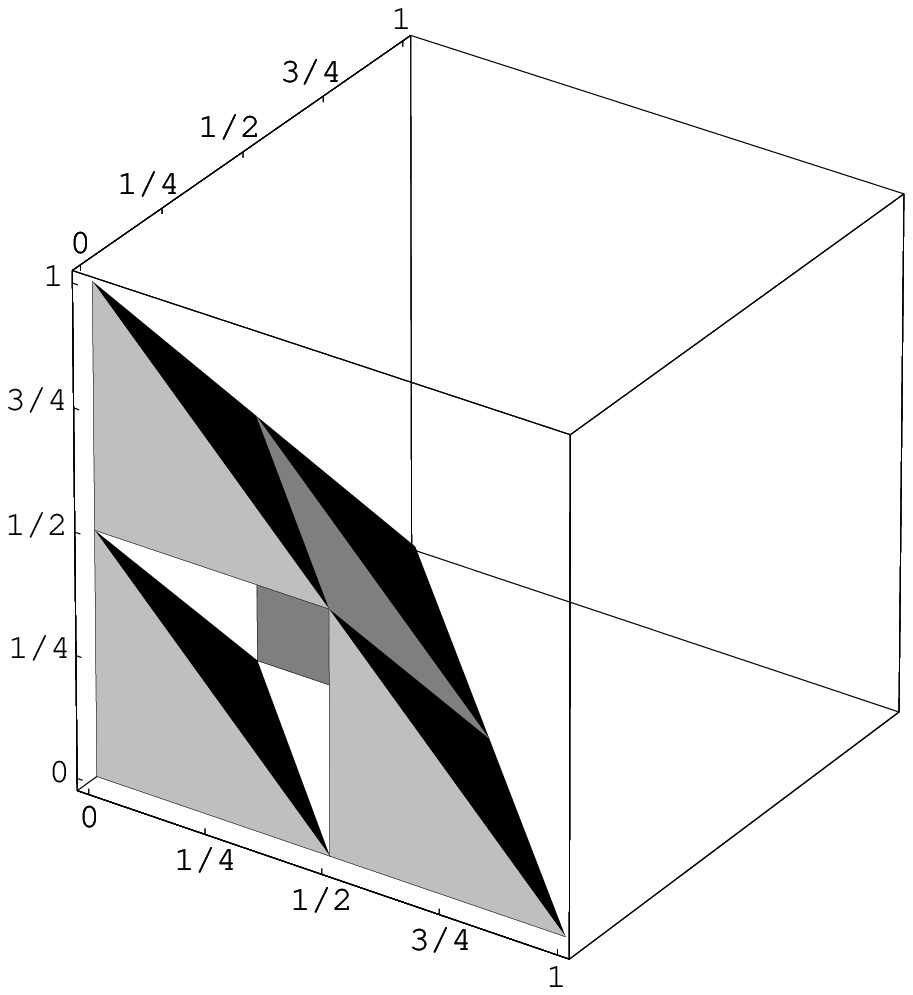,width=3in}\ \\
First iteration\\[\baselineskip]
\ \psfig{file=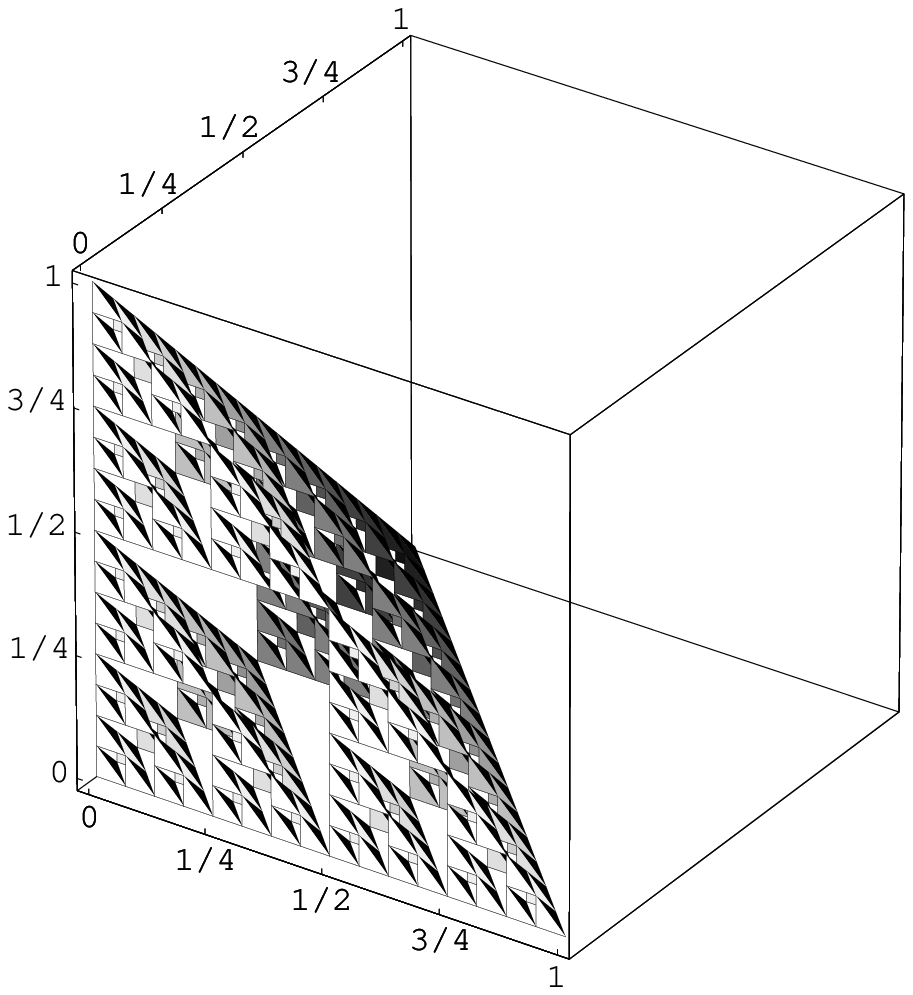,width=3in}\ \\
Fourth iteration (shading increases with depth)
\end{center}
\caption{The Eiffel Tower (Example \ref{ExaNew7.3})}
\label{eiffel}
\end{figure}%
%

\begin{example}
\label{ExaNew7.3}The conditions are satisfied in the following example (%
\cite[Example 7.4]{JoPe96}, the Eiffel Tower, see Figure \ref{eiffel}): 
\begin{equation}\label{eq7pound1}\begin{aligned}
\nu & =3,   \\
R& =\begin{pmatrix}2&0&0\\0&2&0\\0&0&2\end{pmatrix},  \\
B& =\left\{ \begin{pmatrix}0\\0\\0\end{pmatrix},
\begin{pmatrix}\frac12\\0\\0\end{pmatrix},
\begin{pmatrix}0\\\frac12\\0\end{pmatrix},
\begin{pmatrix}0\\0\\\frac12\end{pmatrix}\right\} ,   \\
L& =\left\{ \begin{pmatrix}0\\0\\0\end{pmatrix},
\begin{pmatrix}1\\1\\0\end{pmatrix},
\begin{pmatrix}1\\0\\1\end{pmatrix},
\begin{pmatrix}0\\1\\1\end{pmatrix}\right\} . 
\end{aligned}\end{equation}%
%
The
natural candidate for a
subset $P\subset \mathbb{R}^{3}$ such that $\left\{
e_{\lambda }:\lambda \in P\right\} $ is an orthonormal basis in $L^{2}\left(
\mu \right) $ is $P=\left\{ l_{0}+2l_{1}+2^{2}l_{2}+%
\dots%
%
:l_{i}\in L\text{, finite sums}\right\} $. If $\lambda \in P$, the three
coordinates $\lambda =\left( 
\begin{smallmatrix}a\\b\\c\end{smallmatrix}%
%
\right) $ are all in $\mathbb{N}_{0}$. 
One of the metrics on the corresponding space 
$W_{0}$ will be 
\begin{equation*}
d\left( Q_{1},Q_{2}\right) :=\sum_{j=1}^{\nu }\sup_{\overline{V^{\ast }}%
\mathstrut }\left| \frac{\partial \,}{\partial t_{j}}\left(
Q_{1}-Q_{2}\right) \right| .
\end{equation*}
\end{example}

The contractivity proof generalizes that of $\nu =1$ in Section \ref{S5}.
Note (in the general case $\nu >1$) 
\begin{equation}
\sum_{l\in L}\frac{\partial \,}{\partial t_{j}}\left| \chi _{B}\left(
t-l\right) \right| ^{2}=0,  \label{eq28}
\end{equation}
since 
\begin{equation}
\sum_{l\in L}\left| \chi _{B}\left( t-l\right) \right| ^{2}=1,\quad t\in 
\mathbb{R}^{\nu }.  \label{eq29}
\end{equation}
Pick a term in (\ref{eq28}) which is positive, e.g., 
\begin{equation*}
\frac{\partial \,}{\partial t_{1}}\left| \chi _{B}\left( t\right) \right|
^{2}\geq 0;
\end{equation*}
then 
\begin{equation*}
\sum_{l\in L\diagdown \left\{ 0\right\} }\frac{\partial \,}{\partial t_{1}}%
\left| \chi _{B}\left( t-l\right) \right| ^{2}=-\frac{\partial \,}{\partial
t_{1}}\left| \chi _{B}\left( t\right) \right| ^{2}.
\end{equation*}
(The significance of (\ref{eq28}) and (\ref{eq29}) will become
clear in Lemma \ref{LemNew7.5} below.)
Make line integrals connecting $0$ to $R^{*\,-1}\left( t-l\right) $ for each 
$l\in L\diagdown \left\{ 0\right\} $, or even curve integrals (if $\overline{%
V^{*}}$ isn't convex). Hence, replacing $Q_{1}-Q_{2}$ by $Q$, we find that 
\begin{equation}
\label{eq7pound2}
\sum_{l\in L}\frac{\partial \,}{\partial t_{1}}\left| \chi _{B}\left(
t-l\right) \right| ^{2}Q\left( \sigma _{l}^{*}\left( t\right) \right)
\end{equation}
is estimated by a constant times 
\begin{equation*}
\sum_{j=1}^{\nu }\left\| \frac{\partial \,}{\partial t_{j}}Q\right\|
_{\infty }=\left\| \nabla Q\right\| _{\infty \overline{V^{*}}}.
\end{equation*}

The property (\ref{eq29}) is important for the argument, and we include the
(easy) proof. It is based on the assumption (\ref{eq12}) placed on the two
given finite subsets $B$ and $L$ in $\mathbb{R}^{\nu }$. We now calculate
from the left-hand side of (\ref{eq29}): 
\begin{align*}
\sum_{l\in L}\left| \chi _{B}\left( t-l\right) \right| ^{2}&
=N^{-2}\sum_{l\in L}\left| 
\vphantom{\sum}\smash{
\sum_{b\in B}}e_{t-l}\left( b\right) \right| ^{2}
\\
& =N^{-2}\sum_{l\in L}\sum_{b\in B}\sum_{b^{\prime }\in B}e_{t-l}\left(
b-b^{\prime }\right)  \\
& =N^{-2}\sum_{b\in B}\sum_{b^{\prime }\in B}e_{t}\left( b-b^{\prime
}\right) \sum_{l\in L}e_{l}\left( b^{\prime }-b\right)  \\
& =N^{-1}\sum_{b\in B}\sum_{b^{\prime }\in B}e_{t}\left( b-b^{\prime
}\right) \delta _{b,b^{\prime }} \\
& =N^{-1}\sum_{b\in B}e_{t}\left( 0\right) =N^{-1}N=1,
\end{align*}
which is the claim in (\ref{eq29}).

\begin{remark}
\label{Rem6.something}For the Eiffel Tower example we may take 
\begin{equation}
\overline{V^{*}}=\left\{ \left( t_{1},t_{2},t_{3}\right) :-a\leq t_{i}\leq
0\right\}  \label{eq6.something.1}
\end{equation}
where $a\in \mathbb{R}_{+}$ is determined to get the contractive property on 
$Q\mapsto C\left( Q\right) \left( t\right) =\sum_{l\in L}\left| \chi
_{B}\left( t-l\right) \right| ^{2}Q\left( \sigma _{l}^{*}\left( t\right)
\right) $. Set 
\begin{equation*}
d\left( Q_{1},Q_{2}\right) =\sum_{i=1}^{3}\left| \frac{\partial \,}{\partial
t_{i}}\left( Q_{1}-Q_{2}\right) \right| _{\infty \overline{V_{a}^{*}}},
\end{equation*}
and 
\begin{equation*}
W_{0}\left( a\right) =\left\{ Q\in C^{1}\left( -a\leq t_{i}\leq 0\right)
:Q\left( 0\right) =1\right\} .
\end{equation*}
Then 
\begin{equation*}
\frac{\partial \,}{\partial t_{j}}C\left( Q\right) \left( t\right) =\frac{1}{%
2}\sum_{l\in L}\left| \chi _{B}\left( t-l\right) \right| ^{2}\frac{\partial
\,}{\partial t_{j}}Q\left( \sigma _{l}^{*}\left( t\right) \right) +\text{
second term,}
\end{equation*}
and therefore 
\begin{equation*}
d\left( C\left( Q_{1}\right) ,C\left( Q_{2}\right) \right) \leq \frac{1}{2}%
d\left( Q_{1},Q_{2}\right) +\text{ second term.}
\end{equation*}
The ``second term'' comes from evaluating $\frac{\partial \,}{\partial t_{j}}%
\left| \chi _{B}\left( t-l\right) \right| ^{2}$ when 
\begin{equation}
\label{eq7pound3}
\chi _{B}\left( t-l\right) =\frac{1}{4}
\vphantom{1\underset{\text{two minus
signs}}{\underbrace{\pm e^{i\pi t_{1}}\pm e^{i\pi t_{2}}\pm e^{i\pi t_{3}}}}}
\left( \vphantom{1\pm e^{i\pi t_{1}}\pm e^{i\pi t_{2}}\pm e^{i\pi t_{3}}}
\smash{1\underset{\text{two minus
signs}}{\underbrace{\pm e^{i\pi t_{1}}\pm e^{i\pi t_{2}}\pm e^{i\pi t_{3}}}}}
\right) .
\end{equation}
We use the estimate 
\begin{equation*}
\left| \frac{\pi }{4}
\vphantom{\sum}\smash{
\sum_{j=1}^{3}}\int_{-a}^{0}\frac{\partial \,}{\partial
t_{j}}Q\right| \leq \frac{\pi }{4}a\sum_{j=1}^{3}\left| \frac{\partial \,}{%
\partial t_{j}}Q\right| _{\infty }
\end{equation*}
and complete the proof of the contractive property of the operator $Q\mapsto
C\left( Q\right) $: We get a strict contraction if $a\in \mathbb{R}_{+}$ is
chosen (so small that) $\frac{1}{2}+\frac{a\pi }{4}<1$, i.e., $a<\frac{2}{%
\pi }$. Note however that the set (\ref{eq6.something.1}) will not be
invariant under the affine dual maps $t\mapsto R^{-1}\left( t-l\right) $, $%
l\in L$, unless $R$ is expanded by a certain scale, to be specified in the
following section. We show there that strict contractivity of $Q\mapsto
C\left( Q\right) $ may always be obtained by such a ``rescaling''. 
(See the details below in Proposition \ref{ProNew7.6}.)
We also
give a different metric on the functions $Q$ which yields sharper results.
Proposition \ref{ProNew7.6}
implies that the contractivity is indeed strict for the
matrix $R=\left( 
\begin{smallmatrix}
r & 0 & 0 \\
0 & r & 0 \\
0 & 0 & r
\end{smallmatrix}
\right) $, $r\geq 3$,
of Example \ref{ExaNew7.3}. We then get $Q\equiv 1$ by the Banach
fixed point principle. Therefore $\frac{\partial ^{m}\;}{\partial
t_{1}^{m_{1}\mathstrut }\partial t_{2}^{m_{2}\mathstrut }\partial
t_{3}^{m_{3}\mathstrut }}Q\left( t\right) |_{t=0}=0$ if $m=m_{1}+m_{2}+m_{3}%
\geq 1$. Hence $x_{1}^{m_{1}}x_{2}^{m_{2}}x_{3}^{m_{3}}\in H_{2}\left(
P\left( L\right) ,\mu \right) $ and $H_{2}\left( P\left( L\right) ,\mu
\right) =L^{2}\left( \mu \right) $ as claimed.
\end{remark}

For use both in the detailed estimates in Example \ref{ExaNew7.3}, and also
in Sections \ref{S8} and \ref{S9} below, we need the following lemma. Its
statement involves some notation: $\left\| \,\cdot \,\right\| _{2}$ for the $%
l^{2}$-norm in $\nu $ dimension, $\left\| \,\cdot \,\right\| _{op}$ for the
corresponding \emph{operator norm} on $\nu \times \nu $ matrices, i.e.,
\begin{equation*}
\left\| R\right\| _{op}:=\max \left\{ \left\| Rx\right\| _{2}:x\in %
\mathbb{R}^{\nu },\;\left\| x\right\| _{2}=1\right\} .
\end{equation*}
Finally, let $\left| \,\cdot \,\right| _{\infty ,\mathbf{Y}}$ be the \emph{%
supremum-norm} applied to functions on $\mathbf{Y}$, i.e., taking supremum
over some simplex $\mathbf{Y}\subset \mathbb{R}^{\nu }$.

\begin{lemma}
\label{LemNew7.5}Let the setting be as in Scholium \textup{\ref{Lem6.1}},
corresponding to a given expansive \textup{(}see \textup{(\ref{eq3}))}
scaling matrix $R$, and a finite set $L$
of translation vectors in $\mathbb{R}^{\nu }$, such that $0\in L$. Define $%
L^{\ast }=L\diagdown \left\{ 0\right\} $.

\begin{enumerate}
\item  \label{LemNew7.5(1)}Then, for every set of coefficients
$\left\{ c_{l}\right\} _{l\in L}$%
, $c_{l}\in \mathbb{C}$, such that $\sum_{l\in L}c_{l}=0$, we have
\begin{equation*}
S_{c}\left( t\right) =\sum_{l\in L}c_{l}Q\left( \rho _{l}\left( t\right)
\right) =\sum_{l\in L^{\ast }}c_{l}\left( Q\left( R^{-1}\left( t-l\right)
\right) -Q\left( R^{-1}t\right) \right) 
\end{equation*}
where $\rho _{l}\left( t\right) :=R^{-1}\left( t-l\right) $.

\item  \label{LemNew7.5(2)}Since
\begin{equation*}
Q\left( R^{-1}\left( t-l\right) \right) -Q\left( R^{-1}t\right)
=-\int_{0}^{1}\left( R^{-1}\left( \nabla _{l}\right) Q\right) \left(
R^{-1}\left( t-sl\right) \right) \,ds,
\end{equation*}
we get the estimate
\begin{equation*}
\left| S_{c}\left( t\right) \right| \leq \sum_{l\in L^{\ast }}\left|
c_{l}\right| \int_{0}^{1}\left| R^{-1}\left( \nabla _{l}\right) Q\left(
R^{-1}\left( t-sl\right) \right) \right| \,ds,
\end{equation*}
and, if $\mathbf{Y}$ is the convex hull of the attractor $\mathbf{X}$ 
\textup{(}i.e., $\mathbf{X}$ is the solution to $\mathbf{X}=L+R\mathbf{X}$%
\textup{)}, then, for all $t\in \mathbf{Y}$, 
\begin{equation*}
\left| S_{c}\left( t\right) \right| \leq \left\| R^{-1}\right\|
_{op}\sum_{l\in L^{\ast }}\left| c_{l}\right| \left\| l\right\| _{2}\left|
\left\| \nabla Q\right\| _{2}\right| _{\infty ,\mathbf{Y}}.
\end{equation*}

\item  \label{LemNew7.5(3)}In case
$R$ is diagonalizable with a single positive eigenvalue, then 
$\mathbf{Y}$ \textup{(}$=\limfunc{conv}\mathbf{X}$, i.e., the smallest
convex subset in $\mathbb{R}^{\nu }$ invariant under all the $\rho _{l}$'s%
\textup{)} is equal to the simplex $\mathbf{Y}^{\prime }$ in $\mathbb{R}%
^{\nu }$ generated by the points $y_{l}:=-\left( R-I\right) ^{-1}l$, $l\in
L^{\ast }$.
\end{enumerate}
\end{lemma}

\begin{proof}
In the formula for $S_{c}\left( t\right) $, substitute $c_{0}=-\sum_{l\in
L^{\ast }}c_{l}$. Then 
\begin{align*}
S_{c}\left( t\right) &=-\sum_{l\in L^{\ast
}}c_{l}Q\left( R^{-1}t\right) +\sum_{l\in L^{\ast }}c_{l}Q\left(
R^{-1}\left( t-l\right) \right) \\ 
&=\sum_{l\in L^{\ast }}c_{l}\left( Q\left(
R^{-1}\left( t-l\right) \right) -Q\left( R^{-1}t\right) \right) 
\end{align*}
as
claimed. The difference terms inside the sum may then be estimated using the
auxiliary functions $\varphi _{l}\left( s\right) :=Q\left( R^{-1}\left(
t-sl\right) \right) $, $0\leq s\leq 1$, where
\begin{equation*}
\frac{d\,}{ds}\varphi _{l}\left( s\right) =-R^{-1}\left( \nabla _{l}\right)
Q\left( R^{-1}\left( t-sl\right) \right) ,
\end{equation*}
and noting that, for each $s\in \left[ 0,1\right] $, $l\in L$, and for $t\in 
\mathbf{Y}$, the points $R^{-1}\left( t-sl\right) $ are all in $\mathbf{Y}$.

To prove the asserted invariance of the simplex $\mathbf{Y}^{\prime }$, let $%
\mathbf{Y}$ denote the convex hull (in $\mathbb{R}^{\nu }$) of the attractor 
$\mathbf{X}$ of the affine system $\left\{ \rho _{l}:l\in L\right\} $ (i.e.,
the unique $\mathbf{X}$ ($\subset \mathbb{R}^{\nu }$) satisfying $\mathbf{X}%
=L+R\mathbf{X}$). Note that, for each $l\in L$, $y_{l}=-\left( R-I\right)
^{-1}l$ is the unique fixed point of $\rho _{l}$, including the case $l=0$,
where $\rho _{0}\left( 0\right) =0$.

For $n\in \mathbb{N}$, and $\left( l_{1},l_{2},\dots ,l_{n}\right) \in 
\underset{n\text{ times}}{\underbrace{L\times L\times \dots \times L}}$, let 
$y\left( l_{1},\dots ,l_{n}\right) $ be the fixed point of $\rho
_{l_{1}}\circ \rho _{l_{2}}\circ \dots \circ \rho _{l_{n}}$, and let $%
\mathbf{X}_{n}$ be the set of all the fixed points $y\left( l_{1},\dots
,l_{n}\right) $ as $\left( l_{1},\dots ,l_{n}\right) $ varies over $\bigcross%
_{1}^{n}L$. Then it is known (see \cite[Theorem 1]{Hut81}) that
$\mathbf{X}
=\overline{\bigcup_{n\in \mathbb{N}}\mathbf{X}_{n}}$,
and therefore
\begin{equation}
\mathbf{Y}=\limfunc{conv}\mathbf{X}
=\overline{\bigcup_{n\in \mathbb{N}}\limfunc{conv}%
\mathbf{X}_{n}}  \label{eq7ins1}
\end{equation}
where $\limfunc{conv}$ denotes the convex hull. Since $\rho _{0}\left(
0\right) =0$, $\limfunc{conv}\mathbf{X}_{1}$ is the \emph{simplex} $\mathbf{Y%
}^{\prime }$ generated by $\left\{ y_{l}:l\in L^{\ast }\right\} $ where $%
y_{l}=-\left( R-1\right) ^{-1}l$, i.e.,
\begin{equation*}
\mathbf{Y}^{\prime }=\limfunc{conv}\mathbf{X}_{1}=
\vphantom{\sum_{l\in L^{\ast}}}
\left\{ \vphantom{\sum}\smash{\sum_{l\in L^{\ast}}}
s_{l}y_{l}:s_{l}\geq 0,\;\vphantom{\sum}\smash{\sum_{l\in L^{\ast}}}
s_{l}\leq 1\right\} .
\end{equation*}
Using now the contractivity property of the affine maps $\left\{ \rho
_{l}:l\in L\right\} $ implied by (\ref{eq3}), we infer by induction that, if 
$R$ is diagonalizable with a single positive eigenvalue, then
\begin{equation}
\mathbf{X}_{n}\subset \mathbf{Y}^{\prime }\text{\quad for all }n\in %
\mathbb{N},  \label{eq7ins2}
\end{equation}
and, therefore, by (\ref{eq7ins1}), that $\mathbf{Y}\subset \mathbf{Y}%
^{\prime }$. A further volume consideration shows that in fact $\mathbf{Y}=%
\mathbf{Y}^{\prime }$. The invariance of the simplex follows, as claimed.

To prove (\ref{eq7ins2}), proceed by induction, checking first that $\mathbf{%
X}_{2}\subset \limfunc{conv}\mathbf{X}_{1}$ when $R$ has the form $R=\left( %
\begin{smallmatrix}r&0&{\textstyle \cdots } &0 \\ 
0&r&{\textstyle \cdots } &0 \\ 
\raisebox{0pt}[10pt]{$\vdots$} 
&\raisebox{0pt}[10pt]{$\vdots$} 
&\raisebox{0pt}[10pt]{$\mkern-1.5mu\ddots\mkern-1.5mu$} 
&\raisebox{0pt}[10pt]{$\vdots$}\\ 
0&0&{\textstyle \cdots } &r \end{smallmatrix}
\right) $ for some $r\in \mathbb{N}$%
. Every $y\in \mathbf{X}_{2}$ has the representation $y=\left(
I-R^{2}\right) ^{-1}\left( l_{0}+Rl_{1}\right) $ for some $l_{0},l_{1}\in L$%
. It follows that $y=\left( R+I\right) ^{-1}x_{0}+R\left( R+I\right)
^{-1}x_{1}$ where $x_{i}:=-\left( R-I\right) ^{-1}l_{i}\in \mathbf{X}_{1}$
for $i=0,1$. So when the scaling matrix $R$ has the stated form, it follows
that $y=\frac{1}{r+1}x_{0}+\frac{r}{r+1}x_{1}\in \limfunc{conv}\mathbf{X}%
_{1}$. 
Hence $\mathbf{X}_{2}\subset \limfunc{conv}\mathbf{X}_{1}$.
Now suppose (\ref{eq7ins2}) has been proved for $n$, and let $y\in 
\mathbf{X}_{n+1}$. Then $y=\left( I-R^{n+1}\right) ^{-1}\left(
l_{0}+Rl_{1}+\dots +R^{n}l_{n}\right) $ for some $l_{i}\in L$, and it
follows that there are $z\in \mathbf{X}_{n}$ ($z=( I-R^{n})
^{-1}( l_{0}+Rl_{1}+\dots +R^{n-1}l_{n-1}) $), and $x\in \mathbf{X%
}_{1}$, such that $y=\frac{I-R^{n}}{I-R^{n+1\mathstrut }}z+\frac{R^{n}}{%
I+R+\dots +R^{n\mathstrut }}x$. Clearly $\frac{1-r^{n}}{1-r^{n+1\mathstrut }}%
+\frac{r^{n}}{1+r+\dots +r^{n\mathstrut }}=1$, so that, when $R$ has the
stated form, it follows from the inductive hypothesis that $y\in \limfunc{%
conv}\left( \mathbf{X}_{n}\cup \mathbf{X}_{1}\right) \subset \mathbf{Y}%
^{\prime }$, concluding the proof.
\end{proof}

Note that $\mathbf{Y}$ is not in general contained in $\limfunc{conv}\left( 
\mathbf{X}_{1}\right) $, even if $R$ is in diagonal form.
For examples when $\limfunc{conv}\left( \mathbf{X}\right) $
is of a more complex nature, see
\cite[Sections 9--10 and the graphics]{BrJo96b}.

\emph{The first application} is to the contractivity of the $C$-operator in
Example \ref{ExaNew7.3}: Let the $R$ in (\ref{eq7pound1}) be variable; i.e., 
$R=\left( \begin{smallmatrix}r&0&0\\0&r&0\\0&0&r\end{smallmatrix}\right) $,
rather than the $r=2$ version (\ref{eq7pound1}), and assume $r\in \mathbb{N}$%
, $r\geq 2$. Then, for each $r$, we have a corresponding simplex (as in the
lemma). It is given by parameters $s_{i}\geq 0$, $i=1,2,3$, such that $%
s_{1}+s_{2}+s_{3}\leq 1$. Specifically,
\begin{equation}
\textstyle \mathbf{Y}\left( r\right) =
\vphantom{\sum\limits_{i=1}^{3}}
\left\{ t=\left( 
\begin{smallmatrix}t_{1}\\t_{2}\\t_{3}
\end{smallmatrix}\right) 
\in \mathbb{R}^{3}:t_{1}=
\frac{s_{2}+s_{3\mathstrut }}{1-r},t_{2}=
\frac{s_{1}+s_{3\mathstrut }}{1-r},t_{3}=
\frac{s_{1}+s_{2\mathstrut }}{1-r},s_{i}\geq
0\;\forall i,
\vphantom{\sum}\smash{\sum\limits_{i=1}^{3}}
s_{i}\leq 1\right\} .  \label{eq7ins3}
\end{equation}
(See Figure \ref{simplex}.) Recall that $\mathbf{Y}:=\mathbf{Y}\left(
r\right) $ is (by the lemma) then invariant under the four affine
transformations,
\begin{equation*}
t\longmapsto \frac{t-l}{r}\quad (l\in L),
\end{equation*}
and also $\frac{t-sl}{r}\in \mathbf{Y}$, whenever $t\in \mathbf{Y}$ and $%
s\in \left[ 0,1\right] $.

Let $B$ and $L\subset \mathbb{R}^{3}$ be as
in (\ref{eq7pound1}) of Example \ref{ExaNew7.3}. For
$r\in \mathbb{N}$, let
\begin{align*}
P_{r}\left( L\right) &=  \rlap{$\displaystyle
\left\{ l_{0}+rl_{1}+r^{2}l_{2}+
\dots
:l_{i}\in L\text{, finite sums}\right\} ,$}  & &   \\
C_{r}\left( f\right) \left( t\right) &=\sum _{l\in L}
\left| \chi _{B}\left( t-l\right) \right| ^{2}
f\left( \frac{t-l}{r}\right) , 
& f &\in C^{\infty }\left( \mathbb{R}^{3}\right) ,\;t\in \mathbb{R}^{3}, \\
\intertext{and}
Q_{r}\left( t\right) &:=\sum _{\lambda \in P_{r}\left( L\right) }
\left| \widehat{\mu _{r}}\left( t-\lambda \right) \right| ^{2}, 
& t &\in \mathbb{R}^{3}.
\end{align*}
We then have the following 
conclusions concerning this example.
\begin{enumerate}
\item  \label{firstapplicationconclusion(1)}If $r$ is even, then
(\ref{eq12bis}) holds, and we then have the eigenvalue equation
$C_{r}\left( Q_{r}\right) =Q_{r}$ satisfied (by
Lemma \ref{Lem4.1}).

\item  \label{firstapplicationconclusion(2)}If
$r\geq 3$, then $C_{r}$
is contractive in the supremum-gradient 
norm (Proposition \ref{ProNew7.6} below). Hence,
if $r\geq 4$ and even, then 
the exponentials
\begin{equation*}
E_{r}\left( L\right) :=
\left\{ e_{\lambda }:\lambda \in P_{r}\left( L\right) \right\}
\end{equation*}
form an orthonormal basis 
in $L^{2}\left( \mu _{r}\right) $.

\item  \label{firstapplicationconclusion(3)}But 
(\ref{firstapplicationconclusion(1)})--(\ref{firstapplicationconclusion(2)})
do \emph{not} 
carry over to the case
when $r$ in $\mathbb{N}$ is odd. For
example, it can be checked
directly that $C_{3}\left( Q_{3}\right) \neq Q_{3}$,
and that $E_{3}\left( L\right) $ is \emph{not}
orthogonal in $L^{2}\left( \mu _{3}\right) $. (The
points $\left( 0,0,0\right) $ and $\left( 4,4,0\right) $ are
both in $P_{3}\left( L\right) $, but not
orthogonal in $L^{2}\left( \mu _{3}\right) $. Specifically,
$\widehat{\mu _{3}}\left( 4,4,0\right) =
\widehat{\mu _{3}}\left( \frac{4}{3},\frac{4}{3},0\right) 
\neq 0$ by a direct check.)
\end{enumerate}

\TeXButton{figure4}{
\begin{figure}[tbp]
\hspace*{12pt}\psfig{file=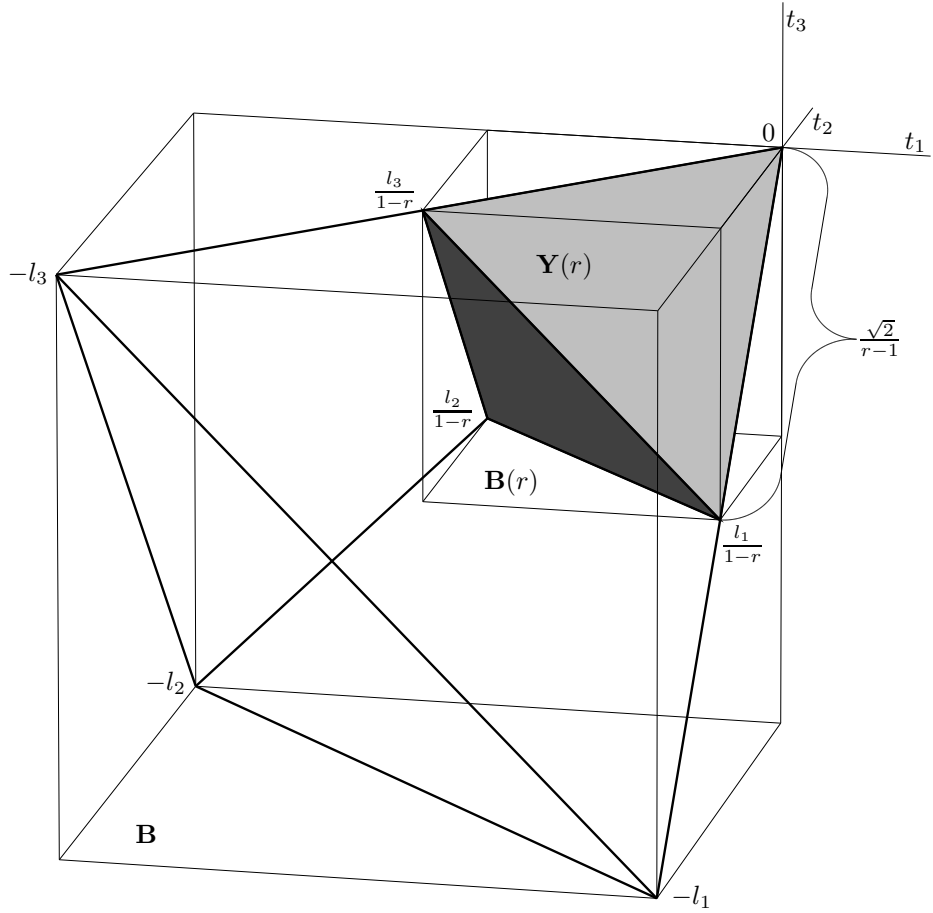,width=360pt}
\makebox[360pt]{%
\raisebox{3pt}[0pt][0pt]{%
\hspace*{8pt}\makebox[0pt]{\raisebox{256pt}[0pt][0pt]{$-l_{3}$}\hss}%
\hspace*{46pt}\makebox[0pt]{\raisebox{43pt}[0pt][0pt]{$\mathbf{B}$}\hss}%
\hspace*{6pt}\makebox[0pt]{\raisebox{101pt}[0pt][0pt]{$-l_{2}$}\hss}%
\hspace*{92pt}%
\makebox[0pt]{\hss\raisebox{287pt}[0pt][0pt]{$\frac{l_{3\mathstrut}}{1-r}$}}%
\hspace*{22pt}%
\makebox[0pt]{\hss\raisebox{204pt}[0pt][0pt]{$\frac{l_{2\mathstrut}}{1-r}$}}%
\hspace*{16pt}\makebox[0pt]{\raisebox{177pt}[0pt][0pt]{$\mathbf{B}(r)$}\hss}%
\hspace*{20pt}\makebox[0pt]{\raisebox{259pt}[0pt][0pt]{$\mathbf{Y}(r)$}\hss}%
\hspace*{49pt}\makebox[0pt]{\raisebox{20pt}[0pt][0pt]{$-l_{1}$}\hss}%
\hspace*{19pt}%
\makebox[0pt]{\raisebox{153pt}[0pt][0pt]{$\frac{l_{1\mathstrut}}{1-r}$}\hss}%
\hspace*{12pt}\makebox[0pt]{\raisebox{309pt}[0pt][0pt]{$0$}\hss}%
\hspace*{10pt}\makebox[0pt]{\raisebox{352pt}[0pt][0pt]{$t_{3}$}\hss}%
\hspace*{10pt}\makebox[0pt]{\raisebox{313pt}[0pt][0pt]{$t_{2}$}\hss}%
\hspace*{20pt}%
\makebox[0pt]{\raisebox{231pt}[0pt][0pt]{$\frac{\sqrt{2}}{r-1}$}\hss}%
\hspace*{15pt}\makebox[0pt]{\raisebox{306pt}[0pt][0pt]{$t_{1}$}\hss}%
\hspace*{15pt}}}
\caption{The Simplex $\protect\mathbf{Y}\left( r\right) $ and Its Vertices}
\label{simplex}
\end{figure}
}
\begin{proposition}
\label{ProNew7.6}In the Eiffel Tower example \textup{(}Example
\textup{\ref
{ExaNew7.3})}, the operator $Q\mapsto C_{r}\left( Q\right) $ is
contractive when
$r
\geq 3$ on the space of functions
$Q\in C^{\infty }\left( \mathbb{R}^{3}\right) $ such that
$Q\left( 0\right) =0$, and relative to the
norm
\begin{equation*}
\left\| Q\right\| _{1,\,\mathbf{Y}\left( 3\right) }:=\sup_{t\in \mathbf{Y}%
\left( 3\right) }\left\| \nabla Q\left( t\right) \right\| 
\end{equation*}
where $\left\| \nabla Q\right\| :=\sum_{j=1}^{3}\left| \frac{\partial \,}{%
\partial t_{j}}Q\right| $.
\end{proposition}

\begin{proof}
The first term in $C_{r}\left( Q\right) $ of Example \ref{ExaNew7.3} is (\ref
{eq7pound2}), and we will restrict, for simplicity, to functions $Q\in
C^{\infty }\left( \mathbb{R}^{3}\right) $ such that $Q\left( 0\right) =0$.
For the derivative terms $\frac{\partial \,}{\partial t_{1}}\left| \chi
_{B}\left( t-l\right) \right| ^{2}$ of (\ref{eq7pound2}), we have
\begin{align*}
\sum_{l\in L^{\ast }}\left| \frac{\partial \,}{\partial t_{1}}\left| \chi
_{B}\left( t-l\right) \right| ^{2}\right| 
&=\frac{\pi }{2}\sum_{l\in L^{\ast
}}\left| \func{Im}\chi _{B}\left( t-l\right) \right|  \\
&=\frac{\pi }{8}
\sum_{l\in L^{\ast }}
\vphantom{\sum_{k=1}^{3}}
\left| \vphantom{\sum}\smash{\sum_{k=1}^{3}}
\left( \pm \right) \sin \left( \pi
t_{k}\right) \right| 
\end{align*}
with the sign-convention (and vanishing of one of the
three terms in the sum) from (\ref{eq7pound3}) above.

Hence
\begin{equation}
\sum_{l\in L^{\ast }}\left| \frac{\partial \,}{\partial t_{1}}\left| \chi
_{B}\left( t-l\right) \right| ^{2}\right| \leq \frac{3\pi }{8}%
\sum_{k=1}^{3}\left| \sin \left( \pi t_{k}\right) \right| \leq \frac{3\pi }{4%
}\sin \left( \frac{\pi }{r-1}\right)   \label{eq7ins4}
\end{equation}
when $t\in \mathbf{Y}\left( r\right) $. (If $r=2$, we bound $\left| \sin
\left( \pi t_{k}\right) \right| $ by $1$.) The last estimate is from
evaluation of the maximum of the $\sum_{k=1}^{3}$-term when $t$ is
restricted to the simplex $\mathbf{Y}\left( r\right) $ from Figure \ref
{simplex}; see also (\ref{eq7ins3}). (The maximum is attained at one of the
vertices.)

Returning to the estimate on the sum 
(\ref{eq7pound2}), we have (using Lemma 
\ref{LemNew7.5}):
\begin{equation}
\left| \text{\textup{(\ref{eq7pound2})}}
\right| \leq \frac{3\pi }{4}\sin \left( \frac{%
\pi }{r-1}\right) \frac{1}{r}\sum_{l\in L^{\ast }}\int_{0}^{1}\left| \nabla
_{l}Q\left( \frac{t-sl}{r}\right) \right| \,ds  \label{eq7ins5}
\end{equation}
where estimate (\ref{eq7ins4}) was used. But when $t$ is restricted to $%
\mathbf{Y}\left( r\right) $, then $\frac{t-sl}{r}\in \mathbf{Y}\left(
r\right) $, and the term under the integral may in that case be estimated by
a supremum over $\mathbf{Y}\left( r\right) $. The corresponding supremum of $%
\sum_{j=1}^{3}\left| \frac{\partial \,}{\partial t_{j}}Q\right| $ will be
denoted simply $\left\| \nabla Q\right\| $, but a more specific terminology
would be $\left| \left\| \nabla Q\right\| _{1}\right| _{\infty ,\,\mathbf{Y}%
\left( r\right) }$. Continuing estimate (\ref{eq7ins5}), we get, for $t\in 
\mathbf{Y}\left( r\right) $, 
\begin{equation*}
\left| \text{\textup{(\ref{eq7pound2})}}
\right| \leq \frac{3\pi }{2}\sin \left( \frac{%
\pi }{r-1}\right) \frac{1}{r}\limfunc{vol}\nolimits_{3}\left( \mathbf{Y}%
\left( r\right) \right) \left\| \nabla Q\right\| .
\end{equation*}
Since $\limfunc{vol}\nolimits_{3}\left( \mathbf{Y}\left( r\right) \right) =%
\frac{1}{3\left( r-1\right) ^{3\mathstrut }}$, we then get
\begin{equation*}
\left| \text{\textup{(\ref{eq7pound2})}}
\right| \leq \frac{\pi }{2r\left( r-1\right)
^{3\mathstrut }}\sin \left( \frac{\pi }{r-1}\right) \left\| \nabla Q\right\|
.
\end{equation*}
The other term in the estimate of $\left| \frac{\partial \,}{\partial t_{1}}%
C_{r}\left( Q\right) \left( t\right) \right| $ is $\frac{1}{r}\left|
C_{r}\left( \frac{\partial \,}{\partial t_{1}}Q\right) \right| $, which, for 
$t\in \mathbf{Y}\left( r\right) $, is bounded above by $\frac{1}{r}\left| 
\frac{\partial \,}{\partial t_{1}}Q\right| _{\infty ,\,\mathbf{Y}\left(
r\right) }$. The identity (\ref{eq29}) was used a second time for that. The
same estimates work also for the other partial derivatives $\left| \frac{%
\partial \,}{\partial t_{j}}C_{r}\left( Q\right) \left( t\right) \right| $,
and summing the estimates on $\sum_{j=1}^{3}\left| \frac{\partial \,}{%
\partial t_{j}}C_{r}\left( Q\right) \right| $, we get
\begin{equation*}
\left\| \nabla C_{r}\left( Q\right) \right\| \leq \frac{1}{r}\left\| \nabla
Q\right\| \cdot \left( 1+\frac{3\pi }{2\left( r-1\right) ^{3\mathstrut }}%
\sin \left( \frac{\pi }{r-1}\right) \right) \text{.}
\end{equation*}

\begin{summary}
With respect to this gradient-norm, $C_{r}$
has $\smash[b]{\frac{1}{r}\cdot \left( 1+\frac{3\pi }{2\left( r-1\right)
^{3\mathstrut }}\sin \left( \frac{\pi }{r-1}\right) \right) }$
as operator bound.
\end{summary}

This estimate
does not get us contractivity of $Q\mapsto
C_{r}\left( Q\right) $ (on functions $Q\in C^{\infty }$, restricted to $%
\mathbf{Y}\left( r\right) $, and satisfying $Q\left( 0\right) =0$), when $r=2
$; but, when $r=3$, we get the number $\frac{1}{3}\left( 1+\frac{3\pi }{%
2\cdot \left( 3-1\right) ^{3\mathstrut }}\right) \approx .53<1$ as an upper
bound on the contractivity constant.
\end{proof}

\begin{remark}
\label{RemNew7.7}Note that, rather than using the simplex $\mathbf{Y}\left(
r\right) $ in the estimates for the operator $C_{r}$ of Example \ref
{ExaNew7.3}, we could alternatively have used simply the enveloping unit
cube $\mathbf{B}$ positioned as follows: $-1\leq t_{j}\leq 0$, $j=1,2,3$, in 
$\mathbb{R}^{3}$. (See Figure \ref{simplex}.) But its volume is one, as
opposed to
$\limfunc{vol}\nolimits_{3}\left( \mathbf{Y}\left( r\right) \right) =
\smash[b]{\frac{1}{3\left( r-1\right) ^{3\mathstrut }}}$, which gives
the present much better
contractivity constant for $C_{r}$.
\end{remark}

\section{\label{S8}A Scaling Property for Matrix Dilation}

Let the system $\left( R,B,L\right) $ be as described, i.e., $R$ is a $\nu $%
-by-$\nu $ integral matrix which is strictly expansive. The subsets $B$ and $%
L$ in $\mathbb{R}^{\nu }$ will be assumed to satisfy (\ref{eq11})--(\ref
{eq12}), specifically $R\left( B\right) \subset \mathbb{Z}^{\nu }$ and $%
L\subset \mathbb{Z}^{\nu }$. If $R$ is scaled by some $r\in \mathbb{N}$,
then the new system $\left( rR,B,L\right) $ satisfies the same conditions.
We introduce the transformations 
\begin{equation}
\sigma _{r,b}\left( x\right) :=\left( rR\right) ^{-1}x+b,\quad x\in \mathbb{R%
}^{\nu },  \label{eq8.1}
\end{equation}
and the corresponding measure $\mu _{r}$ (of course also depending on $B$).
The fixed-point property of $\mu _{r}$ is 
\begin{equation}
\mu _{r}=\frac{1}{N}\sum_{b\in B}\mu _{r}^{{}}\circ \sigma _{r,b}^{-1},
\label{eq8.2}
\end{equation}
where $N=\#\left( B\right) $ ($=\#\left( L\right) $), and we recall that $%
\mu _{r}=\mu $ for $r=1$. The corresponding formula for the Fourier
transform $\widehat{\mu _{r}}$ is 
\begin{equation*}
\widehat{\mu _{r}}\left( t\right) =\chi _{B}\left( t\right) \widehat{\mu _{r}%
}\left( \left( rR\right) ^{-1}t\right) ,
\end{equation*}
where, as before, 
\begin{equation*}
\chi _{B}\left( t\right) =\frac{1}{N}\sum_{b\in B}e^{i2\pi b\cdot t},\quad
t\in \mathbb{R}^{\nu }.
\end{equation*}
The transformation corresponding to (\ref{eq20}) is 
\begin{equation}
C_{r}\left( Q\right) \left( t\right) =\sum_{l\in L}\left| \chi _{B}\left(
t-l\right) \right| ^{2}Q\left( r^{-1}R^{\ast \,-1}\left( t-l\right) \right) .
\label{eq8.3}
\end{equation}

\begin{example}
\label{Exa8.1}We have considered the following two measures, $\mu _{1}$ and $%
\mu _{2}$, in one dimension, arising from $B=\left\{ 0,\frac{1}{2}\right\} $%
, $L=\left\{ 0,1\right\} $, and $R=2$. Then $rR=2$ or $4$, for $r=1$ and $2$%
, and we showed that $\mu _{1}$ was Lebesgue measure on $\left[ 0,1\right] $
while $\mu _{2}$ was the measure with $d_{H}\left( \mu _{2}\right) =\frac{1}{%
2}$ from Theorem \ref{Thm3.4}. We considered the eigenvalue problem $%
C_{r}\left( Q\right) =Q$, and showed that $Q\equiv 1$ is the only solution
normalized to $Q\left( 0\right) =1$ when $r=2$, but not when $r=1$.

We will now show that this is general, also for the examples in higher
dimension. We show that the eigenvalue problem $C_{r}\left( Q\right) =Q$
always will have a unique (normalized) solution provided only that $r$ is
taken large enough.
\end{example}

The proof is based on two technical points:

\begin{lemma}
\label{Lem8.2}For the directional derivative $\nabla _{a}$, $a\in \mathbb{R}%
^{\nu }$, of $\left| \chi _{B}\left( t\right) \right| ^{2}$, we have 
\begin{equation*}
\nabla _{a}\left| \chi _{B}\left( t\right) \right| ^{2}=-2\pi
N^{-2}\sum_{b,b^{\prime }\in B}a\cdot \left( b-b^{\prime }\right) \sin
\left[ 2\pi \left( b-b^{\prime }\right) \cdot t\right] ,\quad t\in \mathbb{R}%
^{\nu }.
\end{equation*}

\begin{proof}%
%
An elementary computation. Recall the notation $\nabla _{a}:=a\cdot \nabla
=\sum_{j=1}^{\nu }a_{j}\frac{\partial \,}{\partial t_{j}}$ relative to the
canonical coordinate system in $\mathbb{R}^{\nu }$.%
\end{proof}%
%
\end{lemma}

With the assumptions from above on the system $\left(
R,B,L\right) $ in $\mathbb{R}^{\nu }$, we have the following estimate for
the contractivity constant $\gamma $ of 
\begin{equation*}
C\left( Q\right) \left( t\right) =\sum_{l\in L}\left| \chi _{B}\left(
t-l\right) \right| ^{2}Q\left( R^{*\,-1}\left( t-l\right) \right) .
\end{equation*}
Let $\mathbf{Y}$ denote the convex hull of the attractor $\mathbf{X}\left(
L\right) $ of the maps 
\begin{equation}
\rho _{l}\left( t\right) :=R^{*\,-1}\left( t-l\right) ,\quad t\in \mathbb{R}%
^{\nu },  \label{eq8.5}
\end{equation}
indexed by $l\in L$. 
%

We replace the norm
$\left\| \,\cdot \, \right\| _{\infty ,\,\mathbf{Y}}$ from Section \ref
{S6} by the following alternative norm: 
\begin{equation}
\left\| Q\right\| _{1,\,\mathbf{Y}}
:=\left\| \left| \nabla Q\right| _{2}\right\|
_{1}=\int_{\mathbf{Y}}\left| \nabla Q\right| _{2}\,dm,
  \label{eq:1-norm}
\end{equation}
where $m$ denotes Lebesgue measure%
.

\begin{theorem}[$L^{1}$-Theorem]
\label{Thm8.3}
Let $(R,B,L)$ be a system in $\mathbb{R}^{\nu }$ satisfying
\textup{(\ref
{B=L})--(\ref{Hadamard})}, $0\in L$. Let
$C$ be the operator given by \textup{(\ref
{eq21}
)}, and let $\mathbf{Y}$ denote the convex hull of the attractor
$\mathbf{X}_{\rho
}$%
. Let $\left\| Q\right\| _{1,\,\mathbf{Y}}$ be given by
\textup{(%
\ref{eq:1-norm})}, and let 
\begin{equation*}
\beta :=2\pi \func{diam}(B)\max_{\substack{b,b^{\prime }\in B \\ l\in
L}}\left\|
\sin (2\pi (b-b^{\prime })(\,\cdot \, -l))\right\| _{\infty }.
\end{equation*}
Let
$\left\| T\right\| _{op}$ be
the operator norm, and $\left\| T\right\|
_{hs}:=\left( \sum
_{j,k=1}^{\nu }\left| t_{j,k}\right| ^{2}\right)
^{1/2}$ 
the Hilbert-Schmidt norm for a
$\nu \times \nu $ matrix $T$. Then we
have 
\begin{equation}
\label{eq8pound1}
\left\| CQ\right\| _{1,\,\mathbf{Y}}
\leq \left| \det R\right| \left[ \left(
1-N^{-1}\right) \beta \left\| R^{-1}\right\| _{op}\max_{l\in L}\left|
l\right| _{2}+N\left\| R^{-1}\right\| _{hs}\right] \left\| Q\right\|
_{1,\,\mathbf{Y}},
\end{equation}
for any $C^{1}$-function $Q$ such that
$Q(0)=0$. If further $\mathbf{Y}\cap (\mathbf{Y}-l)$ is
a
set of Lebesgue measure zero for any $l\in L$ with
$l\neq 0$, then we have
the sharper estimate 
\begin{equation}
\label{eq8pound2}
\left\| CQ\right\| _{1,\,\mathbf{Y}}
\leq \left| \det R\right| \left[ \left(
1-N^{-1}\right) \beta \left\| R^{-1}\right\| _{op}\max_{l\in L}\left|
l\right| _{2}+\left\| R^{-1}\right\| _{hs}\right] \left\| Q\right\|
_{1,\,\mathbf{Y}},
\end{equation}
for any $C^{1}$-function $Q$ such that $Q(0)=0$.
\end{theorem}

\begin{remark}
\label{Rem8ins1}
Note that this result is pointless in one dimension
($\nu =1$),
for then
\begin{equation*}
\left| \det R\right| \left\| R^{-1}\right\| _{hs}=1,
\end{equation*}
while this product gets small in higher dimensions.
In our \cite{JoPe96}
setup, the $L^{1}$-Theorem can only be used if $%
\left| \det R\right| \left\| R^{-1}\right\| _{hs}<1$ is possible for
some
invertible integer matrices $R$. Actually this may work better if we take
the generalized setup 
(\ref{eq12bis})
seriously, and, e.g., take $R=\frac{1}{2}I_{\nu }$, since
then $%
\left| \det R\right| \left\| R^{-1}\right\| _{hs}\rightarrow 0$ as $%
\nu \rightarrow \infty $. 
However, it is
not yet clear
if there are any
nontrivial
examples
allowing us to use only the
$L^{1}$-Theorem. 
\end{remark}

\begin{proof}
We have 
\begin{align*}
\left\| CQ\right\| _{1,\,\mathbf{Y}} 
&=
\vphantom{\sum_{l\in L}}
\left\| \left| \nabla \vphantom{\sum}\smash{\sum
_{l\in
L}}\left| \chi _{B}(\,\cdot \, -l)\right| ^{2}Q\circ \rho _{l}\right|
_{2}\right\| _{1} \\
&=
\vphantom{\sum_{l\in L}}
\left\| \left| \vphantom{\sum}\smash{\sum
_{l\in L}}\left( \nabla \left| \chi _{B}(%
\,\cdot \, -l)\right| ^{2}\right) Q\circ \rho _{l}+
\vphantom{\sum}\smash{\sum
_{l\in
L}}\left|
\chi _{B}(\,\cdot \, -l)\right| ^{2}\nabla \left( Q\circ \rho _{l}\right)
\right| _{2}\right\| _{1} \\
&\leq 
\vphantom{\sum_{l\in L}}
\left\| \left| \vphantom{\sum}\smash{\sum
_{l\in L}}\left( \nabla \left| \chi
_{B}(%
\,\cdot \, -l)\right| ^{2}\right) Q\circ \rho _{l}\right| _{2}\right\|
_{1}+
\vphantom{\sum_{l\in L}}
\left\| \left| \vphantom{\sum}\smash{\sum
_{l\in L}}\left| \chi _{B}(\,\cdot \, %
-l)\right| ^{2}\nabla \left( Q\circ \rho _{l}\right) \right| _{2}\right\|
_{1}.
\end{align*}
As before we
estimate the two terms in this sum, starting
with the first term. Note that Lemma \ref{LemNew7.5} implies 
\begin{equation*}
\sum
_{l\in L}\left( \nabla \left| \chi _{B}(t-l)\right|
^{2}\right)
=0
\end{equation*}
for all $t$. Setting $L^{*}:=L\diagdown \{0\}$, it follows that 
\begin{multline*}
\sum
_{l\in L}\left( \nabla \left| \chi _{B}(t-l)\right|
^{2}\right) Q\circ \rho _{l}(t) 
=\sum
_{l\in L^{*}}\left( \nabla \left| \chi _{B}(t-l)\right|
^{2}\right) \left( Q\circ \rho _{l}(t)-Q\circ \rho _{0}(t)\right)  \\
=-\sum
_{l\in L^{*}}\left( \nabla \left| \chi _{B}(t-l)\right|
^{2}\right) \int_{0}^{1}\left( R^{*\,-1}l\right) \cdot \left( \nabla Q\left(
R^{*\,-1}t-sR^{*\,-1}l\right) \right) \,ds.
\end{multline*}
for all $t$. So we have 
\begin{multline*}
\vphantom{\sum_{l\in L}}
\left\| \left| \vphantom{\sum}\smash{\sum
_{l\in L}}\left( \nabla \left| \chi _{B}(%
\,\cdot \, -l)\right| ^{2}\right) Q\circ \rho _{l}\right| _{2}\right\| _{1}
\\
\begin{aligned}
{}&\leq 
\vphantom{\sum_{l\in L^{*}}}
\left\| \left| \vphantom{\sum}\smash{\sum
_{l\in L^{*}}}\left( \nabla \left| \chi
_{B}(\,\cdot \, -l)\right| ^{2}\right) \right| _{2}\left| \int_{0}^{1}\left(
R^{*\,-1}l\right) \cdot \left( \nabla Q\left( R^{*\,-1}\,\cdot \,
-sR^{*\,-1}l\right)
\right) \,ds\right| \right\| _{1} \\
&\leq 
\vphantom{\sum_{l\in L^{*}}}
\left\| \left| \vphantom{\sum}\smash{\sum
_{l\in L^{*}}}\left( \nabla \left| \chi
_{B}(\,\cdot \, -l)\right| ^{2}\right) \right| _{2}\left| R^{*\,-1}l\right|
_{2}\int_{0}^{1}\left| \nabla Q\left( R^{*\,-1}\,\cdot \, -sR^{*\,-1}l\right)
\right| _{2}\,ds\right\| _{1} \\
&\leq \max_{l\in L^{*}}\left\| \left| \nabla \left| \chi _{B}(\,\cdot \, %
-l)\right| ^{2}\right| _{2}\right\| _{\infty }\max_{l\in L^{*}}\left|
R^{*\,-1}l\right| _{2}
\vphantom{\sum_{l\in L^{*}}}
\left\| \vphantom{\sum}\smash{\sum
_{l\in L^{*}}}\int_{0}^{1}\left|
\nabla Q\left( R^{*\,-1}\,\cdot \, -sR^{*\,-1}l\right) \right| _{2}\,ds\right\|
_{1}.
\end{aligned}
\end{multline*}
We continue the estimate by estimating the last factor in this product. 
\begin{multline*}
\vphantom{\sum_{l\in L^{*}}}
\left\| \vphantom{\sum}\smash{\sum
_{l\in L^{*}}}\int_{0}^{1}\left| \nabla Q\left(
R^{*\,-1}t-sR^{*\,-1}l\right) \right| _{2}\,ds\right\| _{1}  \\
\begin{aligned}
{}&=\int_{\mathbf{Y}}\sum
_{l\in L^{*}}\int_{0}^{1}\left| \nabla Q\left(
R^{*\,-1}t-sR^{*\,-1}l\right) \right| _{2}\,ds\,dt \\
&=\int_{0}^{1}\sum
_{l\in L^{*}}\int_{\mathbf{Y}}\left| \nabla Q\left(
R^{*\,-1}t-sR^{*\,-1}l\right) \right| _{2}\,dt\,ds.
\end{aligned}
\end{multline*}
Making the substitution $u=R^{*\,-1}t-sR^{*\,-1}l=\rho _{s,l}(t)$ we have $%
dt=\left| \det R\right| \,du$, so that 
\begin{align*}
\vphantom{\sum_{l\in L^{*}}}
\left\| \vphantom{\sum}\smash{\sum
_{l\in L^{*}}}\int_{0}^{1}\left| \nabla Q\left(
R^{*\,-1}t-sR^{*\,-1}l\right) \right| _{2}\,ds\right\| _{1} 
&=\int_{0}^{1}\sum
_{l\in L^{*}}\int_{\rho _{s,l}(\mathbf{Y})}\left|
\nabla
Q\left( u\right) \right| _{2}\left| \det R\right| \,du\,ds \\
&\leq \int_{0}^{1}\sum
_{l\in L^{*}}\int_{\mathbf{Y}}\left| \nabla Q\left(
u\right) \right| _{2}\left| \det R\right| \,du\,ds \\
&\leq (N-1)\left| \det R\right| \left\| \left| \nabla Q\right|
_{2}\right\|
_{1}.
\end{align*}
The second to last estimate used 
the fact
that $\rho _{s,l}(\mathbf{Y})\subset \mathbf{Y}$, 
which
in
turn is a consequence of the fact that $\left( \rho _{s,l}(t)\right)
_{0\leq %
s\leq 1}$ is the line connecting $\rho _{l}(t)$ and $\rho _{0}(t),$ the
invariance $\rho _{l}\mathbf{Y}\subset \mathbf{Y}$,
and the convexity of $\mathbf{Y}$. So we have
shown 
\begin{multline*}
\vphantom{\sum_{l\in L}}
\left\| \left| \vphantom{\sum}\smash{\sum
_{l\in L}}\left( \nabla \left| \chi _{B}(%
\,\cdot \, -l)\right| ^{2}\right) Q\circ \rho _{l}\right| _{2}\right\| _{1}
\\
\leq \max_{l\in L^{*}}\left\| \left| \nabla \left| \chi _{B}(\,\cdot \, %
-l)\right| ^{2}\right| _{2}\right\| _{\infty }\max_{l\in L^{*}}\left|
R^{*\,-1}l\right| _{2}(N-1)\left| \det R\right| \,\left\| \left| \nabla
Q\right| _{2}\right\| _{1}.
\end{multline*}
Estimating the $\left\| R^{*\,-1}l\right\| _{2}$ term,
\begin{equation*}
\vphantom{\sum_{l\in L}}
\left\| \left| \vphantom{\sum}\smash{\sum
_{l\in L}}\left( \nabla \left| \chi
_{B}(\,\cdot \, %
-l)\right| ^{2}\right) Q\circ \rho _{l}\right| _{2}\right\| _{1}\leq %
(1-N^{-1})\beta \left\| R^{-1}\right\| _{op}\max_{l\in L}\left| l\right|
_{2}\left| \det R\right| \left\| Q\right\| _{1,\,\mathbf{Y}}.
\end{equation*}
The 
other 
estimate 
is 
\begin{align*}
\vphantom{\sum_{l}}
\left\| \vphantom{\sum}\smash{\sum
_{l}}\left| \chi _{B}(\,\cdot \, -l)\right| ^{2}\left|
\nabla (Q\circ \rho _{l})\right| _{2}\right\| _{1}
&=\int_{\mathbf{Y}}\sum
_{l}\left| \chi _{B}(t-l)\right| ^{2}\left| \nabla
\left( Q\circ \rho _{l}\right) (t)\right| _{2}\,dt \\
&\leq \sum
_{l}\int_{\mathbf{Y}}\left| \nabla \left( Q\circ \rho
_{l}\right)
(t)\right| _{2}\,dt \\
&\leq \sum
_{l}\int_{\mathbf{Y}}\left| \left( \nabla Q\right) \circ \rho
_{l}(t)\right| _{2}\left\| R^{-1}\right\| _{hs}\,dt \\
&=\left| \det R\right| \sum
_{l}\int_{\rho _{l}\mathbf{Y}}\left| \left(
\nabla Q\right) (u)\right| _{2}\left\| R^{-1}\right\| _{hs}\,du \\
&\leq \left| \det R\right| N\left\| \left| \nabla Q\right| _{2}\right\|
_{1}\left\| R^{-1}\right\| _{hs}.
\end{align*}
The last inequality used
$\rho _{l}\mathbf{Y}\subset \mathbf{Y}$; and it can be improved if
$\mathbf{Y}\cap
(\mathbf{Y}-l)$ is a set of Lebesgue measure zero for every
$l\in L^{*}$. In that
case, we have 
\begin{equation*}
\vphantom{\sum_{l}}
\left\| \vphantom{\sum}\smash{\sum
_{l}}\left| \chi _{B}(\,\cdot \, -l)\right| ^{2}\left|
\nabla (Q\circ \rho _{l})\right| _{2}\right\| _{1}\leq \left| \det
R\right|
\left\| \left| \nabla Q\right| _{2}\right\| _{1}\left\| R^{-1}\right\|
_{hs}.
\end{equation*}
It is now easy to complete the proof. 
\end{proof}

\begin{remark}
\label{RemNew8.4}
The following additional
arguments may be used in
a sequence of steps (to be described
below) leading to the
elimination of the factor
$\left| \det \left( R\right) \right| $ in the two
estimates (\ref{eq8pound1}) and (\ref{eq8pound2}) in
the conclusion of Theorem \ref{Thm8.3}:
After estimating $\left\| \nabla Q\right\| _{L^{1}\left( \mathbf{Y}\right) }$
by a constant times
\begin{equation}
\label{eq8ins1}
\sum _{l\in L}\int _{\mathbf{Y}}\left| \left( \nabla Q\right) \circ
\rho _{l}\right| \,dt
\end{equation}
for $Q\in C^{\infty }\left( \mathbb{R}^{\nu }\right) $ such that
$Q\left( 0\right) =0$,
we may change variables in the
$\int _{\mathbf{Y}}\left( \cdots \right) \,dt$
terms in the summation (\ref{eq8ins1}).
Let $\mathbf{Y}$ have the non-overlap
property (\ref{eq5}), i.e., suppose
that the Lebesgue measure of
the possible overlap in the union
\begin{equation}
\label{eq8ins2}
T\left( \mathbf{Y}\right) :=\bigcup _{l\in L} \rho _{l}
\left( \mathbf{Y}\right) 
\end{equation}
is zero. For the sum in (\ref{eq8ins1}),
we then have
\begin{equation*}
\sum _{l\in L}\int _{\mathbf{Y}}\left| \nabla Q\right| \circ 
\rho _{l}\,dt=\left| \det \left( R\right) \right| 
\int _{T\left( \mathbf{Y}\right) }\left| \nabla Q\right| \,dt
\end{equation*}
and
\begin{equation}
\label{eq8ins3}
\begin{aligned}
\int _{T\left( \mathbf{Y}\right) }\left| \nabla Q\right| \,dt
&\leq 
\sup _{T\left( \mathbf{Y}\right) }\left| \nabla Q\right| 
\limfunc{vol}\nolimits_{\nu }\left( T\left( \mathbf{Y}\right) \right)  \\
&=
\sup _{T\left( \mathbf{Y}\right) }\left| \nabla Q\right| 
\frac{N\limfunc{vol}\nolimits_{\nu }\left( \mathbf{Y}\right) }
{\left| \det \left( R\right) \right| } .
\end{aligned}
\end{equation}
Then combining estimates (\ref{eq8ins1}) and (\ref{eq8ins3}),
we get
\begin{equation*}
\sum _{l\in L}\int _{\mathbf{Y}}
\left| \left( \nabla Q\right) \circ 
\rho _{l}\right| \,dt\leq 
N\limfunc{vol}\nolimits_{\nu }\left( \mathbf{Y}\right) 
\sup _{T\left( \mathbf{Y}\right) }\left| \nabla Q\right| .
\end{equation*}
When combined in turn with the
estimate from Theorem \ref{Thm8.3},
we get
\begin{equation*}
\left\| \nabla C\left( Q\right) \right\| _{L^{1}\left( \mathbf{Y}\right) }
\leq N\sup _{T\left( \mathbf{Y}\right) }\left| \nabla Q\right| 
F_{\nu }\left( R\right) 
\limfunc{vol}\nolimits_{\nu }\left( \mathbf{Y}\right) ,
\end{equation*}
where $F_{\nu }\left( R\right) $ is an expression
(see Theorem \ref{Thm8.3} and Section \ref{S9}) which depends on
the operator norm, and the
Hilbert-Schmidt norm, applied
to $R^{-1}$. In the case of
Example \ref{ExaNew7.3} (``The Eiffel Tower''),
$F_{3}\left[ 
\left( \begin{smallmatrix}r&0&0\\0&r&0\\0&0&r\end{smallmatrix}\right) 
\right] =\func{const.}\cdot r^{-1}$.

When the same method is
applied to $\nabla C^{n}\left( Q\right) $, it
gives a new tool for showing
$\lim _{n\rightarrow \infty }C^{n}\left( Q\right) =0$ when
$Q\left( 0\right) =0$.
Since
\begin{equation}
\label{eq8ins4}
\mathbf{X}=\bigcap _{n}T^{n}\left( \mathbf{Y}\right) ,
\end{equation}
and $\mathbf{X}$ has Lebesgue measure zero when
$N<\left| \det R\right| $, the limits can be
estimated further.

Recall that the $T^{n}$-terms in
the intersection (\ref{eq8ins4}) are defined
from iteration of (\ref{eq8ins2}), i.e.,
\settowidth{\qedskip}{$\displaystyle
\bigcup 
_{\left( l_{1},\dots ,l_{n}\right) \in \bigcross _{1}^{n}L}\bigcup $}
\begin{equation*}
T^{n}\left( \mathbf{Y}\right) =
\makebox[0.5\qedskip][l]{$\displaystyle
\bigcup 
_{\left( l_{1},\dots ,l_{n}\right) \in \bigcross _{1}^{n}L}
$\hss }\,
\rho_{l_{1}}\circ \dots \circ \rho_{l_{n}}\left( \mathbf{Y}\right) ,
\end{equation*}
and also that we are finished
after showing $C^{n}\left( Q_{0}\right) \rightarrow 0$
for all $Q_{0}\in C^{\infty }\left( \mathbb{R}^{\nu }\right) $
such that $Q_{0}\left( 0\right) =0$,
because 
the fixed point
$C\left( Q_{1}\right) =Q_{1}$ ($Q_{1}=\mathbf{1}-Q_{0}$)
will then be unique
subject to $Q_{1}\left( 0\right) =1$.
\end{remark}

\begin{remark}
\label{Rem8.4}The constant $\beta $ in (\ref{eq8pound1}) may be easily
calculated for the three-di\-men\-sion\-al example in Figure \ref{eiffel},
i.e., Example \ref{ExaNew7.3}. We have $\beta =2\pi \limfunc{diam}\left(
B\right) =2\pi \frac{1}{\sqrt{2}}=\pi \sqrt{2}$. 
Calculating 
$\det \left( R\right)
$ in this example, 
it follows from the formula for
$F_{3}\left( R\right) $ 
that
the number $%
\left( \pi \sqrt{2}%
+1\right) F_{3}\left( R\right) $ ($<1$)
is a contractive upper bound for the constant $\gamma $
when $R=\left( 
\begin{smallmatrix}
r & 0 & 0 \\
0 & r & 0 \\
0 & 0 & r
\end{smallmatrix}
\right) $, and $r\in \mathbb{N}$ is taken sufficiently
large. See formulas (\ref{T:estimate:gamma})
and (\ref{T:estimate:boundgradsupnorm})
in Section \ref{S9} for details. We expect $r=2$ will
suffice.
(Hence we get strict contractivity in this new metric as
claimed.) An application of Lemma \ref{Lem3.3} then shows that
$L^{2}\left( \mu _{r}\right) $
of Figure \ref{eiffel} has an orthonormal basis
$\left\{ e_{\lambda }:\lambda \in P_{r}\right\} $
when $P_{r}$ is given in (\ref{eqNew7.1}).
\end{remark}

\begin{corollary}
\label{Cor8.4}When the contractivity constant $\gamma _{r}$ is calculated
for the system $\left( rR,B,L\right) $ as described above, then $%
\lim_{r\rightarrow \infty }\gamma _{r}=0$, and it follows that $C_{r}$ is
strictly contractive when $r$ increases. As a consequence,
if $L$ is not contained in a hyperplane in $\mathbb{R}^{\nu }$, then
the exponentials $%
\left\{ e_{\lambda }:\lambda \in P_{r}\left( L\right) \right\} $ will form
an orthonormal basis for $L^{2}\left( \mu _{r}\right) $ where $P_{r}\left(
L\right) =\left\{ l_{0}+rR^{\ast }l_{1}+%
\dots%
%
:l_{i}\in L\text{, finite sums}\right\} $.
\end{corollary}

\begin{remark}
\label{Rem8.5}The formula 
for the contractivity constant $%
\gamma _{r}$ contains a factor $
r
^{-1}$%
, but we also note that the convex hull $\mathbf{Y}_{r}$ of $\mathbf{X}%
_{r}\left( L\right) $ (= the attractor of $\{\rho _{r,l}:l\in L\}$)
decreases as a function of $r$.
(See the next section for details.)
\end{remark}

\section{\label{S9}Estimates and Duality of Scales: Fractals in the Large}

In this section, we provide
a context in $\mathbb{R}^{\nu }$ (the dimension
$\nu $ arbitrary) which utilizes
two types of duality: first the
traditional time/frequency duality
of harmonic analysis, and secondly
the duality between the iteration
limits when the scale (given
by a relative integral $\nu $-by-$\nu $ matrix $R$)
gets arbitrarily small, versus when
it gets arbitrarily large.

The axioms below are such that
the scaling matrix $R$ is only
required to be integral with
respect to the two given
(finite) translation sets $B$ and
$L$; see (\ref{eq12bis}), or (\ref{Compatibility}) below.
In this general setup, we give
explicit estimates on
the contractivity constant,
making it clear how it depends
on the given data $(R,B,L)$ and
the dimension $\nu $. The estimates are
gradient-supremum-norm ones.

Consider a triple $(R,B,L)$ such that $R$ is an expansive
$\nu \times \nu $ matrix
with real entries and $B$ and $L$ are subsets of $\mathbb{R}^{\nu }$
satisfying:
\begin{align}
N &:=\#B=\#L;  \label{B=L} \\
R^{n}b\cdot l &\in \mathbb{Z},\text{ for any }n\in \Bbb{N},\;b\in B,\;l\in L;
\label{Compatibility} \\
H_{B,L} &:=N^{-1/2}\left( e^{i2\pi b\cdot l}\right) _{b\in B,\;l\in L}\text{
is a unitary }N\times N\text{ matrix.}  \label{Hadamard}
\end{align}
We introduce two dynamical systems 
\begin{align*}
\sigma _{b}(x) &:=R^{-1}x+b; \\
\tau _{l}(x) &:=R^{*}x+l
\end{align*}
and the corresponding inverse functions 
\begin{align*}
\omega _{b}(x) &:=R(x-b); \\
\rho _{l}(x) &:=R^{*\,-1}(x-l)
\end{align*}
where $x\in \mathbb{R}^{\nu }$, $b\in B$, and $l\in L$. We showed in
\cite{JoPe94}
that one may think of $\left( \sigma _{b}\right) $ and $\left( \tau
_{l}\right) $ as dual dynamical systems. It is known that there are unique
non-empty compact sets $\mathbf{X}%
_{\sigma }$ and $\mathbf{X}%
_{\rho }$ such that 
\begin{align}
\mathbf{X}%
_{\sigma }&=\bigcup_{b\in B}\sigma _{b}\left( \mathbf{X}%
_{\sigma }\right)
=B+R^{-1}\mathbf{X}%
_{\sigma }  \label{sigma-recur-dom}  \\
\intertext{and} 
\mathbf{X}%
_{\rho }&=\bigcup_{l\in L}\rho _{l}(\mathbf{X}%
_{\rho })=R^{*\,-1}(\mathbf{X}%
_{\rho }-L);
\label{rho-recur-dom}
\end{align}
and it is easy to see that 
\begin{equation*}
\mathbf{X}%
_{\sigma }=
\vphantom{\sum_{k=0}^{\infty }}
\left\{ \vphantom{\sum}\smash{\sum_{k=0}^{\infty }}
R^{-k}b_{k}:b_{k}\in B\right\} \text{\quad 
and \quad }\mathbf{X}%
_{\rho }=
\vphantom{\sum_{k=0}^{\infty }}
\left\{ -\vphantom{\sum}\smash{\sum_{k=1}^{\infty }}
R^{*\,-k}l_{k}:l_{k}\in L\right\} .
\end{equation*}
By uniqueness, these sets are the closures of the
respective orbits of $0$ under the maps 
$\left( \sigma _{b}\right) $, respectively $\left( \rho _{l}\right) $. In
\cite{JoPe96}, we argued that it is natural to think of $\mathbf{X}%
_{\rho }$ as
a dual of the attractor $\mathbf{X}%
_{\sigma }$. Let $\mu =\mu _{\sigma }$ be the
unique probability measure on $\mathbb{R}^{\nu }$ such that 
\begin{equation}
\mu =N^{-1}\sum_{b\in B}\mu \circ \sigma _{b}^{-1}.  \label{eq:mu}
\end{equation}
Then this measure has
$\mathbf{X}%
_{\sigma }$ as its support. Similarly, $\mathbf{X}%
_{\rho }$ is the
support of the unique probability measure
$\mu ^{\prime }=\mu _{\rho }$ satisfying 
\begin{equation}
\mu ^{\prime }=N^{-1}\sum_{l\in L}\mu ^{\prime }
\circ \rho _{l}^{-1}.  \label{eq:nu}
\end{equation}
The corresponding expansive orbits 
\begin{align}
\mathcal{L} &=\mathbf{X}%
_{\tau }:=
\vphantom{\sum_{k=0}^{n}}
\left\{ \vphantom{\sum}\smash{\sum_{k=0}^{n}}
R^{*\,k}l_{k}:n\in 
\mathbb{N}%
,l_{k}\in L\right\}  \label{eq:Ldef} \\
\mathcal{B} &=\mathbf{X}%
_{\omega }:=
\vphantom{\sum_{k=1}^{n}}
\left\{ -\vphantom{\sum}\smash{\sum_{k=1}^{n}}
R^{k}b_{k}:n\in 
\mathbb{N}%
,b_{k}\in B\right\}  \label{eq:Bdef}
\end{align}
of $0$ under $\tau $ and $\omega $ are also important. It should now be
clear that the pairs $\left( \sigma _{b},\tau _{l}\right) $ and $\left( \rho
_{l},\omega _{b}\right) $ share many properties,
which however we do not always make
explicit in the following.
Notice that the set $\mathcal{L}$
has a recursive property similar to (\ref{eq:mu}):
\begin{equation}
\mathcal{L}=L+R\mathcal{L}.  \label{eq:Lrecurrence}
\end{equation}
We discussed
uniqueness of the corresponding decompositions of elements in $%
\mathcal{L}$ above.
We also investigated when the functions
$\{e_{\lambda }:\lambda \in \mathcal{L}\}$
form an orthonormal basis for $L^{2}(\mu )$.
One may equally well ask when $%
\{e_{\lambda }:\lambda \in \mathcal{B}\}$ is an orthonormal basis for $%
L^{2}(\mu ^{\prime })$, but that is an essentially equivalent question.
It turns out that the measure $\mu ^{\prime }$ is important for our
investigation. We will see that the set $\mathcal{B}$ has properties that
imply that translates of
$\mathbf{X}%
_{\sigma }$ tile $\mathbb{R}^{\nu }$, when $N=\left| \det
R\right| $. The assumptions on $(R,B,L)$ insure that the functions $%
\{e_{\lambda }:\lambda \in \mathcal{L}\}$ are orthogonal, so we really only
have to investigate totality. 

The following result, combined with
\cite[Theorems 1.1 and 1.2]{LaWa96d},
shows that, if $%
N=\left| \det R\right| $, then $\mathbf{X}%
_{\sigma }$ and $\mathbf{X}%
_{\rho }$ are translation
tiles.
We must also check that $\mathcal{L}$
is uniformly discrete, i.e., that
there is an $\varepsilon >0$ such that
$\left\| \lambda -\lambda ^{\prime }\right\| _{2}\geq \varepsilon $ for all
$\lambda $ and $\lambda ^{\prime }$,
$\lambda \neq \lambda ^{\prime }$, in $\mathcal{L}$. That follows
from our assumptions on $\left( R,B,L\right) $,
as can be verified by an induction argument.

\begin{lemma}
\label{LemFHS2.1}Each point
$\lambda \in \mathcal{L}$ has a unique representation $\lambda
=\sum R^{*\,k}l_{k}$ with $l_{k}\in L$. Similarly each point in $\mathcal{B}$
has a unique representation of the form $\sum R^{k}b_{k}$ with $b_{k}\in B$.
\end{lemma}

\begin{proof}
Each point $\lambda \in \mathcal{L}$ has a representation of the form $%
\lambda =l+R^{*}p$ with $l\in L$ and $p\in \mathcal{L}$. Suppose $%
l-l^{\prime }+R^{*}(p-p^{\prime })=0$ where $l,l^{\prime }\in L$ and $%
p,p^{\prime }\in \mathcal{L}$. We will show that $l=l^{\prime }$,
the desired result,
then follows by induction. If $b\in B$, then 
\begin{equation}
\label{oneequalseb}
\begin{aligned}
1 &=e_{b}(l-l^{\prime }+R^{*}(p-p^{\prime })) \\
&=e_{b}(l-l^{\prime })e_{b}(R^{*}(p-p^{\prime })) \\
&=e_{b}(l-l^{\prime })e_{Rb}(p-p^{\prime }) \\
&=e_{b}(l-l^{\prime }),
\end{aligned}
\end{equation}
where the last equality used (\ref{Compatibility}). If
$l\neq l^{\prime }$ then we
have $\sum_{b\in B}e_{b}(l-l^{\prime })=0$ by (\ref{Hadamard}).
But this contradicts (\ref{oneequalseb}) because this sum
would be $N$ if (\ref{oneequalseb}) holds, hence $l=l^{\prime }$ as
we had to show.
\end{proof}

We will show that $\mathcal{L}$ has the basis property by showing that $C$
from (\ref{eq21})
acts contractively on a set of functions including
$Q_{1}$ of (\ref{eq16}) and $\mathbf{1}$
(where $\mathbf{1}$ denotes the constant function)
and
thereby showing $\sum_{\lambda \in P}\left| \hat{\mu}\left( \,\cdot \,-\lambda
\right) \right| ^{2}=\mathbf{1}$. That is, we want to use the uniqueness part
of the Banach fixed-point theorem on a suitably chosen metric space. It is
natural to consider $C$ as an operator on a set of smooth functions $Q$
satisfying $Q(0)=1$. However, the appropriate choice of a metric is not so
obvious.
The present discussion is based on gradient-supremum-norm estimates.

Let $H_{2}(\mathcal{L})$ denote the subspace of $L^{2}(\mu )$ spanned by the
orthonormal set $\{e_{\lambda }:\lambda \in \mathcal{L}\}$. Any $e_{t}$, $%
t\in \mathbb{C}^{\nu }$, is in $L^{2}(\mu )$, so
$H_{2}(\mathcal{L})$ is a subspace
of $L^{2}(\mu )$. We will show that $H_{2}(\mathcal{L})=L^{2}(\mu )$ for
certain systems $(R,B,L)$ satisfying (\ref{B=L})--(\ref{Hadamard}). Let $A$
denote the orthogonal projection of $L^{2}(\mu )$ onto $H_{2}(\mathcal{L})$,
that is 
\begin{equation*}
Af:=\sum_{\lambda \in \mathcal{L}}\left\langle e_{\lambda }\mid
f\right\rangle _{\mu }\,e_{\lambda }
\end{equation*}
for any $f\in L^{2}(\mu )$. Our next result is a useful criterion for when $%
H_{2}(\mathcal{L})=L^{2}(\mu )$.

\begin{lemma}
\label{L:diff}Let $\partial _{j}:=\frac{\partial \,}{\partial t_{j}}$. If $%
0\in L$, then $H_{2}(\mathcal{L})=L^{2}(\mu )$ if and only if $\partial
_{1}^{n_{1}}\cdots \partial _{\nu }^{n_{\nu }}Q_{1}(0)=0$ for all
$n_{j}\geq 0$ with $n_{1}+\cdots +n_{\nu }\geq 1$,
where $Q_{1}\left( t\right) 
=\sum_{\lambda \in P}\left| \hat{\mu}\left( t-\lambda
\right) \right| ^{2}$.
\end{lemma}

\begin{proof}
Restating the definition of $Q_{1}$,
we have $Q_{1}(t)=\left\| Ae_{t}\right\| _{\mu }^{2}$ for all $t$%
. If $H_{2}(\mathcal{L})=L^{2}(\mu )$, then $Ae_{t}=e_{t}$ for all $t$, hence 
$Q_{1}(t)=1$ for all $t\in \mathbb{R}^{\nu }$, in particular, $\partial
_{1}^{n_{1}}\cdots \partial _{\nu }^{n_{\nu }}Q_{1}(0)=0$ for all
$n_{j}\geq 0$ with $%
n_{1}+\cdots +n_{\nu }\geq 1$. Conversely, suppose
$\partial _{1}^{n_{1}}\cdots
\partial _{\nu }^{n_{\nu }}Q_{1}(0)=0$ for all
$n_{j}\geq 0$ with $n_{1}+\cdots
+n_{\nu }\geq 1$. By Lemma \ref{LemNew4.3}
we have a power series expansion $%
Q_{1}(z)=\sum_{n\in \mathbb{N}_{0}^{\nu }}a_{n}z^{n}$;
and, by our assumption, $a_{n}=0$
unless $n=0$. Hence $Q_{1}$ is a constant function. Since $0\in L$, we have $%
Q_{1}(t)=1$ for all $t\in \mathbb{R}^{\nu }$.
It follows that $e_{t}\in H_{2}(\mathcal{%
L})$ for all $t\in \mathbb{R}^{\nu }$,
so an application of the Stone--Weierstrass
theorem implies that $H_{2}(\mathcal{L})$ contains all continuous functions
on $\mathbf{X}%
_{\sigma }$, hence $H_{2}(\mathcal{L})$, being a dense subspace of $%
L^{2}(\mu )$, must be equal to $L^{2}(\mu )$ as we wanted to show.
\end{proof}

Let $\mathbf{Y}%
$ denote the convex hull of the attractor $\mathbf{X}%
_{\rho }$ given by (\ref
{rho-recur-dom}), let $\left\| Q\right\| _{\infty }:=\sup_{y\in \mathbf{Y}%
}\left|
Q(y)\right| $, and let
\begin{equation}
\left\| Q\right\| _{\mathbf{Y}%
}:=\left\| \left| \nabla Q\right| _{2}\right\| _{\infty
}  \label{eq:Y-sup-norm}
\end{equation}
where $\left| z\right| _{2}:=\left( \sum_{j=1}^{\nu }\left| z_{j}\right|
^{2}\right) ^{1/2}$ is the usual Hilbert norm on
$\mathbb{C}^{\nu }$. We begin by
finding an explicit formula for the operator norm of $C$ acting on a
suitable set of smooth functions.

\begin{theorem}
\label{th:Cbound}Let $(R,B,L)$ be a system in
$\mathbb{R}^{\nu}$ satisfying \textup{(\ref
{B=L})--(\ref{Hadamard})}, $0\in L$.
Let $C$ be the operator given by \textup{(\ref
{eq21})}, let $\mathbf{Y}%
$ denote the convex hull of the attractor $\mathbf{X}%
_{\rho }$
given by \textup{(\ref{rho-recur-dom})},
and let $\left\| Q\right\| _{\mathbf{Y}%
}$ be given by
\textup{(\ref{eq:Y-sup-norm})}.
If $L$ spans $\mathbb{R}^{\nu }$, and if there exists $%
\gamma <1$ so that $\left\| CQ\right\| _{\mathbf{Y}%
}\leq \gamma \left\| Q\right\|
_{\mathbf{Y}%
} $ for all $Q$ in a set of $C^{1}$-functions containing $\mathbf{1}-Q_{1}$,
where $Q_{1}\left( t\right) 
=\sum_{\lambda \in P}\left| \hat{\mu}\left( t-\lambda
\right) \right| ^{2}$,
then $H_{2}(\mathcal{L})=L^{2}(\mu ).$
\end{theorem}

\begin{proof}
Since $C\mathbf{1}=\mathbf{1}$ and $CQ_{1}=Q_{1}$ we have 
\begin{equation*}
\left\| \mathbf{1}-Q_{1}\right\| _{\mathbf{Y}%
}=\left\| C(\mathbf{1}-Q_{1})\right\| _{\mathbf{Y}%
}\leq
\gamma \left\| \mathbf{1}-Q_{1}\right\| _{\mathbf{Y}%
}
\end{equation*}
so $Q_{1}(y)=1$ for all $y\in \mathbf{Y}%
$. Since $\mathbf{Y}%
$ is a convex set containing
$0$, and $%
\mathbf{Y}%
 $ is not contained in any proper subspace of
$\mathbb{R}^{\nu}$, it follows
that
$\partial _{1}^{n_{1}}\cdots \partial _{\nu }^{n_{\nu }}Q_{1}(0)=0$
for all $%
n_{j}\geq 0 $ with $n_{1}+\cdots +n_{\nu}\geq 1$.
Hence an application of
Lemma \ref{L:diff} completes the proof.
\end{proof}

\begin{remark}
\label{RemFHS3.1}In the statement
and proof of Theorem \ref{th:Cbound},
$\mathbf{Y}%
$ could be replaced
by any set containing $\mathbf{X}%
_{\rho }$.
But $\mathbf{Y}%
$ will be chosen so as to produce optimal estimates.

The restricting assumption in Theorem \ref{th:Cbound}
(and also in the corresponding results
in Sections \ref{S5} and \ref{S6}) \emph{that the set
$L\subset \mathbb{R}^{\nu }$ not be contained in
a hyperplane} is necessary. In
\cite[Example 7.3]{JoPe96}, we considered a
triple $\left( R,B,L\right) $ in $\mathbb{R}^{2}$
where $\nu =2$, $N=3$, and where the
conditions in Theorem \ref{th:Cbound} all hold,
\emph{except} that $L$, in this \cite{JoPe96} example,
\emph{is} contained in a line
in the plane. Specifically, take
$R=\left( \begin{smallmatrix}6&0\\0&6\end{smallmatrix}\right) $,
$B=\left\{ \left( \begin{smallmatrix}0\\0\end{smallmatrix}\right) , 
\left( \begin{smallmatrix}\frac{1}{2}\\0\end{smallmatrix}\right) , 
\left( \begin{smallmatrix}0\\\frac{1}{2}\end{smallmatrix}\right) \right\} $,
$L=\left\{ \left( \begin{smallmatrix}0\\0\end{smallmatrix}\right) ,
\pm l_{1}\right\} $ where $l_{1}=\frac{2}{3} 
\left( \begin{smallmatrix}1\\-1\end{smallmatrix}\right) $.
Properties (\ref{B=L})--(\ref{Hadamard}) are
clearly satisfied for the example.
Let $\mu $ be the unique probability measure on $\mathbb{R}^{2}$
satisfying
\begin{align*}
\hat{\mu}\left( t\right) &=
\chi _{B}\left( t\right) \hat{\mu}\left( \frac{t}{6}\right) ,
\quad t\in \mathbb{R}^{2},  \\
\intertext{where}
\chi _{B}\left( t\right) &=
\frac{1}{3}\left( 1+e^{i\pi t_{1}}+e^{i\pi t_{2}}\right) .
\end{align*}
Let $P\left( L\right) $, $Q_{1}$, and $C$ be defined (as usual):
\begin{align*}
P\left( L\right) &=
\left\{ l_{0}+6l_{1}+6^{2}l_{2}+
\dots
:l_{i}\in L\text{, finite sums}\right\} , \\
Q_{1}\left( t\right) &:=
\sum _{\lambda \in P\left( L\right) }
\left| \hat{\mu}\left( t-\lambda \right) \right| ^{2}, \\
\intertext{and}
C\left( f\right) \left( t\right) &=
\sum _{l\in L}
\left| \chi _{B}\left( t-l\right) \right| ^{2}
f\left( \frac{t-l}{6}\right) ,\quad t\in \mathbb{R}^{2}.
\end{align*}
From Lemma \ref{Lem4.1}, we get
the eigenvalue identity
$C\left( Q_{1}\right) =Q_{1}$ satisfied, and
conclude that the exponentials
$E\left( P\left( L\right) \right) =
\left\{ e_{\lambda }:\lambda \in P\left( L\right) \right\} $
are orthogonal in $L^{2}\left( \mu \right) $.
Moreover, the argument from
the proof of Proposition \ref{ProNew7.6}
gets us \emph{contractivity} for $C$, i.e.,
\begin{equation}
\label{contractivityforC}
\left\| \nabla C\left( f\right) \right\| _{\mathbf{Y}}
\leq \gamma \left\| \nabla f\right\| _{\mathbf{Y}}
\end{equation}
for some constant $\gamma <1$,
and for all $f\in C^{\infty }\left( \mathbb{R}^{2}\right) $ such that
$f\left( 0\right) =0$; where
\begin{equation*}
\left\| \nabla f\right\| _{\mathbf{Y}} :=
\sup _{t\in \mathbf{Y}}
\vphantom{\sum _{j=1}^{2}}
\left( \vphantom{\sum}\smash{\sum _{j=1}^{2}}
\left| \partial _{j}f\left( t\right) \right| \right) ,
\end{equation*}
and $\mathbf{Y}$ is the $L$-attractor.
But $\mathbf{Y}$ is a finite \emph{line} in $\mathbb{R}^{2}$
by Lemma \ref{LemNew7.5}; in fact,
$\mathbf{Y}=
\left\{ \left( \begin{smallmatrix}t_{1}\\t_{2}\end{smallmatrix}
\right) :t_{2}=-t_{1},\,\left| t_{1}\right| \leq \frac{2}{15}\right\} $.
From the uniqueness of the Banach
fixed point, we get that
$Q_{1}\equiv 1$ on $\mathbf{Y}$, but \emph{not} on
$\mathbb{R}^{2}$.
Indeed, if we had $Q_{1}\equiv 1$ on 
$\mathbb{R}^{2}$, then it would
follow from Lemma \ref{L:diff} that
$H_{2}\left( P\left( L\right) \right) =L^{2}\left( \mu \right) $,
i.e., that
$E\left( P\left( L\right) \right) $ is an orthonormal
basis for $L^{2}\left( \mu \right) $. But it is
not, as the function
$f_{\perp }\left( x\right) =e^{i\left( x_{1}+x_{2}\right) }$ is in
$L^{2}\left( \mu \right) $,
but not in $H_{2}\left( P\left( L\right) \right) $, since
$\left( \begin{smallmatrix}1\\1\end{smallmatrix}\right) 
\perp P\left( L\right) $. That would
make $f_{\perp }$ constant on every line
of slope $45$ degrees
and therefore constant on $\mathbf{X}\left( B\right) $%
. But we
can see that it is not from
the formula for $\mu $ and its
support, the attractor $\mathbf{X}\left( B\right) $.
This attractor has Hausdorff
dimension $\frac{\ln 3}{\ln 6}$: see Figure \ref{smithereens}%
, and the formula for $\mathbf{X}\left( B\right) $,
\begin{equation}
\label{attractorX}
\mathbf{X}\left( B\right) =
\vphantom{\sum _{i=0}^{\infty }}
\left\{ \vphantom{\sum}\smash{\sum _{i=0}^{\infty }}
b_{i}\frac{1}{6^{i}}:b_{i}\in B\right\} .
\end{equation}
From (\ref{contractivityforC}), and Lemma \ref{LemNew4.3},
we conclude that
$Q_{1}$ has the form
$Q_{1}\left(t_{1},t_{2}\right) =F\left(t_{1}+t_{2}\right) $
for some entire analytic
function $F$ of one variable. If
$E\left( P\left( L\right) \right) $ were an orthonormal
\emph{basis} for $L^{2}\left( \mu \right) $, it
would follow similarly that
every $f\in L^{2}\left( \mu \right) $ had the
representation 
$f\left(x_{1},x_{2}\right) =g_{f}\left(x_{1}-x_{2}\right) $
for some function $g_{f}$ of one
variable, which we noted is
not the case: that would make
$f_{\perp }\left( x_{1},x_{2}\right) =e^{i\left( x_{1}+x_{2}\right) }$
constant
on $\mathbf{X}\left( B\right) $ which it is not,
as can be checked by inspection of (\ref{attractorX}).
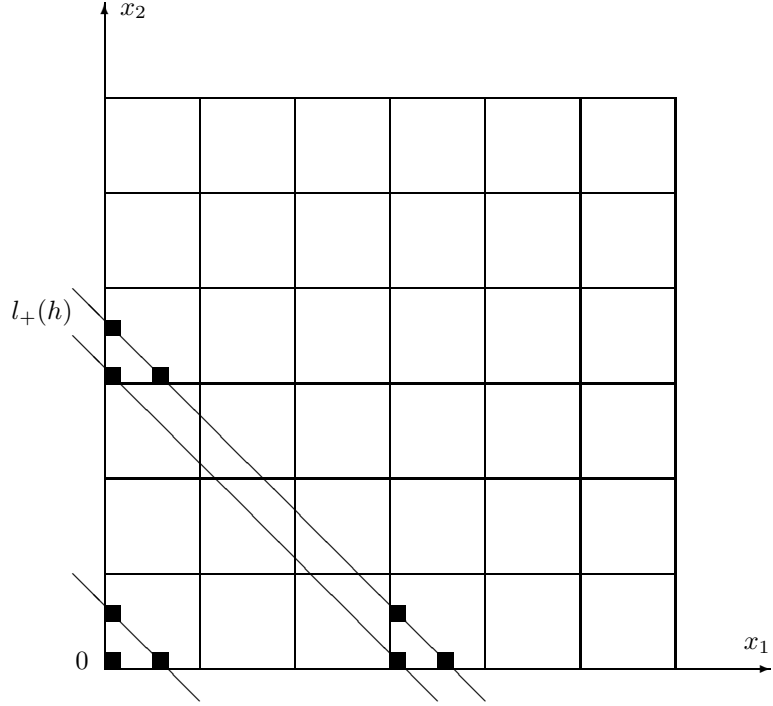
\begin{figure}[tbp]
\setlength{\unitlength}{6pt}
\begin{picture}(48,48)(-6,-6)
\multiput(6,0)(6,0){6}{\line(0,1){36}}
\put(0,0){\vector(0,1){42}}
\put(39,1){\makebox(3,3)[rb]{$x_{1}$}}
\multiput(0,6)(0,6){6}{\line(1,0){36}}
\put(0,0){\vector(1,0){42}}
\put(1,39){\makebox(3,3)[lt]{$x_{2}$}}
\put(-4,0){\makebox(3,3)[rb]{$0$}}
\put(6,-2){\line(-1,1){8}}
\put(21,-2){\line(-1,1){23}}
\put(24,-2){\line(-1,1){26}}
\put(-5,21){\makebox(3,3)[r]{$l_{+}(h)$}}
\linethickness{\unitlength}
\put(0,0.5){\line(1,0){1}}
\put(3,0.5){\line(1,0){1}}
\put(18,0.5){\line(1,0){1}}
\put(21,0.5){\line(1,0){1}}
\put(0,3.5){\line(1,0){1}}
\put(18,3.5){\line(1,0){1}}
\put(0,18.5){\line(1,0){1}}
\put(3,18.5){\line(1,0){1}}
\put(0,21.5){\line(1,0){1}}
\end{picture}
\caption{The Cut and Projection Fractal $\protect\mathbf{X}\left( B\right) $}
\label{smithereens}
\end{figure}

The argument shows that
functions $f\in H_{2}\left( P\left( L\right) \right) $ are
characterized by the following representation: let
\begin{equation*}
\Gamma _{l_{1}}:=
\left\{ \lambda \cdot l_{1}:\lambda \in P\left( L\right) \right\} .
\end{equation*}
Then $f\left( x_{1},x_{2}\right) =
\sum _{\xi \in \Gamma _{l_{1}}}
c\left( \xi \right) e_{\xi }\left( x_{1}-x_{2}\right) $,
where
$\sum _{\xi \in \Gamma _{l_{1}}}
\left| c\left( \xi \right) \right| ^{2}<\infty $.
\end{remark}

Inspection of (\ref{attractorX}) yields the
following representation of points
$x=\left(
\begin{smallmatrix}
x_{1} \\ x_{2}
\end{smallmatrix}
\right) $ in $\mathbf{X}\left( B\right) $:
\begin{equation*}
x=\frac{1}{2}\sum _{i=0}^{\infty }\eta _{i}
\begin{pmatrix}
\varepsilon _{i} \\
1-\varepsilon _{i}
\end{pmatrix}
\frac{1}{6^{i}},
\end{equation*}
where $\varepsilon _{i},\eta _{i}\in \left\{ 0,1\right\} $.
Let $\pi _{\pm}\left( x\right) :=x_{1}\pm x_{2}$,
and consider the lines $l_{\pm}\left( h\right) $
given by $x_{1}\pm x_{2}=h$, respectively,
for $h\in \mathbb{R}$, $\left| h\right| \leq 1$. Expanding
$h$ in base $6$, 
\begin{equation*}
h=\frac{1}{2}\sum _{i=0}^{\infty }h_{i}
\frac{1}{6^{i}},\quad 
h_{i}\in \left\{ 0,\pm 1,\pm 2\right\} ,
\end{equation*}
we have, for the
intersections
$\mathbf{X}_{\pm }\left( h\right) :=
\mathbf{X}\left( B\right) \cap l_{\pm }\left( h\right) $,
that
$\mathbf{X}_{+}\left( h\right) \neq \varnothing $
if and only if
$h_{i}\in \left\{ 0,1\right\} $. In that
case
\begin{align*}
\mathbf{X}_{+}\left( h\right) 
&=\left\{ x\in \mathbf{X}\left( B\right) :
\pi _{+}\left( x\right) =h\right\}   \\
&=
\vphantom{\sum _{i=0}^{\infty }}
\left\{ \frac{1}{2}
\vphantom{\sum}\smash{\sum _{i=0}^{\infty }}
h_{i}
\begin{pmatrix}
\varepsilon _{i} \\
1-\varepsilon _{i}
\end{pmatrix}
\frac{1}{6^{i}}:
\varepsilon _{i}\in \left\{ 0,1\right\} \right\} .
\end{align*}
Hence both $\mathbf{H}:=\left\{ h\in \mathbb{R}:
\mathbf{X}_{+}\left( h\right) \neq \varnothing \right\} $,
and the $\mathbf{X}_{+}\left( h\right) $'s, are fractals
(embedded in one dimension). The
respective Hausdorff dimensions are
$D_{H}\left( \mathbf{H}\right) =\frac{\ln 2}{\ln 6}$, and
$D_{H}\left( \mathbf{X}_{+}\left( h\right) \right) \leq \frac{\ln 2}{\ln 6}$,
with
$D_{H}\left( \mathbf{X}_{+}\left( h\right) \right) =\frac{\ln 2}{\ln 6}$,
generically.

On the other hand, the other
line-intersections are at most
singletons, i.e.,
$\mathbf{X}_{-}\left( h\right) 
=\left\{ x\in \mathbf{X}\left( B\right) :
\pi _{-}\left( x\right) =h\right\} $
is empty, or a singleton.
We have $\mathbf{X}_{-}\left( h\right) \neq \varnothing $ if and only if the
digits $\left( h_{i}\right) $ are restricted to 
$h_{i}\in \left\{ 0,\pm 1\right\} $; now a set of
Hausdorff dimension $\frac{\ln 3}{\ln 6}$. The
singleton intesection is then
\begin{multline*}
\mathbf{X}_{-}\left( h\right) 
=\left\{ 
x\sim \left( \eta _{i},\varepsilon _{i}\right) 
\in \mathbf{X}\left( B\right) :
\eta _{i}=0\text{ if }h_{i}=0,
\vphantom{\text{ and }
\eta _{i}=1,\,\varepsilon _{i}=\frac{1+h_{i}}{2}
\text{ if }h_{i}\in \left\{ \pm 1\right\} }
\right.   \\
\left. 
\vphantom{x\sim \left( \eta _{i},\varepsilon _{i}\right) 
\in \mathbf{X}\left( B\right) :
\eta _{i}=0\text{ if }h_{i}=0,}
\text{ and }
\eta _{i}=1,\,\varepsilon _{i}=\frac{1+h_{i}}{2}
\text{ if }h_{i}\in \left\{ \pm 1\right\} \right\} . 
\end{multline*}
\begin{summary}
The Hardy subspace
$H_{2}\left( P\left( L\right) \right) $ in $L^{2}\left( \mu \right) $,
for the $\mathbf{X}\left( B\right) $
fractal in $\mathbb{R}^{2}$, is canonically
associated with the following
one-dimensional pair:
$\left( \pi _{-}\left( \mathbf{X}\left( B\right) \right) ,
\Gamma _{l_{1}}\right) $, where
$\pi _{-}\left( \mathbf{X}\left( B\right) \right) $
is the above mentioned fractal in $\mathbb{R}^{1}$
\textup{(}of the same Hausdorff
dimension $\frac{\ln 3}{\ln 6}$\textup{)}, and
$\Gamma _{l_{1}}$ is the
associated spectrum.
\end{summary}

Our next result shows that one can find a bound for the constant $\gamma $
in Theorem \ref{th:Cbound}, the bound depending only on the data
$(R,B,L)$ in
a very explicit way.
Recall,
$\left\| T\right\|
_{op}$
is the operator norm, and $\left\| T\right\| _{hs}:=\left(
\sum _{j,k=1}^{\nu }\left| t_{j,k}\right| ^{2}\right) ^{1/2}$
is the Hilbert-Schmidt norm for a $\nu \times \nu $ matrix $T$.

\begin{theorem}
\label{T:estimate}Let $(R,B,L)$ be a system in
$\mathbb{R}^{\nu }$ satisfying \textup{(\ref
{B=L})--(\ref{Hadamard})}, $0\in L$. Let $C$ be the operator given by
\textup{(\ref
{eq21})}, and let $\mathbf{Y}%
$ denote the convex hull of the attractor $\mathbf{X}%
_{\rho }$
given by \textup{(\ref{rho-recur-dom})}.
If $\left\| Q\right\| _{\mathbf{Y}%
}$ is given by \textup{(\ref
{eq:Y-sup-norm})} and 
\begin{equation}
\label{T:estimate:gamma}
\gamma :=\left( N-1\right) 
^{2}N^{-1}
\beta 
\left\| R^{-1}\right\| _{op}
\max_{l\in L}\left| l\right| _{2}+\left\| R^{-1}\right\| _{hs},
\end{equation}
where 
\begin{equation*}
\beta :=2\pi \func{diam}\left( B\right) \max_{
\substack{b,b^{\prime }\in B \\
l\in L}}
\left\|
\sin (2\pi (b-b^{\prime })(\,\cdot \, -l))\right\| _{\infty },
\end{equation*}
$\left\| R\right\| _{op}$ is the operator norm of $R$, and $\left\|
T\right\| _{hs}:=\left( \sum
_{j,k=1}^{\nu }\left| t_{j,k}\right|
^{2}\right) ^{1/2}$ for a $\nu \times \nu $ matrix $T$, then we have 
\begin{equation}
\label{T:estimate:boundgradsupnorm}
\left\| CQ\right\| _{\mathbf{Y}%
}\leq \gamma \left\| Q\right\| _{\mathbf{Y}%
},
\end{equation}
for any $C^{1}$-function $Q$ such that $Q(0)=0$.
\end{theorem}

\begin{proof}
Begin by observing that 
\begin{align*}
\left| \chi _{B}(t-l)\right| ^{2} &=N^{-2}\sum_{b,b^{\prime }}e^{i2\pi %
(b-b^{\prime })\cdot (t-l)} \\
&=N^{-2}\sum_{b,b^{\prime }}\cos (2\pi (b-b^{\prime })(t-l)) \\
\intertext{so}
\left| \nabla \left( \left| \chi _{B}(t-l)\right| ^{2}\right) \right|
_{2} &=N^{-2}\sum_{b,b^{\prime }}2\pi \left| b-b^{\prime }\right| _{2}\left|
\sin (2\pi (b-b^{\prime })(t-l))\right| .
\end{align*}
Since the $N$ terms with $b=b^{\prime }$ are $0$, it follows that
\begin{equation}
\left\| \left| \nabla \left( \left| \chi _{B}(\,\cdot \,-l)\right| ^{2}\right)
\right| _{2}\right\| _{\infty }\leq (1-N^{-1})\beta
\label{eq:sinSumEstimate}
\end{equation}
for each $l\in L$. 
Since $\rho _{l}\mathbf{X}%
_{\rho }\subset \mathbf{X}%
_{\rho }$ and $\rho _{l}$ is affine we
have $\rho _{l}\mathbf{Y}%
\subset \mathbf{Y}%
$. If $0\leq Q(t)\leq 1$, then clearly we have
the following estimate:
\begin{multline}
\label{eq:estimateC}
\left\| CQ\right\| _{\mathbf{Y}%
}=\left\| \left| \nabla CQ\right| _{2}\right\|
_{\infty }   \\
\begin{aligned}
{}&=
\vphantom{\sum_{l}}
\left\| \left| \nabla \vphantom{\sum}\smash{\sum
_{l}}\left| \chi _{B}(\,\cdot \,
-l)\right| ^{2}Q(\rho _{l}t)\right| _{2}\right\| _{\infty } \\
&=
\vphantom{\sum_{l}}
\left\| \left| \vphantom{\sum}\smash{\sum
_{l}}\nabla \left( \left| \chi _{B}(\,\cdot \,
-l)\right| ^{2}\right) Q\circ \rho _{l}+
\vphantom{\sum}\smash{\sum
_{l}}\left| \chi
_{B}(\,\cdot \, -l)\right| ^{2}\nabla (Q\circ \rho _{l})\right| _{2}\right\|
_{\infty }   \\
&\leq
\vphantom{\sum_{l}}
\left\| \left| \vphantom{\sum}\smash{\sum
_{l}}\nabla \left( \left| \chi _{B}(\,\cdot \,
-l)\right| ^{2}\right) Q\circ \rho _{l}
\right| _{2}\right\| _{\infty } 
+
\vphantom{\sum_{l}}
\left\| \left| 
\vphantom{\sum}\smash{\sum
_{l}}\left| \chi
_{B}(\,\cdot \, -l)\right| ^{2}\nabla (Q\circ \rho _{l})
\right| _{2}\right\| _{\infty }   \\
&\leq 
\underset{\text{first term}}{\underbrace{
\vphantom{\sum_{l}}
\left\| \vphantom{\sum}\smash{\sum
_{l}}\left| \nabla \left( \left| \chi
_{B}(\,\cdot \, -l)\right| ^{2}\right) \right| _{2}Q\circ \rho _{l}\right\|
_{\infty }
}}+
\underset{\text{second term}}{\underbrace{
\vphantom{\sum_{l}}
\left\| \vphantom{\sum}\smash{\sum
_{l}}\left| \chi _{B}(\,\cdot \, -l)\right|
^{2}\left| \nabla (Q\circ \rho _{l})\right| _{2}\right\| _{\infty }
}}.
\end{aligned}
\end{multline}
Let us begin by estimating the
\emph{first term} in this last sum. 
Using
Lemma \ref{LemNew7.5} we have
\begin{multline*}
\sum
_{l
}
\nabla 
\left( 
\left| \chi _{B}(t-l)\right|
^{2}\right) 
Q
\left( 
\rho _{l}(t)
\right)  \\
=-\sum
_{l\in L^{*}}
\nabla 
\left( 
\left| \chi _{B}(t-l)\right|
^{2}\right) \int_{0}^{1}\left( R^{*\,-1}l\right) \cdot 
\nabla Q
\left(
R^{*\,-1}
\left( 
t-s
l\right) \right) \,ds.
\end{multline*}
Since $(R^{*-1}(t-sl))_{0\leq s \leq 1}$ is the line connecting
$\rho_0(t)$ and $\rho_l(t)$ it follows from (\ref{eq:sinSumEstimate})
that the right hand
side is dominated by
\begin{multline}
\label{eq:estimate_1}
\left\| \vphantom{\sum}\smash{\sum_{l\neq 0}}
\left| \nabla \left( \left| \chi _{B}(\,\cdot \,
-l)\right| ^{2}\right) \right| _{2}Q\circ \rho _{l}\right\| _{\infty } \\
\leq \left( N-1\right) 
\left( 1-N^{-1}\right) \beta \left\| R^{-1}\right\| _{op}
\max_{l\in L}\left| l\right| _{2}
\left\|
Q
\right\| _{
\mathbf{Y}}.  
\end{multline}
The $\left( N-1\right) $ 
factor is from an application of Lemma \ref{LemNew7.5}.
In the estimates above
we use%
d
the
invariance $\rho _{l}
\mathbf{
Y
}\subset 
\mathbf{
Y
}$
and the convexity of $\mathbf{Y}$.

To bound the
\emph{second term} in (\ref{eq:estimateC}),
we begin by noticing that
Lemma \ref{Lem4.1} and (\ref{eq29}) lead to 
\begin{equation}
\vphantom{\sum_{l}}
\left\| \vphantom{\sum}\smash{\sum
_{l}}\left| \chi _{B}(\,\cdot \, -l)\right| ^{2}\left|
\nabla (Q\circ \rho _{l})\right| _{2}\right\| _{\infty }\leq \max_{l\in
L}\left\| \left| \nabla (Q\circ \rho _{l})\right| _{2}\right\| _{\infty }.
\label{eq:estimate_2}
\end{equation}
If $c_{j}$ denotes the \emph{j}th column in
$R^{*\,-1}$, then 
\begin{equation*}
\frac{\partial \,}{\partial t_{j}}Q\circ \rho _{l}=\left( \left( \nabla
Q\right) \circ \rho _{l}\right) \cdot c_{j}
\end{equation*}
so that 
\begin{align*}
\left| \nabla (Q\circ \rho _{l})\right| _{2} &=
\vphantom{\sum_{j=1}^{\nu }}
\left(
\vphantom{\sum}\smash{\sum
_{j=1}^{\nu }}\left| \left( \left( \nabla Q\right) \circ \rho
_{l}\right) \cdot c_{j}\right| ^{2}\right) ^{1/2} \\
&\leq \left| \left( \left( \nabla Q\right) \circ \rho _{l}\right) \right|
_{2}
\vphantom{\sum_{j=1}^{\nu }}
\left(
\vphantom{\sum}\smash{\sum
_{j=1}^{\nu }}\left| c_{j}\right| _{2}^{2}\right) ^{1/2}
\\
&\leq \left\| Q\right\| _{\mathbf{Y}%
}\left\| R^{-1}\right\| _{hs},
\end{align*}
and it follows that 
\begin{equation*}
\vphantom{\sum_{l}}
\left\| \vphantom{\sum}\smash{\sum
_{l}}\left| \chi _{B}(\,\cdot \, -l)\right| ^{2}\left|
\nabla (Q\circ \rho _{l})\right| _{2}\right\| _{\infty }\leq \left\|
R^{-1}\right\| _{hs}\left\| Q\right\| _{\mathbf{Y}%
},
\end{equation*}
hence the proof is complete.
\end{proof}

\begin{corollary}
Let $(R,B,L)$ satisfy
\textup{(\ref{B=L})--(\ref{Hadamard})},
and, for $r\in \mathbb{N}$,
let 
\begin{equation*}
\mathcal{L}_{r}:=
\vphantom{\sum_{k=0}^{n}}
\left\{ \vphantom{\sum}\smash{\sum_{k=0}^{n}}
\left( rR^{*}\right) ^{k}l_{k}:n\in 
\mathbb{N},\;l_{k}\in L\right\}
\end{equation*}
and let $\mu _{r}$ be the probability measure solving 
\begin{equation*}
\mu _{r}=N^{-1}\sum_{b\in B}\mu _{r}\circ \sigma _{r,b}^{-1}
\end{equation*}
where $\sigma _{r,b}(x):=\left( rR\right) ^{-1}x+b$. 
Then 
$\{e_{\lambda }:\lambda \in \mathcal{L}_{r}\}$ is an
orthonormal basis for $L^{2}(\mu _{r})$ for $r$ sufficiently large.
\end{corollary}

\begin{proof}
This is an immediate consequence of 
Lemma \ref{Lem3.3}
and Theorem \ref
{T:estimate} 
since $\left\| rR\right\| _{op}^{-1}\rightarrow 0$ and $\left\|
\left( rR\right) ^{-1}\right\| _{hs}\rightarrow 0$ as $r\rightarrow \infty $.
\end{proof}

In one dimension, with $R$ an integer, and $N=2$,
the general situation is
summarized in the following strengthening
of Theorem \ref{ThmNew6.something} above.

\begin{corollary}
\label{CorNew9.7}
Suppose $\nu =1$, $N=2$, $B=\{0,b\}$, with
$b\in \mathbb{R}\diagdown \{0\}$, $R$
is an integer with $\left| R\right| \geq 2$, and $\mu $ is given by 
\textup{(\ref{determinemu})}.
If $R$ is odd, then $L^{2}(\mu )$ does not have a basis of
exponentials for any $b\in \mathbb{R}\diagdown \{0\}$.
If $R$ is even, and $%
\left| R\right| \geq 4$, then $L^{2}(\mu )$ has a basis of exponentials
for
all $b\in \mathbb{R}\diagdown \{0\}$.
If $R\geq 4$ and even, then $L$ may be
chosen in such a way that $L^{2}(\mu )$ is identified as the boundary
values of analytic functions in the unit disc.
\end{corollary}

\begin{proof}
The argument for odd $R$ is in the proof of Theorem \ref{ThmNew6.1}.

Suppose $R$ is even, and
let $L:=\{0,\frac{1}{|2b|}\}$.
Then (\ref{B=L})-(\ref{Hadamard})
are
satisfied, and 
we can easily find
the convex hull, $\mathbf{Y}$,
of the attractor $\mathbf{X}_{\rho }$. 
(See Lemma \ref{LemNew7.5}(\ref{LemNew7.5(3)}).)
In fact,
if $R>0$,
then $\mathbf{Y}$ is the closed interval with endpoints
$\frac{-1}{(|2b|(R-1))}$ and $0$;
while if $R<0$, then
$\mathbf{Y}$ is the closed inteval with endpoints $%
\frac{-1}{(|2b|(R^{2}-1))}$ and
$\frac{-1}{(|2bR|(R^{2}-1))}$. Using the notation from
Theorem \ref{T:estimate},
we have $\left\| R\right\| _{op}=\left| R\right|
$, 
$\left\| R^{-1}\right\| _{hs}=\left| R^{-1}\right| $, and 
\begin{equation*}
\beta =2\pi |b|\left\| \sin 2\pi b\,\cdot \,\right\| _{\infty }.
\end{equation*}
From the characterizations of $\mathbf{Y}$,
it follows 
that 
\begin{align*}
\mkern36mu \beta &=2\pi |b|\left| \sin \frac{\pi }{R-1}\right| , 
& &\text{if }R>0,\mkern36mu  \\
\intertext{and}
\beta &=2\pi |b|\sin \frac{\pi }{R^{2}-1}, & &\text{if }R<0. \\
\intertext{Using $N=2$, and $\max _{l\in L}b\left| l\right| _{2}=
\frac{1}{\left| 2b\right| }$, 
we conclude that }
\gamma &=\frac{\pi }{2R}\left| \sin \frac{\pi }{R-1}\right| +\frac{1}{R},
& &\text{if }R>0, \\
\intertext{and}
\gamma &=\frac{\pi }{2\left| R\right| }\left| \sin \frac{\pi }{R^{2}-1}
\right| +\frac{1}{\left| R\right| }, & &\text{if }R<0.
\end{align*}
Thus $\gamma <1$ if $\left| R\right| \geq 4$
, but not if $%
R=-2,2
$. 
\end{proof}

If $R=2$, then
$\{e_{\lambda }:\lambda \in \mathbb{Z}\}$ is an orthonormal
basis
for $L^{2}(\mu )$, and clearly
$\{e_{\lambda }:\lambda \in \mathcal{L}\}$
cannot be an orthonormal basis for $L^{2}(\mu )$ for any choice of $%
L=\{0,l\} $ with $l\neq 0$; in fact, depending on the sign of $l$, the set
$%
\mathcal{L}$ will consist solely of non-negative real numbers, or solely of
non-positive real numbers. We do not know what happens for general
$b\neq 0$.

If $R=-2$, then 
we do not 
know
if $L^{2}(\mu )$ has a basis of exponentials
for some $b\in \mathbb{R}\diagdown \{0\}$ .

\begin{acknowledgements}
The authors are very grateful to Brian Treadway for an expert typesetting
job, and for programming the Mathematica-generated graphics.
Helpful comments from R.S. Strichartz on a first version
are also gratefully acknowledged.
\end{acknowledgements}

\bibliographystyle{bftalpha}
\bibliography{jorgen}

\end{document}